\documentclass[fullpage,12pt]{amsart}
\usepackage{amssymb}
\usepackage{comment}
\usepackage{cite}
\usepackage{tikz}
\usepackage[aligntableaux=center]{ytableau}
\newcommand{\id}{\mbox{\rm id}}
\newcommand{\Ind}{\mbox{\rm Ind}}
\newcommand{\Res}{\mbox{\rm Res}}
\newcommand{\End}{\mbox{\rm End}}
\newcommand{\cont}{\mbox{\rm cont}}
\newcommand{\ch}{\mbox{\rm ch}}
\newcommand{\Ch}{\mbox{\rm Ch}}
\newcommand{\minich}{\mbox{\rm \scriptsize ch}}
\newcommand{\coch}{\mbox{\rm coch}}
\newcommand{\Stab}{\mbox{\rm Stab}}
\newcommand{\shape}{\mbox{\rm shape}}
\newcommand{\smax}{\mbox{\rm \scriptsize max}}
\newcommand{\smin}{\mbox{\rm \scriptsize min}}

\newcommand{\triv}{\mbox{\rm \scriptsize triv}}
\newcommand{\gr}{\mbox{\rm gr}}

\newtheorem{thm}{Theorem}[section]
\newtheorem{lemma}{Lemma}[section]
\newtheorem{prop}{Proposition}[section]
\newtheorem{cor}{Corollary}[section]

\newtheorem{defn}{Definition}[section]
\newtheorem{rem}{Remark}[section]

\title{From Young's Lattice to Coinvariants}
\author{Eugene Stern}
\address{Bloomberg L.P., 731 Lexington Ave., New York, NY 10022}
\email{eugene.stern@gmail.com}

\begin{document}
\ytableausetup{boxsize=1.1em}

\begin{abstract}
We extend Vershik and Okounkov's \cite{bib_VO} inductive spectral approach from irreducible representations of the symmetric group to the left and right regular representation. By following induced representations along paths in Young’s lattice, we find a rigid orthonormal weight basis for $\mathbb{C}[S_n]$ indexed by pairs of standard tableaux. Our main finding is that this weight basis already carries an implicit grading, given by the charge statistic on the tableau that records the induction path. More explicitly, after realizing $\mathbb{C}[S_n]$ as an $S_n$-bimodule inside the polynomial ring $\mathbb{C}[z_1, \ldots , z_n]$, we find that each weight basis vector has a natural minimal degree, which corresponds exactly to the charge of the induction tableau. We use this to define a degree-preserving isomorphism, which we call the {\em charge map}, from $\mathbb{C}[S_n]$ to the ring of coinvariants, showing that the usual graded view of the regular representation of $S_n$ can be derived from the branching alone, without appealing to geometric constructions. This exposes the structure behind the results of Ariki, Terasoma, and Yamada \cite{bib_ATY} on higher Specht polynomials. The proof that the charge map is an isomorphism is based on an algebraic connection between the charge and the action of adjacent transpositions on weight vectors in the seminormal representation of $S_n$.
\end{abstract}

\maketitle
\section{Introduction}
This work is inspired by Vershik and Okounkov’s paper \cite{bib_VO}, which made for wonderful reading during an isolated pandemic winter. The main idea of \cite{bib_VO} is that virtually the entire classical representation theory of the symmetric group, including the combinatorics of Young tableaux, can be derived from the branching rules for restriction and induction of representations between $S_n$ and $S_{n+1}$ as encoded in {\em Young’s lattice}. We extend this philosophy into what is typically considered the geometric side of the representation theory of $S_n$. Our main finding is that the {\em grading} commonly associated with the regular representation, which is typically obtained through geometric constructions and described combinatorially by Kostka polynomials via the charge statistic, can also be derived using only induction paths in Young’s lattice.

To trace the chain of ideas from \cite{bib_VO} to this paper in a bit more detail, Vershik and Okounkov derive the seminormal form of the irreducible representations of $S_n$ essentially from first principles, based on two main ideas:
\begin{enumerate}
\item Representations of $S_n$ can be treated from an inductive spectral point of view, in terms of weight vectors, with the Jucys-Murphy elements playing a role similar to that of a Cartan subalgebra in a Lie algebra.
\item The eigenvalues of Jucys-Murphy elements acting on weight vectors record how those vectors restrict to representations of nested subgroups of $S_n$, which can also be encoded via paths down Young’s lattice, or, equivalently, by standard Young tableaux.
\end{enumerate}
Our starting point was to ask how this picture changes when we replace restriction (from $S_n$ to $S_{n-1}$, $S_{n-2}$ and so on) with induction. This expands the natural setting from an irreducible $S_n$-module to the regular representation $\mathbb{C}[S_n]$. There, we have these counterparts to the points above:
\begin{enumerate}
\setcounter{enumi}{2}
\item The regular representation has a unique weight basis decomposition when we consider the left and right actions of $S_n$ (specifically, of the Jucys-Murphy elements inside $S_n$) simultaneously. In this decomposition, eigenspaces with a common weight for the right action are irreducible representations for the left action, and vice versa.
\item Pairs of standard Young tableaux (of the same shape) correspond to paths all the way up (from $S_1$ to $S_n$) and then all the way back down Young’s lattice. These label weight vectors in $\mathbb{C}[S_n]$, with each pair of tableaux representing the corresponding vector’s induction path for the right action and restriction path for the left action, or vice versa. 
\end{enumerate}
Thus, the left and right regular representation decomposes as
\begin{equation} \label{left_right_tableaux}
\mathbb{C}[S_n] = \bigoplus_{T, T'} \langle v_{T, T'} \rangle,
\end{equation}
where the sum is over pairs $T$ and $T'$ of standard tableaux having $n$ boxes and the same shape, and $v_{T, T'}$ is the corresponding left and right weight vector. Fixing $T$ in (\ref{left_right_tableaux}), we get an irreducible right $S_n$-module $V_{T, \bullet} = \bigoplus_{T'} \langle v_{T, T'} \rangle$ of isomorphism type $\lambda = \shape(T)$ (dimension $=$ number of standard tableaux of shape $\lambda$). Similarly, fixing $T'$, we get an irreducible left $S_n$-module $V_{\bullet, T'} = \bigoplus_{T} \langle v_{T, T'} \rangle$. This gives a canonical decomposition of $\mathbb{C}[S_n]$ into irreducible left and right $S_n$-modules.

The left regular representation of $S_n$ has a graded version, the space (or ring) of {\em coinvariants}, that typically arises via the geometry of configuration spaces or flags. Algebraically, this can be defined as the quotient
\begin{equation}
R[z_1, \ldots , z_n]_{S_n} = R[z_1, \ldots , z_n] / \langle e_1 , \ldots , e_n \rangle,
\end{equation}
where for us, $R = \mathbb{C}$ (though we could use any commutative ring $R$), $\langle e_1 , \ldots , e_n \rangle$ represents the ideal generated by positive degree symmetric polynomials, and the $n$ elementary symmetric polynomials $e_1, \ldots , e_n$ form one minimal generating set. This is a left graded $S_n$-module under the standard action of permuting the variables $z_i$, and the homogeneous components are submodules.

By a result typically credited to Stanley \cite{bib_S} and Lusztig, $R[z_1, \ldots , z_n]_{S_n}$ is isomorphic as an $S_n$-module to $\mathbb{C}[S_n]$. Further, the grading can be described as follows: the multiplicity of $V_{\lambda}$, the irreducible representation corresponding to shape $\lambda$, in the degree $d$ component, is equal to the number of standard tableaux of shape $\lambda$ having {\em major index} equal to $d$. This number, in turn, also counts standard tableaux of shape $\lambda$ having {\em (co)charge}\footnote{Charge and cocharge are complementary statistics that always add up to $\binom{n}{2}$, making it a matter of convention which to call which. In our setting, it is simplest to align charge with the major index, so that the two have the same generating function over standard tableaux, although they need not agree on each individual one.} $d$, the latter being a statistic that keeps track of the relative position of adjacent pairs of numbers inside a standard tableau.

The surprise, and the main result of this paper, is that the branching rules for $S_n$, as expressed through the decomposition (\ref{left_right_tableaux}), already determine this grading. More precisely, we show that there is a natural mapping $\mathbb{C}[S_n] \to R[z_1, \ldots , z_n]_{S_n}$, which we call the {\em charge map}, which is degree-preserving if we define the degree of $V_{\bullet, T'}$ to be the charge of $T'$. A more ad hoc variant of this mapping was given in the papers of Ariki, Terasoma, and Yamada \cite{bib_ATY,bib_TY} on {\em higher Specht polynomials}, so this result can be viewed as uncovering the intrinsic structure behind their work.

This paper is organized as follows. In Section \ref{sec_weight_basis}, we use induced representations to derive the weight basis for $\mathbb{C}[S_n]$ and work out its basic properties (orthogonality, duality between the left and right actions and induction and restriction, projection formulas, equivalence to matrix units). It turns out to be easiest to start by analyzing the diagonal of $\mathbb{C}[S_n]$, which is the span of weight vectors whose left and right weights are equal, and then extend across the entire regular representation using the left and right actions of $S_n$. Then, in Section \ref{sec_young_sym}, we relate the weight basis to the better known Specht-style basis defined in terms of Young symmetrizers, recovering the triangularity between these two bases.

In Section \ref{sec_grading_iso}, we use the weight basis decomposition to relate the regular representation to the space of coinvariants. We realize the regular representation $\mathbb{C}[S_n]$ in homogeneous polynomials in $n$ variables having multi-degree $(0, 1, \ldots , n-1)$, with monomials corresponding to (delta functions of) permutations, the left action of $S_n$ corresponding to permuting variables in terms of their indices as usual, and the right action corresponding to permuting variables in terms of their degrees. The key observation is that {\em charge is related to the action of adjacent transpositions in $S_n$ on the weight basis.} This is because when $(i \ i+1)$ acts on a weight vector, the sign of the leading term is determined by the relative positions of $i$ and $i+1$ in the corresponding tableau, and charge tracks precisely this. This connection determines the natural degree associated with each weight vector in the space of coinvariants.

We formalize this by the {\em charge map}, which maps $V_{\bullet, T'}$ into the homogeneous component of the coinvariant ring whose multi-degree is given by the {\em charge tableau} of $T'$. Most of the work of Section \ref{sec_grading_iso} is in proving that the weight basis of $\mathbb{C}[S_n]$ remains independent after this degree transformation and modulo symmetric functions. To do this, we adapt some of the tools in \cite{bib_ATY}, while working more directly in terms of the inductive structure on which our initial construction is based.

\section{Weight Basis for $\mathbb{C}[S_n]$} \label{sec_weight_basis}

\subsection{Inductive View of $S_n$}
Fix a positive integer $n$. The symmetric group $S_n$ of permutations of an $n$-element set (usually taken to be $\{1, \ldots , n \}$) contains a distinguished copy of $S_{n-1}$ as the subgroup of permutations that fix the last ($n$-th) element. This gives us a chain of inclusions
\begin{displaymath}
\{ \id \} = S_1 \subseteq S_2 \subseteq S_3 \subseteq \cdots \subseteq S_{n-1} \subseteq S_n.
\end{displaymath}

We can find generators for $S_n$ that respect this chain. To begin with, any permutation can be written in terms of {\em transpositions}, which are elements $(i \ j)$ of $S_n$ that interchange two elements $i$ and $j$ and leave the rest fixed. For example, a rearrangement $i_1 \ i_2 \ \ldots \ i_n$ of $1 \ 2 \ \ldots \ n$ can be realized by transposing $i_n$ with $n$, then transposing $i_{n-1}$ with $n-1$, and so on. (Note that by the time we transpose $i_k$ with $k$, it may have moved to a new location, $i_k'$, so the specific transposition we need to apply is not $(i_k \ k)$ but $(i_k' \ k)$.) Each transposition, in turn, can be written as a product of {\em adjacent transpositions} $(i \ i+1)$ by conjugating repeatedly: if $i < j$, we can first write
\begin{eqnarray*}(i \ i+2) &=& (i+1 \ i+2) \cdot (i \ i+1) \cdot (i+1 \ i+2) \\
(i \ i+3) &=& (i+2 \ i+3) \cdot (i \ i+2) \cdot (i+2 \ i+3), \\
\end{eqnarray*}
and so on, finally reaching
\begin{displaymath}
(i \ j) = (j-1 \ j) \cdot (i \ j-1) \cdot (j-1 \ j).
\end{displaymath}

The adjacent transpositions $\sigma_i = (i \ i+1)$ are often called {\em Coxeter generators} for the symmetric group. They respect the inclusion chain above, in the sense that $\sigma_1, \ldots , \sigma_{k-1}$ generate $S_k$ (or, in terms of algebras, the group algebra $\mathbb{C}[S_k]$), and adding in $\sigma_k$ extends this to generate $S_{k+1}$ (or $\mathbb{C}[S_{k+1}]$).

Non-adjacent $\sigma_i$ commute ($\sigma_i \sigma_j = \sigma_j \sigma_i$ if $j \neq i \pm 1$) but adjacent ones do not; instead we have $\sigma_i \sigma_{i+1} \sigma_i = \sigma_{i+1} \sigma_i \sigma_{i+1}$, which is known as the braid relation. However, we can partially express the extension of $S_{n-1}$ to $S_n$ in a more commutative way as follows. Write the following sums in the group algebra $\mathbb{C}[S_n]$:
\begin{displaymath}
X_k = (1 \ k) + (2 \ k) + \cdots + (k-1 \ k).
\end{displaymath}
By convention, we have $X_1 = 0$, and by definition, $X_2 = \sigma_1$. The $X_k$ are commonly known as the {\em Jucys-Murphy elements}, though, like everything else in the representation theory of the symmetric group, they have also been associated with Alfred Young. $X_k$ is a kind of symmetrization of $\sigma_{k-1}$ across $S_{k-1}$, which gives rise to commutativity properties: $X_k$ commutes with $S_{k-1}$ and the $X_k$ all commute with each other. This enables us to use the subalgebra generated by the $X_k$ in a manner analogous to the Cartan subalgebra in Lie theory.

We have the simple conjugation equation
\begin{displaymath}
\sigma_k \cdot X_k \cdot \sigma_k = X_{k+1} - \sigma_k.
\end{displaymath}
Multiplying on the left and on the right by $\sigma_k$, we obtain the equivalent commutation relations
\begin{equation}
\sigma_k X_{k+1} - X_k \sigma_k = X_{k+1} \sigma_k - \sigma_k X_k = 1. \label{jm-hecke-rel}
\end{equation}
The main idea is to analyze representations of $S_n$ in terms of a {\em weight basis}, i.e., a common eigenbasis of the $X_k$. Then the equations above allow us to derive the action of $\sigma_{n-1} \in S_n$ from the action of $X_{n-1}$ and $X_n$. This is what we meant by saying that the Jucys-Murphy elements partially express the extension of $S_{n-1}$ to $S_n$, even though adding $X_n$ to $\mathbb{C}[S_{n-1}]$ does not generate all of $\mathbb{C}[S_n]$.

\subsection{Degenerate Affine Hecke Algebra} \label{sec_dAHA}
The relations (\ref{jm-hecke-rel}) turn out to be powerful enough that it is useful to generalize them beyond the Jucys-Murphy elements. Add to $\sigma_1, \ldots , \sigma_{n-1}$ a set of formal generators $Y_1, \ldots , Y_n$ that commute with each other and satisfy
\begin{equation}
\sigma_k Y_{k+1} - Y_k \sigma_k = Y_{k+1} \sigma_k - \sigma_k Y_k = 1. \label{hecke-rel}
\end{equation}
The algebra generated over $\mathbb{C}$ by the $\sigma_i$ and $Y_j$ is called the {\em degenerate affine Hecke algebra}, which we will write as $H_n$. It has $\mathbb{C}[S_n]$ as a subalgebra, and also as a {\em quotient}, the latter through the mapping $Y_1 \mapsto 0$. The commutation relations (\ref{hecke-rel}) then imply that $Y_2, Y_3, \ldots , Y_n$ must map to the Jucys-Murphy elements $X_2, X_3, \ldots , X_n$. To see this, we can transform (\ref{hecke-rel}) back to the conjugation equation $Y_{k+1} = \sigma_k \cdot Y_k \cdot \sigma_k + \sigma_k$, and then $Y_1 = 0$ implies
\begin{eqnarray*}
Y_2 &=& \sigma_1 = (1 \ 2) = X_2, \\
Y_3 &=& \sigma_2 \sigma_1 \sigma_2 + \sigma_2 = (1 \ 3) + (2 \ 3) = X_3,
\end{eqnarray*}
and so on.

\subsection{Irreducible Representations and Young’s Lattice} \label{basic-rep-theory}
The inclusion chain of symmetric groups allows us to build up $\mathbb{C}[S_n]$, the regular representation of $S_n$, via a sequence of induced representations
\begin{eqnarray*}
\mathbb{C}[S_n] &=& \Ind_{S_{n-1}}^{S^n} \mathbb{C}[S_{n-1}] \\
&=& \Ind_{S_{n-1}}^{S_n} \left( \Ind_{S_{n-2}}^{S_{n-1}} \mathbb{C}[S_{n-2}] \right) \\
&=& \Ind_{S_{n-1}}^{S_n} \left( \Ind_{S_{n-2}}^{S_{n-1}} \left( \Ind_{S_{n-3}}^{S_{n-2}} \left( \cdots \left( \Ind_{S_1}^{S_2} \mathbb{C}[S_1]  \right) \cdots \right) \right) \right).
\end{eqnarray*}
The irreducible representations $V_{\lambda}$ of $S_n$ are indexed by partitions $\lambda$ of $n$. We will write partitions in the usual way: $\lambda = (\lambda_1, \ldots , \lambda_l)$, where $\lambda_i \geq \lambda_{i+1}$ and the $\lambda_i$ add up to the number being partitioned. We will also represent partitions by their {\em Young diagrams}, which are left-justified arrays of boxes where, if we count rows from top to bottom, the $i$-th row has $\lambda_i$ boxes. A {\em Young tableau} of shape $\lambda$ is a filling of the $n$ boxes of $\lambda$ by the numbers $1, \ldots , n$, each appearing once, and a tableau is {\em standard} if its entries increase reading along rows and columns, reading to the right and down.

We will describe two ways to associate an irreducible representation $V_{\lambda}$ to a partition $\lambda$. Section \ref{sec_young_sym} will go into the relationship between the two.

The first way is known as {\em Young’s natural representation} and can be expressed as follows. For any ring $R$, $S_n$ acts on the left on the ring $R[z_1, z_2, \ldots , z_n]$ of polynomials in $z_1, z_2, \ldots , z_n$ by permuting variables in terms of their indices: the transposition $(i \ j)$ interchanges $z_i$ and $z_j$. Given a (not necessarily standard) tableau $T$ of shape $\lambda$, define the {\em Specht polynomial}
\begin{displaymath}
f_T(z_1, \ldots , z_n) = \prod (z_j - z_i),
\end{displaymath}
where the product is over pairs $i < j$ such that $i$ and $j$ lie in the same column of $T$. For example, if $T = \ytableaushort{124,35}$, then $f_T(z_1, \ldots , z_5) = (z_3 - z_1) (z_5 - z_2)$.

We can break up $f_T$ into factors corresponding to each column:
\begin{equation} \label{eq-specht-poly}
f_T(z_1, \ldots , z_n) = \Delta_{\lambda_1'} \cdot \Delta_{\lambda_2'} \cdots \Delta_{\lambda_{l'}'}.
\end{equation}
where:
\begin{itemize} 
\item $\lambda_1', \lambda_2', \ldots , \lambda_{l'}'$ are the lengths of the columns of $\lambda$ (equivalently, the rows of the {\em conjugate partition} $\lambda^* = (\lambda_1', \lambda_2', \ldots , \lambda_{l'}')$, where we conjugate a partition by exchanging rows and columns in its Young diagram).
\item For a variable set $\{z_{i_1}, \ldots , z_{i_r}\} \subseteq \{z_1, \ldots , z_n\}$, the polynomial $\Delta_r (z_{i_1}, \ldots , z_{i_r})$ is the {\em Vandermonde determinant} in the variables $\{z_{i_1}, \ldots , z_{i_r}\}$:
\begin{displaymath}
\Delta_r (z_{i_1}, \ldots , z_{i_r}) = \prod_{1 \leq j < k \leq r} (z_{i_k} - z_{i_j}).
\end{displaymath}
\item $\Delta_{\lambda_j'}$ has as its arguments those variables $z_i$ whose indices $i$ lie in the $j$-th column of $T$.
\end{itemize}

The left action of $S_n$ on polynomials takes one Specht polynomial $f_T$ to another as the arguments of each $\Delta_{\lambda_j'}$ in the product (\ref{eq-specht-poly}) change, but the structure of the product, corresponding to the diagram $\lambda$, remains the same. The resulting $S_n$-module $V_{\lambda}$ spanned by the Specht polynomials is called a {\em Specht module}. It is known to be irreducible, with the $f_T$ that correspond to {\em standard} tableaux $T$ forming a basis, so that the number of standard tableaux of shape $\lambda$ is the dimension of $V_{\lambda}$.

The second way of associating irreducible $S_n$-modules with partitions is called {\em Young’s seminormal form} and has been described roughly as follows by Vershik and Okounkov in \cite{bib_VO}. Given an irreducible representation $V_{\lambda}$ of $S_n$, for each box $b$ in the Young diagram of $\lambda$, define the {\em content} of $b$, written $\cont(b)$, by
\begin{displaymath}
\cont(b) = \mbox{\rm column index of } i - \mbox{\rm row index of } i.
\end{displaymath}
For example, the first box in the first row has content $0$ because its row index and column index are both equal to $1$. The second box in the first row has content $1$, and the first box in the second row has content $-1$.

Then $V_{\lambda}$ has a weight basis $\{v_T\}$ of eigenvectors of $X_1, \ldots , X_n$, indexed by standard tableaux $T$ of shape $\lambda$, where the $X_i$ act on $v_T$ as follows:
\begin{equation}
(X_1, \ldots , X_n) \cdot v_T = (\cont(b(T, 1)), \ldots , \cont(b(T, n)) ) \cdot v_T, \label{eq-diagonal}
\end{equation}
where $b(T, i)$ is the box of $T$ containing $i$. (This recovers the dimension of $V_{\lambda}$ as the number of standard tableaux $T$ of shape $\lambda$.)

For example, since $1$ is always in the box in the first row and column, we have $\cont(b(T, 1)) = 0$ for every $T$, corresponding to $X_1 = 0$. We have $\cont(b(T, 2)) = 1$ or $-1$, depending on whether $2$ lies to the right of or below $1$ in $T$. Since $X_2 = (1 \ 2)$, the generator of $S_2$, we have $X^2 = 1$, which implies that $X_2$ must indeed act as $1$ (trivial representation of $S_2$) or $-1$ (sign representation of $S_2$) on any eigenvector.

More generally, the relations (\ref{jm-hecke-rel}) allow us to derive the action of the $\sigma_i$ on the weight basis from the diagonal action given by (\ref{eq-diagonal}). Following \cite{bib_VO}, consider the subalgebra $H_{i,i+1}$ of $H_n$ generated by $X_i$, $X_{i+1}$, and $\sigma_i$. Given a representation of $S_n$ or, more generally, of $H_n$ (just replace the $X_i$ by $Y_i$), consider any $H_{i,i+1}$-invariant subspace. Since $X_i$ and $X_{i+1}$ commute, they have a common eigenvector $v$ in this subspace. Let $X_i \cdot v = av$ and $X_{i+1} \cdot v = bv$. Then by (\ref{jm-hecke-rel}), we have
\begin{eqnarray*}
X_i \cdot \sigma_i v &=& \sigma_i X_{i+1} v - v = b \cdot \sigma_i v - v, \\
X_{i+1} \cdot \sigma_i v &=& \sigma_i X_i v + v = a \cdot \sigma_i v + v.
\end{eqnarray*}
Since $X_i$, $X_{i+1}$, and $\sigma_i$ all preserve the span of $v$ and $\sigma_i v$, this subspace is invariant under the action of $H_{i,i+1}$.

Now consider the vectors $v \pm \sigma_i v$. By construction, they are eigenvectors for $\sigma_i$, with eigenvalues $\pm 1$. Let us compute the action of $X_i$ and $X_{i+1}$:
\begin{eqnarray*}
X_1 \cdot (v \pm \sigma_i v) &=& (a \mp 1) v \pm b \cdot \sigma_i v, \\
X_2 \cdot (v \pm \sigma_i v) &=& (b \pm 1) v \pm a \cdot \sigma_i v.
\end{eqnarray*}
We have the following possibilities:
\begin{enumerate}
\item $b = a+1$, $v$ and $\sigma_i v$ are dependent. $v - \sigma_i v$ is an eigenvector of $X_i$ and $X_{i+1}$ with eigenvalues $a+1$ and $a$. Since $\sigma_i^2 = 1, \sigma_i v = \pm v$. If we had $\sigma_i v = -v$, then $v - \sigma_i v = 2v$, which forces $X_i \cdot (v - \sigma v)$ to equal $av$ and $(a+1)v$ simultaneously. Hence $\sigma_i v = v$. \label{case_row}

\item $b = a-1$, $v$ and $\sigma_i v$ are dependent. $v + \sigma_i v$ is an eigenvector of $X_i$ and $X_{i+1}$ with eigenvalues $a-1$ and $a$. Since $\sigma_i^2 = 1, \sigma_i v = \pm v$. If we had $\sigma_i v = v$, then $v + \sigma_i v = 2v$,which forces $X_i \cdot (v + \sigma v)$ to equal $av$ and $(a-1)v$ simultaneously. Hence $\sigma_i v = -v$. \label{case_column}

\item $b = a+1$, $v$ and $\sigma_i v$ are independent. Then $w = v - \sigma_i v$ is an eigenvector of $\sigma_i$ with eigenvalue $-1$, and we have $X_i \cdot w = (a+1) w$, $X_{i+1} \cdot w = aw$. This is equivalent to case (\ref{case_column}) above, with $a$ substituted for $a+1$ and $v$ for $w$.

\item $b = a-1$, $v$ and $\sigma_i v$ are independent. Then $w = v + \sigma_i v$ is an eigenvector of $\sigma_i$ with eigenvalue $1$, and we have $X_i \cdot w = (a-1) w$, $X_{i+1} \cdot w = aw$. This is equivalent to case (\ref{case_row}) above, with $a+1$ substituted for $a$ and $v$ for $w$. (The point here is that if $b = a \pm 1$, then we can always find an eigenvector of $X_i$ and $X_{i+1}$ on which $\sigma_i$ acts by $\pm 1$, even if it is not the one we initially picked.)

\item $b \neq a \pm 1$. In this case, the span of $v$ and $\sigma_i v$ is irreducible as an $H_{i,i+1}$-module. The other eigenvector of $X_i$ and $X_{i+1}$ it contains is $v + (a-b) \sigma_i v$, with eigenvalues $b$ and $a$. If we write $v = v_{(a, b)}$ to represent its spectrum $(a, b)$ under $X_i$, $X_{i+1}$, we should write $v + (a-b) \sigma_i v = v_{(b, a)}$. Then we can write the action of $\sigma_i$ as
\begin{displaymath}
\sigma_i \cdot v_{(a,b)} = \frac{1}{b-a} \left( v_{(a,b)} - v_{(b,a)} \right) .
\end{displaymath}
\end{enumerate}
The first case on this list corresponds to $i+1$ lying directly to the right of $i$ in $T$ ($\sigma_i$ restricts to the trivial representation), the second to $i+1$ lying directly below $i$ in $T$ ($\sigma_i$ restricts to the sign representation), and the last to $i$ and $i+1$ being in different rows and columns, so $\sigma_i$ can exchange them. More precisely, the case-by-case analysis above boils down to:
\begin{prop} \label{prop_local_action}
Let $T$ be a standard tableau of shape $\lambda$ with $\cont(b(T, i)) = a$, $\cont(b(T, i+1)) = b$. Let $v_T$ be the corresponding weight vector in $V_{\lambda}$.
\begin{enumerate}
\item If $i$ and $i+1$ are not immediately adjacent in $T$, then $\sigma_i \cdot T$ is also a standard tableau, with $\cont(b(\sigma_i \cdot T, i)) = b$, $\cont(b(\sigma_i \cdot T, i+1)) = a$. Write $v_{\sigma_i T}$ for the corresponding weight vector. Then we can scale $v_{\sigma_i T}$ so that
\begin{displaymath}
\sigma_i \cdot v_T = \frac{1}{b-a} \left( v_T - v_{\sigma_i T} \right).
\end{displaymath}
\item If $i$ and $i+1$ are immediately adjacent in $T$, then $b = a \pm 1$ and
\begin{displaymath}
\sigma_i \cdot v_T = \frac{1}{b-a} \cdot v_T.
\end{displaymath}
\end{enumerate}
\end{prop}

We highlight two aspects of Proposition \ref{prop_local_action} that underlie the rest of this paper. The first point is that it implies that the action of $S_n$ on the weight basis is {\em local}, in the sense that $\sigma_i$ acts only on $i$ and $i+1$ in each standard tableau $T$. This implies that the weight basis is aligned with restriction from $S_n$ to $S_{n-1}$. The second point is that the {\em sign} of the $v_T$ term in $\sigma_i \cdot v_T$ records whether $i+1$ lies northeast (plus sign) or southwest (minus sign) of $i$ in $T$. It will take a bit of time to realize the full value of the second point, but we can begin using the first now.

Let $V_{\lambda}$ be an irreducible representation of $S_n$, and let $T$ be a standard tableau of shape $\lambda$. Let us define
\begin{displaymath}
T^{(i)} = \mbox{\rm subtableau of $T$ containing $1,2, \ldots , i$}.
\end{displaymath}
Consider the restriction $\Res_{S_{n-1}}^{S_n} V_{\lambda}$. Define an {\em outer box} of the Young diagram $\lambda$ as one that has no box below or to the right, so that removing it leaves another Young diagram with $n-1$ boxes. I.e., outer boxes of $\lambda$ are in bijection with $n-1$-box Young diagrams $\lambda'$ with $\lambda' \subseteq \lambda$. (Here $\lambda'$ does {\em not} mean conjugate partition.) Equivalently, an outer box of $\lambda$ is one that can contain $n$ in a standard tableau of shape $\lambda$, and removing the $n$ from such a tableau leaves a standard tableau of shape $\lambda'$.

Given the weight decomposition
\begin{displaymath}
V_{\lambda} = \bigoplus_T \langle v_T \rangle,
\end{displaymath}
where the sum is over all standard tableaux $T$ of shape $\lambda$, let $\lambda' \subseteq \lambda$ have $n-1$ boxes, and write
\begin{displaymath}
V_{\lambda'} = \bigoplus_T \langle v_T \rangle,
\end{displaymath}
where this sum is over standard tableaux $T$ with $1, \ldots , n-1$ making up a subtableau of shape $\lambda'$, and $n$ written in the remaining outer box that we remove from $\lambda$ to obtain $\lambda'$. Then we have
\begin{equation}
V_{\lambda} = \bigoplus_{\lambda' \subset \lambda} V_{\lambda'} \label{eq-young-rule}
\end{equation}
as a direct sum of vector spaces, taken over all $\lambda'$ with $n-1$ boxes contained in $\lambda$.

The local nature of the action of the $\sigma_i$ implies that each $V_{\lambda'}$ is an $S_{n-1}$-submodule, as only $\sigma_{n-1}$ can move $n$ to another (outer) box. The weight description of the irreducible representations, applied to $S_{n-1}$, implies that the $V_{\lambda'}$ are irreducible and distinct. Hence we have

\begin{prop} \label{prop-young-rule}
Let $V_{\lambda}$ be an irreducible $S_n$-module. $\Res_{S_{n-1}}^{S_n} V_{\lambda}$ decomposes into irreducible $S_{n-1}$-modules according to (\ref{eq-young-rule}), and this decomposition is multiplicity-free.
\end{prop}

We can apply this proposition repeatedly to the weight vector $v_T \in V_{\lambda}$ corresponding to a standard tableau $T$ of shape $\lambda$. Restricting from $S_n$ to $S_{n-1}$, we find that $v_T$ lies in the unique irreducible submodule of $V_{\lambda}$ isomorphic to $V_{\lambda'}$, where $\lambda'$ is the shape of $T^{(n-1)}$ (i.e., of $T$ with the box containing $n$ removed). Restricting from $S_{n-1}$ to $S_{n-2}$, we find that $v_T$ lies in the unique irreducible submodule of $V_{\lambda'}$ isomorphic to $V_{\lambda''}$, where $\lambda''$ is the shape of $T^{(n-2)}$ (i.e., of $T$ with the boxes containing $n$ and $n-1$ removed). As we continue to restrict from $S_k$ to $S_{k-1}$ in this way, the isomorphism type of the $S_{k-1}$-module in which $v_T$ lies is unique at each step. This generates a well-defined restriction path terminating at an irreducible $S_1$-module, which is necessarily $1$-dimensional, hence nothing more than the span of $v_T$.

\begin{figure}[htbp]
    \centering
    \ytableausetup{boxsize=0.8em}
    \begin{tikzpicture}[
        scale=1.3,
        node distance=1.8cm,
        edge/.style={draw, thick, ->, shorten >=2pt, shorten <=2pt, color=gray!80}
    ]

        \node (1) at (0,-1.5) {$\ytableaushort{\ }$};

        \node (2) at (-1,-3) {$\ytableaushort{\ \ }$};
        \node (11) at (1,-3) {$\ytableaushort{\ , \ }$};

        \node (3) at (-2.5,-4.5) {$\ytableaushort{\ \ \ }$};
        \node (21) at (0,-4.5) {$\ytableaushort{\ \ , \ }$};
        \node (111) at (2.5,-4.5) {$\ytableaushort{\ , \ , \ }$};

        \node (4) at (-3.5,-6.5) {$\ytableaushort{\ \ \ \ }$};
        \node (31) at (-1.7,-6.5) {$\ytableaushort{\ \ \ , \ }$};
        \node (22) at (0,-6.5) {$\ytableaushort{\ \ , \ \ }$};
        \node (211) at (1.7,-6.5) {$\ytableaushort{\ \ , \ , \ }$};
        \node (1111) at (3.5,-6.5) {$\ytableaushort{\ , \ , \ , \ }$};

        \draw[edge] (1) -- (2);
        \draw[edge] (1) -- (11);

        \draw[edge] (2) -- (3);
        \draw[edge] (2) -- (21);
        \draw[edge] (11) -- (21);
        \draw[edge] (11) -- (111);

        \draw[edge] (3) -- (4);
        \draw[edge] (3) -- (31);
        \draw[edge] (21) -- (31);
        \draw[edge] (21) -- (22);
        \draw[edge] (21) -- (211);
        \draw[edge] (111) -- (211);
        \draw[edge] (111) -- (1111);

    \end{tikzpicture}
\end{figure}

\ytableausetup{boxsize=1.1em}

Pictorially, we can represent the restriction path as a walk down {\em Young’s lattice}, which depicts how we obtain bigger Young diagrams from smaller ones by adding boxes one at a time, or, equivalently, how we obtain smaller diagrams from bigger ones by taking away boxes one at a time. Paths in the lattice correspond to different ways of doing either. They also correspond to standard Young tableaux, as requiring the entries of a tableau to increase down rows and columns is exactly what guarantees that removing the box with the largest entry generates another legal diagram. In terms of the lattice, our discussion of the chain of restrictions shows that:
\begin{prop}
Weight vectors for $V_{\lambda}$ correspond to paths down Young’s lattice, in the sense that $v_T$ can be characterized as the unique (up to scalars) weight vector in $V_{\lambda}$ whose restriction path from $S_n$ to $S_1$, as defined by the isomorphism classes of the irreducible $S_k$-modules along the restriction path in which it sits, is specified by the tableau $T$.
\end{prop}
We will sometimes refer to $T$ as the {\em restriction tableau} of the weight vector $v_T$.

\subsection{Weight Basis for $\mathbb{C}[S_n]$}
Our first extension of the Vershik-Okounkov picture is to change direction in Young’s lattice. Let $\lambda$ be a partition of $n$, let $\lambda'$ be a partition of $n-1$, let $\chi_{\lambda}$ denote the character of $V_{\lambda}$, and let $\langle \cdot , \cdot \rangle$ denote the inner product of characters. Then Frobenius reciprocity states that
\begin{equation}
\langle \Res_{S_{n-1}}^{S_n} \chi_{\lambda}, \chi_{\lambda'} \rangle = \langle \chi_{\lambda}, \Ind_{S_{n-1}}^{S_n} \chi_{\lambda'}. \rangle \label{eq_frob_rec}
\end{equation}
By Proposition \ref{prop-young-rule}, the left hand side is $1$ if $\lambda' \subset \lambda$, otherwise $0$, hence the right hand side must be so as well. In other words, Young’s lattice determines the structure of induction from $S_{n-1}$ to $S_n$ as well as restriction from $S_n$ to $S_{n-1}$: for any irreducible $S_{n-1}$-module $V_{\lambda'}$, the decomposition of $\Ind_{S_{n-1}}^{S_n} V_{\lambda'}$ into irreducibles is multiplicity-free and the $S_n$-modules that appear correspond to diagrams $\lambda$ obtained from $\lambda'$ by adding a box.

This means that we can interpret a standard tableau $T$ of shape $\lambda$ as a path up Young’s lattice, as follows. Let $V_{T^{(1)}} = \mathbb{C}[S_1]$ be the $1$-dimensional trivial representation of $S_1$. Consider the induced module $\Ind_{S_1}^{S_2} V_{T^{(1)}} = \mathbb{C}[S_2]$, and define $V_{T^{(2)}} \subseteq \mathbb{C}[S_2]$ to be the unique irreducible component of $\mathbb{C}[S_2]$ whose isomorphism class is given by the shape of $T^{(2)}$.

Since $V_{T^{(2)}} \subseteq \mathbb{C}[S_2]$, we can consider $\Ind_{S_2}^{S_3} V_{T^{(2)}} \subseteq \Ind_{S_2}^{S_3} \mathbb{C}[S_2] = \mathbb{C}[S_3]$. We define $V_{T^{(3)}} \subseteq \mathbb{C}[S_3]$ to be the unique irreducible component of $\Ind_{S_2}^{S_3} V_{T^{(2)}}$ whose isomorphism class is given by the shape of $T^{(3)}$. Continuing in this way, we obtain an irreducible submodule $V_T \subseteq \mathbb{C}[S_n]$ of isomorphism type $\lambda$, the shape of $T$.

\begin{prop} \label{direct_sum_gz}
The regular representation $\mathbb{C}[S_n]$ decomposes into irreducible $S_n$-modules as
\begin{displaymath}
\mathbb{C}[S_n] = \bigoplus_{T : \ |T|=n} V_T,
\end{displaymath}
where the sum is over all standard tableaux $T$ with $n$ boxes (of all shapes).
\end{prop}

\begin{proof}
This is clear by induction on $n$, as we have
\begin{eqnarray}
\mathbb{C}[S_n] &=& \Ind_{S_{n-1}}^{S_n} \mathbb{C}[S_{n-1}] \\
&=& \Ind_{S_{n-1}}^{S_n} \left( \bigoplus_{T' : \ |T'|=n-1} V_{T'} \right) \label{eq_ind_hyp} \\
&=& \bigoplus_{T' : \ |T'|=n-1} \left( \Ind_{S_{n-1}}^{S_n} V_{T'} \right) \\
&=& \bigoplus_{T' : \ |T'|=n-1} \left( \bigoplus_{T \supset T', \ |T|=n} V_T \right) \label{eq_box_add} \\
&=& \bigoplus_{T : \ |T|=n} V_T.
\end{eqnarray}
Here (\ref{eq_ind_hyp}) is the inductive hypothesis, (\ref{eq_box_add}) is a restatement of (\ref{eq_frob_rec}) together with the definition of $V_T$, and the last step follows by collecting together all standard tableaux $T$ with $n$ boxes and forgetting what $T^{(n-1)}$ happens to be.
\end{proof}

Borrowing from \cite{bib_VO}, we will refer to $V_T$ as the {\em Gelfand-Tsetlin module} corresponding to the standard tableau $T$. We will also call $T$ the {\em induction tableau} of $V_T$. The direct sum decomposition makes plain the well known fact that the multiplicity of $V_{\lambda}$ in the regular representation is equal to its dimension, the number of standard tableaux of shape $\lambda$.

Each $V_T$ has a weight basis, indexed by the standard tableaux $T'$ of the same shape as $T$. Let us write this weight basis as $\{ v_{T', T} \}$; here $T'$ plays the role of the restriction tableau and $T$ is a fixed induction tableau. Since $V_T$ is spanned by $\{ v_{T', T} \}$, with $T'$ varying and $T$ fixed, we could also write $V_T = V_{\bullet, T}$, implicitly thinking of each (left) Gelfand-Tsetlin module as the column of a matrix. This is not needed right now, but will be a useful way to keep track of things later on.

Collecting the weight bases over all $V_T$, or all induction tableaux $T$, we can conclude by Proposition \ref{direct_sum_gz} that the set
\begin{displaymath}
\{ v_{T', T} \}_{\mbox{\rm \scriptsize $T$, $T'$ standard of the same shape}}
\end{displaymath}
is a weight basis for $\mathbb{C}[S_n]$. Geometrically, the pair $(T, T')$ represents a circuit up and down Young’s lattice, from $1$ to $n$ and back again.

\subsection{Orthogonality}
The regular representation $\mathbb{C}[G]$ of any finite group $G$ has a standard $G$-invariant inner product $\langle \cdot , \cdot \rangle$ in which the individual group elements $g \in G$ form an orthonormal basis: $\langle g, g' \rangle = \delta_{g, g'}$. We have
\begin{prop} \label{prop-GZ-orthogonality}
The weight basis for $\mathbb{C}[S_n]$ is orthogonal under the standard $S_n$-invariant inner product.
\end{prop}

\begin{proof}
We show that
\begin{enumerate}
\item The weight basis inside any irreducible $S_n$-module is orthogonal under any $S_n$-invariant inner product, and
\item Two different irreducible Gelfand-Tsetlin modules $V_T$ and $V_{T'}$ inside $\mathbb{C}[S_n]$ are orthogonal.
\end{enumerate}

For the first part, let $V$ be an irreducible $S_n$ module. The key point is that restricting $V$ to $S_{n-1}$ is multiplicity-free. So let $W \subseteq V$ be an irreducible $S_{n-1}$-submodule. The orthogonal complement of $W$ is also invariant under $S_{n-1}$, because for $\sigma \in S_{n-1}$, $w \in W$, $w' \in W^{\perp}$, we have
\begin{displaymath}
\langle w, \sigma \cdot w' \rangle = \langle \sigma^{-1} \cdot w , \sigma^{-1} \sigma \cdot w' \rangle = \langle \sigma^{-1} \cdot w , w' \rangle = 0.
\end{displaymath}
Replacing $V$ by $W^{\perp}$, picking an irreducible $S_{n-1}$-submodule inside $W^{\perp}$, taking the orthogonal complement of that, and continuing, we finally obtain a decomposition of $V$ into mutually orthogonal irreducible $S_{n-1}$-modules. Since $\Res_{S_{n-1}}^{S_n} V$ is multiplicity-free, this matches the usual decomposition (\ref{eq-young-rule}). So (\ref{eq-young-rule}) is an orthogonal decomposition. Then the decomposition of each irreducible $S_{n-1}$-module into $S_{n-2}$-irreducibles is mutually orthogonal, and so on. Continue the process until it terminates at the weight vectors $\{v_T\}$, which must therefore be mutually orthogonal.

For the second part, we work by induction, in both senses of the word, on $k < n$. Given $k < n$ and an irreducible $S_k$-module $V \subseteq \mathbb{C}[S_k] \subseteq \mathbb{C}[S_n]$, we know that $\Ind_{S_k}^{S_{k+1}} V$ is also multiplicity-free. The argument from the first part (taking $\sigma \in S_{k+1}$ and aligning orthogonal complements with irreducibles) shows that the irreducible components of $\Ind_{S_k}^{S_{k+1}} V$ are mutually orthogonal. To conclude, we also need to show that induction preserves orthogonality, i.e., if $V$ and $V'$ are mutually orthogonal $S_k$-modules inside $\mathbb{C}[S_k]$, then $\Ind_{S_k}^{S_{k+1}} V$ and $\Ind_{S_k}^{S_{k+1}} V'$ are mutually orthogonal $S_{k+1}$-modules inside $\mathbb{C}[S_{k+1}]$.

Let $V$ and $V'$ be as above, and let $v \in \Ind_{S_k}^{S_{k+1}} V$, $v' \in \Ind_{S_k}^{S_{k+1}} V'$. Since we know that $(1 \ k+1), (2 \ k+1), \ldots , (k \ k+1)$ are a full set of left coset representatives for $S_k$ in $S_{k+1}$, by the definition of induced representation, we can write
\begin{eqnarray*}
v &=& (1 \ k+1) \cdot v_1 + (2 \ k+1) \cdot v_2 + \cdots + (k \ k+1) \cdot v_{k+1}, \\
v' &=& (1 \ k+1) \cdot v_1' + (2 \ k+1) \cdot v_2' + \cdots + (k \ k+1) \cdot v_{k+1}'
\end{eqnarray*}
for $v_i \in V$, $v_j' \in V'$, both inside $\mathbb{C}[S_k]$. By bilinearity of $\langle \cdot , \cdot \rangle$, it suffices to show
\begin{equation}
\langle (i \ k+1) \cdot v_i, (j \ k+1) \cdot v_j' \rangle = 0
\end{equation}
for all $i, j \leq k$.
If $i = j$, then by $S_n$-invariance of $\langle \cdot , \cdot \rangle$, we can write
\begin{displaymath}
\langle (i \ k+1) \cdot v_i, (i \ k+1) \cdot v_i' \rangle = \langle v_i, v_i' \rangle = 0,
\end{displaymath}
since $v_i \in V$, $v_i' \in V'$. If $i \neq j$, then expand $v_i$ and $v_j'$ as linear combinations of individual group elements in $\mathbb{C}[S_k]$. These group elements all act on $1, \ldots , k$ only and preserve $k+1$. Then $(i \ k+1) \cdot v_i$ is a linear combination of distinct elements of $S_{k+1}$ that take $k+1$ to $i$, and similarly $(j \ k+1) \cdot v_j$ is a linear combination of distinct elements of $S_{k+1}$ that take $k+1$ to $j'$. These two sets of are mutually distinct, hence mutually orthogonal by the definition of $\langle \cdot , \cdot \rangle$. Hence $\langle (i \ k+1) \cdot v_i, (j \ k+1) \cdot v_j' \rangle = 0$ in this case as well, concluding our proof.
\end{proof}

\subsection{Projection and Duality}
The regular representation of $G$ carries a right action as well as a left action. This means we can repeat our walks up and down Young’s lattice via induction and restriction, only with $S_n$ acting on the right, generating a ``right'' weight basis for $\mathbb{C}[S_n]$ parallel to the ``left'' one just constructed. Our goal now will be to show that this actually produces the {\em same} basis, just with the restriction and induction tableaux switched.

To begin to understand how interchanging restriction and induction in $\mathbb{C}[S_n]$ corresponds to interchanging the left and right actions of $S_n$, consider the following extension of the restriction rule, which can be found in \cite{bib_VO}:

\begin{prop}
Restriction eigen-lemma: Given a partition $\lambda$ of $n$, consider the decomposition of $\Res_{S_{n-1}}^{S_n} V_{\lambda}$ given by (\ref{eq-young-rule}):
\begin{displaymath}
\Res_{S_{n-1}}^{S_n} V_{\lambda} = \bigoplus_{\lambda' \subset \lambda, \ |\lambda'| = n-1} V_{\lambda'}.
\end{displaymath}
For $\lambda'$ with $n-1$ boxes contained in $V_{\lambda}$, we can characterize $V_{\lambda'}$ as the left eigenspace of $X_n$ corresponding to the eigenvalue $\cont(b')$, where $b'$ is the box removed from $\lambda$ to obtain $\lambda'$.
\end{prop}

\begin{proof}
We already know that for a fixed $n-1$-box diagram $\lambda' \subset \lambda$, $V_{\lambda'}$ is spanned by $v_T$ corresponding to those $T$ with $n$ written in the box $b$ removed from $\lambda$ to obtain $\lambda'$. Moreover, we know from (\ref{eq-diagonal}) that $X_n$ acts on all such $v_T$ by the eigenvalue $\cont(b')$. The only thing left is to confirm that the eigenvalues are all distinct, i.e., if $\lambda'$ and $\lambda''$ are two distinct $n-1$-box tableaux obtained from $\lambda$ by removing boxes $b'$ and $b''$, respectively, then $\cont(b') \neq \cont(b'')$. This is clear because level sets of the content function are diagonals running northwest to southeast. Thus, let $\cont(b') = \cont(b'')$, and assume without loss of generality that $b''$ lies in a row below $b'$. Then $b''$ also lies in a column to the right of $b'$. This forces the entry in $b''$ to always be higher than the entry in $b'$ in any standard tableau containing both boxes, meaning that the highest value $n$ could never be written in box $b'$ in a diagram of shape $\lambda$ with $n$ boxes. (More tersely: two boxes on the same NW-SE diagonal can never be outer boxes of the same Young diagram, and thus would not both be allowed to contain $n$.)
\end{proof}

A partial expression of the duality between left and right and between induction and restriction is the following analog of this lemma for induction from $S_{n-1}$ to $S_n$:

\begin{prop} \label{ind-lemma}
Induction eigen-lemma: Given a partition $\lambda'$ of $n-1$, let $V_{\lambda'}$ be any irreducible submodule of $\mathbb{C}[S_{n-1}]$ with isomorphism class $\lambda'$. Consider $\Ind_{S_{n-1}}^{S_n} V_{\lambda'}$ as a submodule of $\Ind_{S_{n-1}}^{S_n} \mathbb{C}[S_{n-1}] = \mathbb{C}[S_n] $. Then in the irreducible decomposition
\begin{displaymath}
\Ind_{S_{n-1}}^{S_n} V_{\lambda'} = \bigoplus_{\lambda \supset \lambda', \ |\lambda| = n} V_{\lambda},
\end{displaymath}
we can characterize $V_{\lambda}$ as the right eigenspace of $X_n$ corresponding to the eigenvalue $\cont(b)$, where $b$ is the box added to $\lambda'$ to obtain $\lambda$.
\end{prop}

To prove the induction eigen-lemma, which will take a little time, we will use orthogonal projectors onto the weight basis. We write these in an explicit inductive form which was already known to Murphy \cite{bib_M}, as follows:

Let $T' \subset T$ be standard tableaux, of shapes $\lambda'$ and $\lambda$, with $k-1$ and $k$ boxes respectively. Define an {\em outer corner} of $T'$ (more properly, of its shape $\lambda'$) to be a location where we could add a box to $\lambda'$ to form a legal Young diagram. For example, if $\lambda'$ is the $4$-box Young diagram $(3,1)$, we can extend it to form the $5$-box Young diagrams $(4,1)$, $(3,2)$, or $(3,1,1)$. Thus $(3,1)$ has $3$ outer corners, corresponding to the $3$ possible boxes we can add.

Now let $b$ be the box containing $k$ in $T$, so that removing $b$ from $T$ gives us $T'$. In terms of our language of ``out''-liers, $b$ is an outer box of $T$, and sits at an outer corner of $T'$. Let $b_1, \ldots , b_j$ represent all the other outer corners of $T'$, {\em not} including $b$. Let us write
\begin{displaymath}
p_{T', T}(X_k) = \prod_{i=1}^j \frac{X_k - \cont(b_i)}{\cont(b) - \cont(b_i)}.
\end{displaymath}
This definition really depends on the shapes $\lambda'$ and $\lambda$ rather than on $T'$ and $T$, but writing it in terms of tableaux makes for more intuitive notation in what follows.

Applying this definition to the example above, with $T'$ a tableau of shape $\lambda' = (3,1)$, we have outer corners in the first, second, and third row, of content $3$, $0$, and $-2$, respectively. If we let $T$ be a tableau of shape $\lambda = (4,1)$, corresponding to adding a box in the outer corner in the first row, we take the product over the two outer corners that remain, to obtain
\begin{displaymath}
p_{T', T}(X_5) = \frac{X_5 - 0}{3 - 0} \cdot \frac{X_5 + 2}{3 + 2} = \frac{1}{15} X_5 (X_5+2).
\end{displaymath}

Next, for $T$ a standard tableau with $n$ boxes, write
\begin{displaymath}
p_T(X_1, \ldots , X_n) = \prod_{k=2}^n p_{T^{(k-1)}, T^{(k)}}(X_k).
\end{displaymath}
For example, if $\lambda = (2,2)$ and $T = \ytableaushort{12,34}$, we have
\begin{eqnarray*}
p_{T^{(1)}, T^{(2)}}(X_2) &=& \frac12(X_2 + 1) , \\
p_{T^{(2)}, T^{(3)}}(X_3) &=& -\frac13(X_3 - 2), \\
p_{T^{(3)}, T^{(4)}}(X_4) &=& p_{T^{(3)}, T}(X_4) = \frac{X_4-2}{-2} \cdot \frac{X_4+2}{2} = -\frac14 (X_4-2)(X_4+2),
\end{eqnarray*}
and thus
\begin{displaymath}
p_T(X_1, X_2, X_3, X_4) = \frac{1}{24}(X_2+1)(X_3-2)(X_4-2)(X_4+2).
\end{displaymath}

We can interpret $p_{T', T}(X_k)$ and $p_T(X_1, \ldots , X_n)$ in several ways: as an element of $\mathbb{C}[S_n]$, as an operator acting on $\mathbb{C}[S_n]$ by multiplication (either on the left or on the right), or as an operator acting on the individual irreducible modules $V_{\lambda}$ (on the left, since we have been considering the modules that way). Considered as an element of $\mathbb{C}[S_n]$, $p_T(X_1, \ldots , X_n)$ is the image of the operator $p_T(X_1, \ldots , X_n)$ acting on the identity element $1 \in \mathbb{C}[S_n]$.

Considering $p_T$ as an operator, we have
\begin{prop} \label{left-proj}
If $T$ is a standard tableau of shape $\lambda$, then $p_T(X_1, \ldots , X_n)$ acts on the left on $V_{\lambda}$ as orthogonal projection onto the weight vector $v_T$.
\end{prop}

\begin{proof}
Since $X_k$ acts on $v_T$ as multiplication by the content of the box containing $k$, the numerator and denominator cancel in each fraction above, and $p_T$ acts on $v_T$ by 1. We show that $p_T \cdot v_{T'} = 0$ if $T \neq T'$. Let $k$ be the first number that appears in a different box in $T$ and in $T'$. Letting $b$ be the box added to $T^{(k-1)}$ to obtain $T^{(k)}$, and looking at the definition of $p_{T^{(k-1)}, T^{(k)}} (X_k) $, we have it that for some $b_i$, which is the box added to $T^{(k-1)}$ to obtain $T'^{(k)}$, we have $X_k = \cont(b_i)$, and hence $X_k - \cont(b_i) = 0$, $p_{T^{(k-1)}, T^{(k)}} (X_k) \cdot v_{T'} $ = 0, and $p_T (X_1, \ldots , X_n) \cdot v_{T'} = 0$. $p_T$ is an orthogonal projector since the $v_T$ form an orthogonal basis.
\end{proof}

Now consider the (left) action of $p_T(X_1, \ldots , X_n)$ on $\mathbb{C}[S_n]$, which we decompose as a direct sum of Gelfand-Tsetlin modules. If $T'$ is a tableau of different shape from $T$, $p_T$ kills $v_{T'}$ by the argument above, in which we looked at the first number that appears in different boxes in $T$ and $T'$. Hence $p_T$ kills the isotypic component of $\lambda'$ for any $\lambda'$ not equal to $\lambda$, the shape of $T$. In the isotypic component of $\lambda$, Proposition \ref{left-proj} tells us that for every $T'$ of shape $\lambda$, $p_T$ takes the Gelfand-Tsetlin module $V_{T'}$ to the weight vector $v_{T, T'}$. The span of all these vectors is the weight space in $\mathbb{C}[S_n]$ on which the $X_i$ act by the content vector $(\cont(b(T, 1)), \ldots , \cont(b(T, n)) )$. Equivalently, it is the space of all vectors in $\mathbb{C}[S_n]$ whose restriction path for the left action of $S_n$ is given by $T$.

Since in this weight space we have $T$ (restriction tableau) fixed and $T'$ (induction tableau) varying, we can denote the weight space by $V_{T, \bullet}$. This indicates that we can implicitly think of it as the row of a matrix, in the same way that we previously described the left Gelfand-Tsetlin module as the column of a matrix.

Describing $V_{T, \bullet}$ as the image of the projection operator $p_T(X_1, \ldots , X_n)$ acting on the left immediately yields
\begin{prop}
$V_{T, \bullet}$ is an invariant subspace for the right action of $S_n$ on $\mathbb{C}[S_n]$.
\end{prop}

\begin{proof}
The left and right actions of $S_n$ commute, so $S_n$ acting on the right commutes with the left action of the $X_i$ and hence of $p_T(X_1, \ldots , X_n)$. In other words, for $\sigma \in S_n$,
\begin{displaymath}
(p_T \cdot v) \cdot \sigma = p_T \cdot (v \cdot \sigma),
\end{displaymath}
which shows that $\sigma$ preserves the image of $p_T$, as asserted.
\end{proof}

Since we are operating in the $\lambda$-isotypic component of $\mathbb{C}[S_n]$, comparing dimensions shows that $V_{T, \bullet}$ must be irreducible.

We mentioned earlier that we can view $p_T(X_1, \ldots , X_n)$ as an element of $\mathbb{C}[S_n]$ as well as an operator. It’s illuminating to identify which element:

\begin{prop} \label{proj-diag}
As an element of $\mathbb{C}[S_n]$, $p_T(X_1, \ldots , X_n) = v_{T,T}$. 
\end{prop}
\begin{proof}
Since we can write $p_T(X_1, \ldots , X_n) = p_T(X_1, \ldots , X_n) \cdot 1$, we know that as an element of $\mathbb{C}[S_n]$, $p_T(X_1, \ldots , X_n)$ is contained in the image of $p_T(X_1, \ldots , X_n)$ viewed as a left multiplication operator acting on $\mathbb{C}[S_n]$. By the above, this means that as an element of $\mathbb{C}[S_n]$, $p_T(X_1, \ldots , X_n)$ lies in the weight space $V_{T, \bullet}$. I.e., it is an eigenvector for left multiplication by $X_1, \ldots , X_n$ with spectrum given by $(\cont(b(T, 1)), \ldots , \cont(b(T, n)) )$.

Now assume by induction that $p_{T^{(n-1)}}(X_1, \ldots , X_{n-1}) = v_{T^{(n-1)}, T^{(n-1)}}$. Then it lies in the (left) Gelfand-Tsetlin module $V_{T^{(n-1)}}$ for $S_{n-1}$. Since we have
\begin{displaymath}
\Ind_{S_{n-1}}^{S_n} V_{T^{(n-1)}} = \mathbb{C}[S_n] \cdot V_{T^{(n-1)}} \subseteq \mathbb{C}[S_n]
\end{displaymath}
and we can write
\begin{displaymath}
p_T(X_1, \ldots , X_n) = p_{T^{(n-1)}, T} (X_n) \cdot p_{T^{(n-1)}}(X_1, \ldots , X_{n-1}),
\end{displaymath}
we can conclude that $p_T(X_1, \ldots , X_n)$ lies in $\Ind_{S_{n-1}}^{S_n} V_{T^{(n-1)}}$.

Now we have just shown that $p_T(X_1, \ldots , X_n)$ has spectrum corresponding to $T$ for the left action of the $X_i$. In the decomposition of $\Ind_{S_{n-1}}^{S_n} V_{T^{(n-1)}}$, $V_T$ is the only component that includes a vector of that spectrum. That vector is $v_{T,T}$, hence $p_T(X_1, \ldots , X_n) = v_{T,T}$.
\end{proof}

Since we had to fix the scaling of $p_T$ to make it an idempotent, but haven’t paid attention to the scaling of $v_{T,T'}$ to this point, we take this opportunity to {\em define} the scaling of $v_{T,T}$ so that Proposition \ref{proj-diag} holds. We can now finish the proof of the induction eigen-lemma (Proposition \ref{ind-lemma}):

\begin{proof}
Since the $X_i$ commute, for any $X_i$ we have
\begin{displaymath}
X_i \cdot p_T(X_1, \ldots , X_n) = p_T(X_1, \ldots , X_n) \cdot X_i,
\end{displaymath}
which means that $p_T(X_1, \ldots , X_n)$ is also an eigenvector for the {\em right} multiplication action of the $X_i$, with the same spectrum. Since the left and right actions of $S_n$ commute, $\mathbb{C}[S_n]$ acting on itself on the left preserves right eigenvectors of the $X_i$ and their spectra. But $\mathbb{C}[S_n] \cdot v_{T,T}$ is all of $V_T$ since $V_T$ is irreducible. Hence all of $V_T$ is an eigenspace for the $X_i$ acting on the right, with spectrum $(\cont(b(T, 1)), \ldots , \cont(b(T, n)) )$.

Fixing an $n-1$-box tableau $T' \subset T$ and looking at the decomposition of $\Ind_{S_{n-1}}^{S_n} V_{T'}$, we see that for all $n$-box tableaux containing $T'$ other than $T$, $n$ appears in a different row and hence must have different content. Thus the condition $X_n = \cont(b(T, n))$ is necessary and sufficient to pick out $V_T$ as a right eigenspace of $X_n$ in this decomposition. This is exactly the statement of the induction eigen-lemma for the Gelfand-Tsetlin module $V_{T'} \subseteq \mathbb{C}[S_{n-1}]$. The general lemma follows since any submodule $V_{\lambda'} \subseteq \mathbb{C}[S_{n-1}]$ of the same isomorphism type can be mapped to a Gelfand-Tsetlin module by an $S_n$-intertwiner, which must preserve the isomorphism type of the components of $\Ind_{S_{n-1}}^{S_n}$.
\end{proof}

Recall that we can think of $V_{T, \bullet}$ as the space of all vectors in $\mathbb{C}[S_n]$ whose restriction path for the left action of $S_n$ is given by $T$. It is now straightforward to prove
\begin{thm} Left-Right / Induction-Restriction Duality:
\begin{enumerate}
\item Each weight space $V_{T, \bullet}$ is also equal to the space obtained by following the induction path specified by $T$ for the {\em right} action of $S_n$.
\item The weight basis for the left action of $S_n$ is also a weight basis for the right action, with the induction and restriction tableaux reversed.
\end{enumerate}
\end{thm}

\begin{proof}
Applying the induction eigen-lemma to each step up the induction chain for the left action, we see that for every $k$ and for any $T$ and $T'$, $v_{T', T}$ is an eigenvector for the right action of $X_k$ with eigenvalue $X_k = \cont(b(T, k))$. In other words, the content vector of $T$ (the induction tableau with respect to the left action of $S_n$) records the spectrum for the right action of $X_1, \ldots , X_n$ on $v_{T', T}$ in the same way that the content vector of $T'$ (the restriction tableau with respect to the left action of $S_n$) recorded the spectrum for the left action.

Now repeat the induction-restriction construction for the right action of $S_n$, generating a right weight basis $v_{T', T}'$. Formulating the above with respect to the right action, the content vector of $T'$ records the spectrum for the left action of $S_n$ on $v_{T', T}'$, and the content vector of $T$ records the spectrum for the right action of $S_n$. But if the contents of $T_1$ and $T_1'$ are equal to the contents of $T_2$ and $T_2'$, then $T_1 = T_1'$ and $T_2 = T_2'$. By the uniqueness of the spectral decomposition, we have $v_{T, T'} = v_{T', T}'$, and both parts of the theorem follow.
\end{proof}

Summarizing the results of this section, we have shown:
\begin{thm} \label{weight-space-thm}
The following are equivalent descriptions of the left Gelfand-Tsetlin module $V_T = V_{\bullet , T} \subseteq \mathbb{C}[S_n]$:
\begin{enumerate}
\item The left $S_n$-module generated by following the induction path corresponding to $T$.
\item The weight space corresponding to $T$ (via content vectors) for $X_1, \ldots , X_n$ acting by right multiplication.
\item The image of $p_T(X_1, \ldots , X_n)$ acting as a right multiplication operator on $\mathbb{C}[S_n]$.
\end{enumerate}
Of course, we obtain a corresponding description of $V_{T, \bullet}$ by switching left and right and switching induction and restriction.
\end{thm}

A good way to think of the equivalence of the first and third items above is to go back to the definition
\begin{displaymath}
p_T(X_1, \ldots , X_n) = \prod_{k=2}^n p_{T^{(k-1)}, T^{(k)}}(X_k),
\end{displaymath}
and note that right multiplication by $p_{T^{(k-1)}, T^{(k)}}(X_k)$ projects $\Ind_{S_{k-1}}^{S_k} V_{T^{(k-1)}}$ onto $V_{T^{(k)}}$. This is precisely the $k-1$-st step in the induction path corresponding to $T$ that generates the Gelfand-Tsetlin module $V_{\bullet, T}$.

\subsection{Endomorphism Algebra}
The simple idea underlying our discussion of duality was to treat induced modules and weight spaces like columns and rows of a matrix, and to use the diagonal $v_{T,T}$ and commutativity of the left and right actions to extract enough information about the spectra to show that the columns and rows were dual to each other. Building on this, the spectral information we now have enables us to make the matrix analogy explicit, by identifying the (two-sided!) weight basis of the regular representation with matrix units in the corresponding endomorphism algebra.

It is well-known that for any finite group $G$, we have a canonical algebra isomorphism
\begin{displaymath}
\mathbb{C}[G] \cong \bigoplus_i \End (V_i),
\end{displaymath}
where $V_i$ ranges over all the irreducible representations of $G$. The isomorphism is given by left multiplication by $g \in G$ on each $\End(V_i)$, and justified by noting this is injective because the regular representation is faithful, then comparing dimensions. We can then identify $\End (V)$ with $V \otimes V^*$ as a $G$-module. Specializing to $G = S_n$, we have
\begin{displaymath}
\mathbb{C}[S_n] \cong \bigoplus_{| \lambda \ = n} \End(V_{\lambda}) \cong \bigoplus_{| \lambda| = n} V_{\lambda} \otimes V_{\lambda}^*,
\end{displaymath}
where in the last decomposition, $S_n$ acts on the left on $V_{\lambda}$ and on the right on $V_{\lambda}^*$, and these actions commute. The right action of $S_n$ on $V_{\lambda}^*$ is given by
\begin{displaymath}
(f \cdot \sigma)(v) = f(\sigma \cdot v),
\end{displaymath}
for $v \in V_{\lambda}$, $f \in V_{\lambda}^*$, $\sigma \in S_n$.

Now let $\{v_T\}$ be the weight basis in $V_{\lambda}$, and let $\{v_T^*\}$ be the dual basis in $V_{\lambda}^*$.

\begin{prop}
$\{v_T^*\}$ is a weight basis for the right action of $X_i \in S_n$ on $V_{\lambda}^*$ with the same weights as $\{v_T\}$, i.e.,
\begin{displaymath}
v_T^* \cdot X_i = \cont(b(T,i)) \cdot v_T^*.
\end{displaymath}
\end{prop}

\begin{proof}
For all $T$ and $T'$, we have
\begin{eqnarray*}
(v_T^* \cdot X_i)(v_{T'}) &=& v_T^*(X_i \cdot v_{T'}) \\
&=& v_T^* (\cont(b(T',i)) \cdot v_{T'}) \\
&=& \cont(b(T',i)) \cdot v_T^* (v_{T'}) \\
&=& \cont(b(T',i)) \cdot \delta_{T, T'} \\
&=& \cont(b(T,i)) \cdot \delta_{T, T'} \\
&=& (\cont(b(T,i)) \cdot v_T^*)(v_{T'}),
\end{eqnarray*}
showing that $v_T^* \cdot X_i = \cont(b(T,i)) \cdot v_T^*$.
\end{proof}

This shows that the basis $\{v_T \otimes v_{T'}^*\}$ for $\bigoplus_{\lambda} V_{\lambda} \otimes V_{\lambda}^*$ is an eigenbasis for the $X_i$ with the same spectrum as the weight basis in $\mathbb{C}[S_n]$, and is therefore the same as that basis. (This is just the basis of matrix units in $\bigoplus_{\lambda} \End(V_{\lambda})$ generated by the weight basis of $V_{\lambda}$.) Since we have an $S_n$-invariant inner product with respect to which the $v_T$ form an orthogonal basis, we can identify $V_{\lambda}^*$ with $V_{\lambda}$ and identify $v_T^*$ with $v_T$. The $S_n$-invariance guarantees that this is an isomorphism of $S_n$-modules, and we recover the weight basis as $\{v_T \otimes v{_T'}\}$. Here $v_T$ and $v_{T'}$ correspond to restriction paths in Young’s lattice for the left and right actions of $S_n$, and our induction-restriction duality analysis above shows that they correspond to induction paths as well.

A natural benefit of this point of view is an easy expression of the identity endomorphism, which is just the sum of the matrix units on the diagonal, summed up over all standard tableaux having $n$ boxes:
\begin{displaymath}
1 = \sum_{|T| = n} v_T \otimes v_T^* = \sum_{|T| = n} v_T \otimes v_T.
\end{displaymath}
In $\mathbb{C}[S_n]$, this becomes
\begin{equation} \label{eq-part-unity}
1 = \sum_{|T| = n} v_{T, T} = \sum_{|T| = n} p_T(X_1, \ldots , X_n).
\end{equation}
We will use this expression later in the paper. One useful consequence to keep in mind is
\begin{cor} \label{cor-orthog}
If $T \neq T'$, then $v_{T, T'}$ is orthogonal to $\id \in \mathbb{C}[S_n]$.
\end{cor}

\begin{rem}
We could also prove equation (\ref{eq-part-unity}) directly, as for all $k \leq n$, and for any $k-1$-box tableau $T'$, we have
\begin{displaymath}
\sum_{T \supset T', \ |T| = k} p_{T', T}(X_k) = 1.
\end{displaymath}
\end{rem}

To sketch this, let $j+1$ be the number of outer corners in $T'$ (one located where we add a box $b$ to $T'$ to obtain $T$, along with $j$ more outer corners $b_1, \ldots , b_j$). We can view the definition of $p_{T', T}(X_i)$ as specifying a polynomial $p_{T',T}(z)$ of degree $j$ that vanishes at $\cont(b_1), \ldots , \cont(b_j)$ and evaluates to $1$ at $\cont(b)$. Summing up $p_{T', T}(z)$ over all $T$, we get a polynomial of degree $j$ that equals $1$ at the $j+1$ points $\cont(b), \cont(b_1), \ldots , \cont(b_j)$ (we have noted previously that these points must be distinct). Hence $\sum_{T \supset T', \ |T| = k} p_{T', T}$ is identically $1$.

Now collect together terms in $\sum_{|T| = n} v_{T, T} = \sum_{|T| = n} p_T$ corresponding to those $T$ where $T^{(n-1)}$ is equal to the same $n-1$-box tableau $T'$. We can factor their contribution to $\sum_{|T| = n} v_{T, T} = \sum_{|T| = n} p_T$ as
\begin{displaymath}
p_{T'} (X_1, \ldots , X_{n-1}) \cdot \sum_{T \supset T', \ |T| = k} p_{T', T} (X_n) = p_{T'} (X_1, \ldots , X_{n-1}) \cdot 1.
\end{displaymath}
This shows that we can reduce
\begin{displaymath}
\sum_{|T| = n} p_T(X_1, \ldots , X_n) = \sum_{|T'| = n-1} p_{T'}(X_1, \ldots , X_{n-1}).
\end{displaymath}
Repeating this reduction $n-1$ times, decreasing the number of boxes by $1$ at each step, yields equation (\ref{eq-part-unity}).

\section{Connection to Young Symmetrizers} \label{sec_young_sym}
Our realization of the Gelfand-Tsetlin module $V_T$ as the image of the orthogonal projector $p_T$ acting on the right on $\mathbb{C}[S_n]$ may be reminiscent of the construction of a version of Young’s natural representations of $S_n$ (see Section \ref{basic-rep-theory}) using {\em Young idempotents}, or {\em Young symmetrizers}, acting on $\mathbb{C}[S_n]$. In this section, we review the basic construction using Young idempotents, as well as the connections between Young’s natural and seminormal forms of the irreducible representations. The results of this section are standard, but it is illuminating to frame them in our setting and to show when the constructions coincide.

A related approach to the connections between the natural and seminormal forms appears in \cite{bib_G}.

\subsection{Another Decomposition of $\mathbb{C}[S_n]$}
Given a standard tableau $T$, define {\em row and column stabilizer subgroups} of $S_n$ by
\begin{eqnarray*}
R_T &=& \Stab(\mbox{\rm rows of }T),\\
C_T &=& \Stab(\mbox{\rm columns of }T).
\end{eqnarray*}
For example, if $T = \ytableaushort{123,45}$, we have
\begin{eqnarray*}
R_T &=& S_{ \{ 1,2,3 \} } \times S_{ \{ 4, 5 \} }, \\
C_T &=& S_{ \{ 1,4 \} } \times S_{ \{ 2,5 \} }.
\end{eqnarray*}
The row and column stabilizer groups are made up of permutations that act within each row or column, but do not interchange elements between rows or columns. Clearly the row and column subgroups are direct products of the stabilizers of each individual row or column, which commute with each other.

Now define the following elements in $\mathbb{C}[S_n]$:
\begin{eqnarray*}
r_T &=& \sum_{\sigma \in R_T} \sigma , \\
c_T &=& \sum_{\sigma \in C_T} (-1)^{\sigma} \sigma,
\end{eqnarray*}
where $(-1)^{\sigma}$ denotes the sign of $\sigma$.
$r_T$ is the {\em row symmetrizer} of $T$, and $c_T$ is the {\em column antisymmetrizer} of $T$. Because the row and column stabilizer subgroups are direct products of subgroups corresponding to the individual rows and columns, $r_T$ and $c_T$ can be factored into mutually commuting terms corresponding to each individual row and column. For example, for $T = \ytableaushort{123,45}$ as above, we can write
\begin{eqnarray*}
r_T &=& ( (1 + (1 \ 2))(1 + (1 \ 3) + (2 \ 3)) ) \cdot (1 + (4 \ 5)), \\
c_T &=& (1 - (1 \ 4)) \cdot (1 - (2 \ 5)).
\end{eqnarray*}

To make a partial analogy with the orthogonal projectors $p_T$, the symmetrizers $r_T$ and $c_T$ can be understood as operators acting on the regular representation that kill certain submodules. To express this, we need to define the {\em lexicographic order} on partitions and tableaux.

\begin{defn}
Given partitions $\lambda$ and $\lambda'$ of $n$, and standard tableaux $T$ and $T'$ of the same shape, we define the {\em lexicographic order} by:
\begin{enumerate}
\item $\lambda > \lambda'$ if $\lambda = (\lambda_1, \ldots , \lambda_l)$, $\lambda' = (\lambda_1', \ldots , \lambda_{l'}')$, and $\lambda_i > \lambda_i'$ for the first $i$ such that $\lambda_i \neq \lambda_i'$.
\item $T > T'$ if $\shape(T^{(i)}) > \shape((T')^{(i)})$ in lexicographic order on partitions for the first $i$ such that $T^{(i)} \neq (T')^{(i)}$.
\end{enumerate}
\end{defn}
Intuitively, $\lambda > \lambda'$ means that the first differing row is longer in $\lambda$ than in $\lambda'$, and $T > T'$ means that the first differing number appears in an earlier (i.e., higher) row in $T$ than in $T'$.

It is straightforward to show that:
\begin{prop}
For a fixed shape $\lambda$, the maximal element in the lexicographic order on standard tableaux of shape $T$ is the tableau $T_{\smax}$ with $1, 2, \ldots , n$ running in order along the rows of $\lambda$, and the minimal element is the tableau $T_{\smin}$ with $1, 2, \ldots , n$ running in order down the columns of $\lambda$.
\end{prop}
For example, for $\lambda = (3,2)$, we have $T_{\smax} = \ytableaushort{123,45}$ and $T_{\smin} = \ytableaushort{135,24}$.

The following is fundamental in the application of Young symmetrizers:
\begin{prop} \label{prop-order}
If $T > T'$, or if $\lambda > \lambda'$ and $T$ and $T'$ are a pair of standard tableaux of shape $\lambda$ and $\lambda'$, then there is a pair of numbers $(i , j)$ such that $i$ and $j$ appear in the same row of $T$ and the same column of $T'$.
\end{prop}

\begin{proof}
Whether $\lambda > \lambda'$ or $T > T'$, the first row of $T'$ cannot be longer than the first row of $T$. Either one of the entries in the first row of $T$ appears in a lower row in $T'$, or else the first rows of $T$ and $T'$ are the same ($T'$ has exactly the same entries as $T$, which must be increasing in both tableaux because $T$ and $T'$ are standard).

In the former case, the first entry that appears in a lower row in $T'$ than in $T$, together with the first row entry above it in $T'$ are the pair we seek, since they are in the same column in $T'$, and both appear in the first row of $T$.

In the latter case, move on to the second row, note that in this case the second row of $T'$ can’t be longer than the second row of $T$, and repeat the argument above. If the second rows of $T$ and $T'$ are the same, move on to the third row, and so on. Either this process eventually terminates, producing the pair we seek, or $T = T'$ (and hence $\lambda = \lambda'$), contradicting one of our assumptions.
\end{proof}

The application of Proposition \ref{prop-order} to Young symmetrizers killing subspaces of $\mathbb{C}[S_n]$ is:
\begin{prop} \label{prop-annihilate}
Let $(i, j)$ be a pair of numbers that appear in the same row of $T$ and the same column of $T'$. Then
\begin{displaymath}
r_T \cdot c_{T'} = c_{T'} \cdot r_T = 0.
\end{displaymath}
\end{prop}

\begin{proof}
We give a ``Hecke algebra-style'' argument. If $(i_1, \ldots , i_k)$ are the entries of a particular row of $T$, and $(j_1, \ldots , j_k)$ is {\em any} reordering of those entries, we can factor the symmetrizer of that row as
\begin{eqnarray*}
\sum_{\sigma \in S_{\{i_1,\ldots , i_k\}}} \sigma &= &\sum_{\sigma \in S_{\{j_1,\ldots , j_k\}}} \sigma \\
&=& (1 + (j_1 \ j_2))(1 + (j_1 \ j_3) + (j_2 \ j_3)) \cdots (1 + (j_1 \ j_n) + \cdots + (j_{n-1} \ j_n)).
\end{eqnarray*}
The terms in this factorization commute. We have a similar factorization for $c_T$, with a minus sign in front of every transposition.

This implies that given a pair $(i, j)$ in the same row of $T$, we can factor the symmetrizer of that row as a product of $(1 + (i \ j))$ with a term that commutes with $(1 + (i \ j))$. Bringing $(1 + (i \ j))$ all the way across to the right, using its commutativity with the symmetrizers of any other row, we find that we can write
\begin{displaymath}
r_T = r' \cdot (1 + (i \ j))
\end{displaymath}
for some $r'$. Similarly, factoring $(1 - (i \ j))$ out of $c_{T'}$ and bringing it all the way over to the left, we find that we can write
\begin{displaymath}
c_{T'} = (1 - (i \ j)) \cdot c'
\end{displaymath}
for some $c'$. Since $(1 + (i \ j)) \cdot (1 - (i \ j)) = 0$, we conclude that
\begin{displaymath}
r_T \cdot c_{T'} = r' \cdot (1 + (i \ j)) \cdot (1 - (i \ j)) \cdot c' = 0.
\end{displaymath}
Repeating this argument but placing $(1 + (i \ j))$ at the left of $r_T$ and placing $(1 - (i \ j))$ at the right of $c_{T'}$ enables them to annihilate each other in the product $c_{T'} \cdot r_T$, showing that $c_{T'} \cdot r_T = 0$ as well.
\end{proof}

Now let $T$ be a standard tableau of shape $\lambda$, and look at the products
\begin{displaymath}
\begin{array}{ccc}
c_T \cdot r_T & \mbox{\rm and} & r_T \cdot c_T.
\end{array}
\end{displaymath}
If $T'$ is another tableau of shape $\lambda$ such that $T' < T$, or if $T'$ is a tableau of a different shape $\lambda'$ such that $\lambda' < \lambda$, then by Proposition \ref{prop-annihilate}, we have
\begin{eqnarray*}
(c_T \cdot r_T) \cdot (c_{T'} \cdot r_{T'}) &=& c_T \cdot (r_T \cdot c_{T'}) \cdot r_{T'} = 0, \\
(r_{T'} \cdot c_{T'}) \cdot (r_T \cdot c_T) &=& r_{T'} \cdot (c_{T'} \cdot r_T) \cdot c_T = 0.
\end{eqnarray*}
Equivalently, if $T'$ is another tableau of shape $\lambda$ such that $T' > T$, or if $T'$ is a tableau of a different shape $\lambda'$ such that $\lambda' > \lambda$, we have
\begin{eqnarray*}
(c_{T'} \cdot r_{T'}) \cdot (c_T \cdot r_T) &=& c_{T'} \cdot (r_{T'} \cdot c_T) \cdot r_T = 0, \\
(r_T \cdot c_T) \cdot (r_{T'} \cdot c_{T'}) &=& r_T \cdot (c_T \cdot r_{T'}) \cdot c_{T'} = 0.
\end{eqnarray*}

We will assume the fundamental result in the theory of Young symmetrizers that the products $c_T r_T$ and $r_T c_T$ project $\mathbb{C}[S_n]$ onto irreducible representations with isomorphism class given by the shape of $T$. That is, $(c_T r_T)^2 = k \cdot c_T r_T$, and $(r_T c_T)^2 = k' \cdot r_T c_T$, for constants $k$ and $k'$ we will not particularly need to pay attention to. Ignoring these constants, right multiplication by $c_T r_T$ and $r_T c_T$ projects $\mathbb{C}[S_n]$ onto irreducible left $S_n$-submodules
\begin{eqnarray*}
W_T &=& \mathbb{C}[S_n] \cdot c_T r_T, \\
W_T' &=& \mathbb{C}[S_n] \cdot r_T c_T,
\end{eqnarray*}
which are not equal but which are both isomorphic to $V_{\lambda}$, where $\lambda$ is the shape of $T$. We can also multiply by $c_T r_T$ and $r_T c_T$ on the left, and this projects $\mathbb{C}[S_n]$ onto irreducible right $S_n$-submodules isomorphic to $V_{\lambda}$. This should be reminiscent of our previous use of the orthogonal projectors $p_T$, and the rest of this section works out some relationships between these various projections.

\begin{rem} \label{rem-specht-proj}
Let $T = T_{\smax}$ be the maximal tableau in the lexicographic order, with $1, \ldots , \lambda_1$ in the first row, $\lambda_1 + 1, \ldots , \lambda_1 + \lambda_2$ in the second row, and so on. There is an isomorphism between $\mathbb{C}[S_n] \cdot c_T r_T$ and the Specht polynomial realization of $V_{\lambda}$, which maps
\begin{eqnarray*}
c_T r_T & \mapsto & c_T r_T \cdot (z_{\lambda_1+1} \cdots z_{\lambda_1 + \lambda_2}) \cdot (z_{\lambda_1 + \lambda_2 +1}^2 \cdots z_{\lambda_1 + \lambda_2 + \lambda_3}^2) \cdot (z_{\lambda_1 + \lambda_2 + \lambda_3 + 1}^3 \cdots ) \cdots \\
&=& \Delta_{\lambda_1'} (z_1, z_{\lambda_1+1}, \ldots , z_{n-\lambda_l+1}) \cdot \Delta_{\lambda_2'} (z_2, z_{\lambda_1+2}, \ldots , z_{n-\lambda_l - \lambda_{l-1}+2}) \cdots .
\end{eqnarray*}
\end{rem}
To make this clearer with an example, if we let $T = \ytableaushort{123,45}$, then we map
\begin{displaymath}
c_T r_T =  (1 - (1 \ 4)) \cdot (1 - (2 \ 5)) \cdot ( (1 + (1 \ 2))(1 + (1 \ 3) + (2 \ 3)) ) \cdot (1 + (4 \ 5))
\end{displaymath}
to
\begin{eqnarray*}
c_T r_T \cdot z_4 z_5 &=& (1 - (1 \ 4)) (1 - (2 \ 5)) ( (1 + (1 \ 2))(1 + (1 \ 3) + (2 \ 3)) ) (1 + (4 \ 5)) \cdot z_4 z_5 \\
&=& 2 \cdot 6 \cdot c_T \cdot z_4 z_5 \\
&=& 12 \cdot (z_4 - z_1)(z_5 - z_2).
\end{eqnarray*}
In the language of Section \ref{basic-rep-theory}, this is just the Specht polynomial $f_T(z_1, \ldots , z_5)$ (up to scalars). Note also that $r_T$ acts as a scalar on $z_4 z_5$, as the variables corresponding to each row have the same degree. In the second half of this paper, we will reinterpret and extend this isomorphism to all of $\mathbb{C}[S_n]$.

In analogy with our analysis of $p_T(X_1, \ldots , X_n)$, which we will write simply as $p_T$ going forward, let us look at the kernels of the projectors $c_T r_T$ and $r_T c_T$ acting on the right on $\mathbb{C}[S_n]$. Let $T$ be a standard tableau of shape $\lambda$. Since we know that $c_T r_T$ and $r_T c_T$ are projections onto an irreducible representation of isomorphism class $\lambda$, by Schur’s lemma their restrictions to the isotypic component of any other isomorphism class must be $0$. I.e., just like $p_T$ in the preceding section, $c_T r_T$ and $r_T c_T$ kill any irreducible submodule of isomorphism type $\lambda'$ for $\lambda' \neq \lambda$. (This can also be shown directly using Proposition \ref{prop-annihilate} and a careful extension of the argument that we will give next.)

Now let’s look inside the isotypic component of $\lambda$ inside $\mathbb{C}[S_n]$. Proposition \ref{prop-annihilate} ensures that if $T' > T$ in lexicographic order, with both having shape $\lambda$, then $c_T$, and hence $c_T r_T$, kills $W_{T'} = \mathbb{C}[S_n] \cdot c_{T'} r_{T'}$. Similarly, if $T' < T$, then $r_T$, and hence $r_T c_T$, kills $W_{T'}' = \mathbb{C}[S_n] \cdot r_{T'} c_{T'}$. This implies that $W_{T} \cap W_{T'} = \{0\}$ for every $T$ and $T'$ of shape $\lambda$: if we pick $T'$ to be the larger one out of $T$ and $T'$ in lexicographic order, then $c_T r_T$ acts as a scalar on $W_T$ and kills $W_{T'}$, so $W_T$ and $W_{T'}$ must intersect trivially. Similarly, $W_{T}' \cap W_{T'}' = \{0\}$ as picking $T$ to be the larger of the two in lexicographic order implies $r_T c_T$ acts as a scalar on $W_T$ and kills $W_{T'}$.

The intersection condition implies that the sums of subspaces
\begin{displaymath}
\begin{array}{ccc}
\sum_T W_T & \mbox{\rm and} & \sum_T W_T'
\end{array}
\end{displaymath}
are direct sums of the individual $W_T$ and $W_T'$ inside $\mathbb{C}[S_n]$. This represents the entire $\lambda$-isotypic component inside $\mathbb{C}[S_n]$ since the multiplicity of the latter is the number of standard tableaux $T$ of shape $\lambda$. This shows that we have direct sum decompositions
\begin{displaymath}
\mathbb{C}[S_n] \cong \bigoplus_{T: \ |T| = n} W_T \cong \bigoplus_{T: \ |T| = n} W_T',
\end{displaymath}
analogous to the direct sum decomposition of $\mathbb{C}[S_n]$ into Gelfand-Tsetlin modules $V_T$ in Proposition \ref{direct_sum_gz}. However, note that here the components $W_T$ and $W_T'$ do not need to be orthogonal. Moreover, unlike $\{p_T\}$, neither $\{c_T r_T\}$ nor $\{r_T c_T\}$ are sets of mutually orthogonal projectors, except for special shapes.

\subsection{Lowest Weight Vector in $V_{\lambda}$}
We will work briefly with the Specht polynomial realization of the irreducible $S_n$-module $V_{\lambda}$. Let $T = T_{\smin}$ be the minimal standard tableau of shape $\lambda$ in the lexicographic order. $T$ is the tableau with $1, \ldots , \lambda_1' $ running down the first column of $\lambda$ ($\lambda_1'$ is the length of the column), $\lambda_1' + 1, \ldots, \lambda_1' + \lambda_2'$ running down the second column, and so on. The associated Specht polynomial is
\begin{displaymath}
f_{T_{\smin}}(z_1, \ldots , z_n) = \Delta_{\lambda_1'} (z_1, \ldots , z_{\lambda_1'}) \cdot \Delta_{\lambda_2'} (z_{\lambda_1'+1}, \ldots , z_{\lambda_1'+\lambda_2'}) \cdots \Delta_{\lambda_l'} (z_{n-\lambda_l'+1}, \ldots , z_n).
\end{displaymath}
For example, if $\lambda = (3,2)$, we have $T_{\smin} = \ytableaushort{135,24}$, and $f_{T_{\smin}} = (z_2-z_1)(z_4-z_3)$.

The basic statement is that for $T = T_{\smin}$, the corresponding weight vector and Specht polynomial (we continue to ignore scalars) coincide:
\begin{prop}
$f_{T_{\smin}}$ is the weight vector $v_{T_{\smin}}$ corresponding to the tableau $T_{\smin}$ in $V_{\lambda}$.
\end{prop}

\begin{proof}
Let $T = {T_{\smin}}$. We will compute the action of the $X_i$ on $f_T$ explicitly, building up $T$ a column at a time. To warm up, if $T$ is a single column containing $1, \ldots , n$, $f_T$ is a completely antisymmetric function, on which any transposition acts by $-1$. Starting with $X_1 = 0$, we appeal to the relation
\begin{equation} \label{eq-ind-hecke}
X_{k+1} = \sigma_k X_k \sigma_k + \sigma_k.
\end{equation}
It implies that $f_T$ is an eigenvector of $X_{k+1}$ if it is an eigenvector of $X_k$, with the eigenvalue dropping by $1$ between $X_k$ and $X_{k+1}$ because of the extra $\sigma_k$ term, which acts by $-1$. This corresponds exactly to the content of the corresponding boxes of $T$: the top left box has content $0$, and then each step down the column increases the row index by $1$, and hence decreases the content by $1$.

Now add a single box in the second column, so that our shape is $\lambda = (2, 1^{n-2})$, and $T$ has $1, \ldots , n-1$ running down in the first column and $n$ in the single box in the second. We have
\begin{eqnarray*}
f_T(z_1, \ldots , z_n) &=& f_{T^{(n-1)}}(z_1, \ldots , z_{n-1}) \\
&=&  f_{(1^{n-1})}(z_1, \ldots , z_{n-1}) \\
&=& \sum_{\sigma \in S_{n-1}} (-1)^{\sigma} \cdot z_2 z_3^2 \ldots z_{n-1}^{n-2},
\end{eqnarray*}
with no antisymmetrization corresponding to the second column as there aren’t enough variables to antisymmetrize. The action of $X_1, \ldots , X_{n-1}$ is just what we computed before. As for $X_n$, we claim that it acts by $1$. To see this, look at the action of the full antisymmetrizer
\begin{displaymath}
\sum_{\sigma \in S_n} (-1)^{\sigma} \sigma = (1 - X_n) \cdot \sum_{\sigma \in S_{n-1}} (-1)^{\sigma} \sigma.
\end{displaymath}
$\sum_{\sigma \in S_n} (-1)^{\sigma} \sigma$ acts on $f_T$ by $0$, as the maximum power of any variable in $f_T$ is $n-2$, which means every monomial is missing two variables ($z_n$ and another one) and thus is killed by the complete antisymmetrizer. $\sum_{\sigma \in S_{n-1}} (-1)^{\sigma} \sigma$ acts on $f_T$ as a positive scalar, using the expression for $f_T$ above and the fact that $\sum_{\sigma \in S_{n-1}} (-1)^{\sigma} \sigma$ is an idempotent up to scalars. Hence $(1 - X_n) \cdot f_T = 0$, so $X_n$ has eigenvalue $1$, which is precisely the content of the box containing $n$.

If $\lambda$ is any two-column diagram, $T = T_{\smin}$ has $1, \ldots , l$ running down the first column and $l+1, \ldots , n$ running down the second. The action of $X_1, \ldots , X_l$ and $X_{l+1}$ is exactly as we computed above, because they do not reach beyond the first box in the second column. For the action of $X_{l+2}, \ldots , X_n$, we just need to go back to equation (\ref{eq-ind-hecke}). For $k \geq l+1$, $\sigma_k$ acts entirely on the second column, which is completely antisymmetric, and hence $\sigma_k$ acts by $-1$. Hence our inductive method of computing eigenvalues for the first column applies to the second as well. Anchoring the second column with $X_{l+1}$ acting by $1$, we conclude $f_T$ is also an eigenvector for each of $X_{l+2}, X_{l+3}, \ldots , X_n$, with each eigenvalue being $1$ less than the previous one as we proceed down the column.

Moving to three columns, add a single box in the third column first. Again, we have $f_T = f_{T^{(n-1)}}$. Say the lengths of the first two columns are $l_1$ and $l_2$, so $n = l_1 + l_2 + 1$. The action of $X_1, \ldots , X_{l_1 + l_2}$ is as computed before, and we just need to compute the action of $X_n$.

Write
\begin{displaymath}
X_n = ((1 \ n) + \cdots (l_1 \ n)) + ((l_1 + 1 \ n) + \cdots + (l_1 + l_2 \ n)) = \overline{X}_1 + \overline{X}_2,
\end{displaymath}
and consider the action of $\overline{X}_1$ and $\overline{X}_2$ on
\begin{displaymath}
f_T(z_1, \ldots , z_n) = \prod_{1 \leq i, j \leq l_1} (z_j - z_i) \cdot \prod_{l_1 + 1 \leq i', j' \leq l_1 + l_2} (z_{j'} - z_{i'}) = \overline{f}_1 \cdot \overline{f}_2,
\end{displaymath}
where we use the last equality as the definition of $\overline{f}_1$ and $\overline{f}_2$. $\overline{X}_1$ doesn’t touch $\overline{f}_2$ and $\overline{X}_2$ doesn’t touch $\overline{f}_1$. Hence we can apply the argument from the single-column-plus-one-box case to conclude that
\begin{eqnarray}
\overline{X}_1 \cdot (\overline{f}_1 \cdot \overline{f}_2) &=& (\overline{X}_1 \cdot \overline{f}_1) \cdot \overline{f}_2 = \overline{f}_1 \cdot \overline{f}_2, \\
\overline{X}_2 \cdot (\overline{f}_1 \cdot \overline{f}_2) &=& \overline{f}_1 \cdot (\overline{X}_2 \cdot \overline{f}_2) = \overline{f}_1 \cdot \overline{f}_2, \\
X_n \cdot f_T &=& (\overline{X}_1 + \overline{X}_2) \cdot (\overline{f}_1 \cdot \overline{f}_2) = 2(\overline{f}_1 \cdot \overline{f}_2) = 2 f_T.
\end{eqnarray}
Now we are essentially done. We apply (\ref{eq-ind-hecke}) to extend the above computation to cover any $3$-column diagram, starting at the top of the column with the eigenvalue $2$ we just computed and decreasing the eigenvalue by $1$ with each step down the column. For diagrams with $m$ columns plus a single box in the $m+1$-st column, with the $i$-th column having length $l_i$ and $n = l_1 + \cdots + l_m + 1$, write
\begin{eqnarray*}
\overline{X}_1 &=& (1 \ n) + \cdots + (l_1 \ n) \\
\overline{X}_2 &=& (l_1 + 1 \ n) + \cdots + (l_1 + l_2 \ n) \\
& \vdots & \\
\overline{X}_m &=& (l_1 + \cdots + l_{m-1}+1 \ n) + \cdots + (n-1 \ n).
\end{eqnarray*}
Then
\begin{equation} \label{eq-column-decomp}
X_n = \overline{X}_1 + \cdots + \overline{X}_m.
\end{equation}
Each $\overline{X}_i$ acts on just the $i$-th column of $T$ as above, and hence acts by $1$ on $f_T$. Hence, $X_n \cdot f_T = m \cdot f_T$, and the extension down the column using (\ref{eq-ind-hecke}) takes care of the rest.
\end{proof}
As $T_{\smin}$ is the minimal standard tableau of shape $\lambda$ with respect to the lexicographic order, we will refer to $v_{T_{\smin}} = f_{T_{\smin}}$ as the {\em lowest weight vector} in $V_{\lambda}$.

We can repeat this analysis inside the regular representation $\mathbb{C}[S_n]$:
\begin{prop}
Let $T = T_{\smin}$ be the minimal tableau corresponding to the partition $\lambda$. Then we have
\begin{equation} \label{eq-lw-sym}
v_{T, T} = c_T r_T c_T
\end{equation}
as elements of $\mathbb{C}[S_n]$.
\end{prop}

\begin{proof}
Compute the action of the $X_i$ on $c_T r_T$ from the left and on $r_T c_T$ from the right. This works the same way as the computation we just did for the lowest weight Specht polynomial. Look at the left action: for a single column, $c_T r_T = \sum_{S_n} (-1)^{\sigma} \sigma$, any transposition acts on the left by $-1$, so (\ref{eq-ind-hecke}) implies $X_{k+1} = X_k - 1$, with $X_1 = 0$ by definition. Adding a single box in the second column as above, we want to show
\begin{displaymath}
(1 -  X_n) \cdot \left( \sum_{\sigma \in S_{n-1}} (-1)^{\sigma} \sigma \right) \cdot (1 + (1 \ n)) = 0.
\end{displaymath}
We have
\begin{displaymath}
(1 -  X_n) \cdot \left( \sum_{\sigma \in S_{n-1}} (-1)^{\sigma} \sigma \right) = \sum_{\sigma \in S_n} (-1)^{\sigma} \sigma,
\end{displaymath}
and $(1 \ n)$ acts on this by multiplication from the right by $-1$, so $1 +  (1 \ n)$ acts by $0$.

The rest follows as before. We treat general two-column diagrams by using (\ref{eq-ind-hecke}) inductively down the second column to show that for $k > l$, where $l$ is the length of the first column, $X_{k+1} = X_k - 1$ again. For general $T$, first treat the case where the last column of $T$, say the $m$-th, has a single box. Write $X_n = \overline{X}_1 + \cdots + \overline{X}_m$ as in (\ref{eq-column-decomp}), where, for $i = 1, \ldots , m$, 
\begin{displaymath}
\overline{X}_i = (l_1 + \cdots + l_{i-1}+1 \ n) + \cdots + (l_1 + \cdots + l_i \ n)
\end{displaymath}
acts on the $i$-th column only. We show each $\overline{X}_i$ acts by $1$. Factor 
\begin{displaymath}
r_T = r_1 \cdots r_l,
\end{displaymath}
where $r_j$ is the symmetrizer of the $j$-th row, and factor $r_1$, the symmetrizer of the first row as 
\begin{displaymath}
r_1 = (1 + (l_1 + \cdots + l_{i-1} + 1 \ n)) \cdot r_1'.
\end{displaymath}
Here $l_1 + \cdots + l_{i-1} + 1$ is the first entry in the $i$-th column of $T$. Denote the stabilizer of the $i$-th column by $C_i$ and the stabilizer of the $i$-th column along with $n$ by $C_{i, n}$. Denote the full antisymmetrizer corresponding to $C_i$ by $c_i$. Then
\begin{displaymath}
c_T = c_1 \cdots c_m,
\end{displaymath}
where the $c_j$ all commute with each other and all $c_j$ with $j \neq i$ commute with $\overline{X}_i$. As in the two-column case, we have
\begin{displaymath}
(1 -  \overline{X}_i) \cdot \left( \sum_{\sigma \in C_i} (-1)^{\sigma} \sigma \right) = \sum_{\sigma \in C_{i,n}} (-1)^{\sigma} \sigma.
\end{displaymath}
(So $c_i$ commutes with $\overline{X}_i$ as well.) Then we find
\begin{eqnarray*}
(1 - \overline{X}_i) \cdot c_T r_T &=& (1 - \overline{X}_i) \cdot c_1 \cdots c_m \cdot (1 + (l_1 + \cdots + l_{i-1} + 1 \ n)) \cdot r_1' \cdot r_2 \cdots r_l \\
&=& c_1 \cdots \hat{c}_i \cdots c_m \cdot c_{i, n} \cdot (1 + (l_1 + \cdots + l_{i-1} + 1 \ n)) \cdot r_1' \cdot r_2 \cdots r_l \\
&=& 0,
\end{eqnarray*}
where the notation $\hat{c}_i$ means we remove $c_i$ from the product, and in the final step we use $c_{i, n} \cdot (1 + (l_1 + \cdots + l_{i-1} + 1 \ n)) = 0$. Extending this to compute eigenvalues down the final column if it has more than one box is a matter of using (\ref{eq-ind-hecke}) as usual.

This calculation proves that $X_i \cdot c_T r_T = \cont(b(T,i)) \cdot c_T r_T$ for all $i$. The same calculation from the right side shows that $r_T c_T \cdot X_i = \cont(b(T,i)) \cdot r_T c_T$ as well. By uniqueness of weights when we consider both the left and right sided actions, the proposition follows.
\end{proof}
Essentially, our proof shows that $c_T r_T \mathbb{C}[S_n]$ is the weight space corresponding to $T$ for the left action of $S_n$, which by Theorem \ref{weight-space-thm} is equal to the right Gelfand-Tsetlin $S_n$-module $V_{T, \bullet}$. Similarly, $\mathbb{C}[S_n] r_T c_T$ is the weight space corresponding to $T$ for the right action of $S_n$, and hence equal to the left Gelfand-Tsetlin $S_n$-module $V_{\bullet, T}$. These intersect in $c_T r_T c_T$, which is therefore equal to $v_{T, T}$.

\subsection{Triangularity} \label{sec-triang}
Any permutation $\sigma \in S_n$ takes a Specht polynomial $f_T$ to another Specht polynomial $f_{\sigma \cdot T}$, as permuting the variables $z_1, \ldots , z_n$ corresponds exactly to permuting the entries of the tableau $T$. Starting with $T = T_{\smin}$ and keeping track of the action of $S_n$ in terms of both the Specht basis and the weight basis, we can relate the two bases to each other beyond the lowest weight vector.

To begin with, let us show how to connect any standard tableau $T$ to the minimal tableau $T_{\smin}$. This can be done using a single permutation $\sigma$, which works well for the Specht basis, but isn’t helpful for computing the action on the weight basis. To compute the latter, we break up $\sigma$ and connect $T$ to $T_{\min}$ using a sequence of adjacent transpositions, which, as we have seen, act {\em locally} on the weight basis.

First, we look at how permutations act directly on tableaux:
\begin{prop}
Let $T$ be any standard tableau of shape $\lambda$. Then there is a sequence of adjacent transpositions $\sigma_{i_1}, \ldots , \sigma_{i_k}$ such that:
\begin{enumerate}
\item $T_{\smin} = \sigma_{i_1} \cdots \sigma_{i_k} \cdot T$,
\item For any $j < k$, $\sigma_{i_j}$ acting on $\sigma_{i_{j+1}} \cdots \sigma_{i_k} \cdot T$ moves $i$ down and $i+1$ up, where $i$ is such that $\sigma_{i_j} = (i \ i+1)$.
\item For any $j \leq k$, $\sigma_{i_j} \cdot (\sigma_{i_{j+1}} \cdots \sigma_{i_k} \cdot T) < \sigma_{i_{j+1}} \cdots \sigma_{i_k} \cdot T$ in lexicographic order. 
\end{enumerate}
\end{prop}

\begin{proof}
Let $i$ be the number written in the last box of the last column of $T$ (where $n$ would appear in $T_{\smin}$). Since there are no columns to the right of $i$ and no boxes immediately below, $i+1$ must lie to the left of $i$, which means it must lie in a lower row. Let the first transposition to act on $T$ switch $i$ with $i+1$, moving $i$ down and $i+1$ up. Since $T$ and $(i \ i+1) \cdot T$ agree on $1, \ldots, i-1$ and $i$ appears in a higher row in $T$ than in $(i \ i+1) \cdot T$, we have $(i \ i+1) \cdot T < T$ in lexicographic order.

Next, switch $i+1$ with $i+2$, then $i+2$ with $i+3$, and so on, terminating when we switch $n-1$ with $n$ to move $n$ to the box where it belongs in $T_{\smin}$. Each step is an adjacent transposition that moves a higher number up and a lower number down, meaning that each new tableau is below the one before it in lexicographic order.

Once $n$ has been placed where it belongs in $T_{\smin}$, remove it from $T$ and repeat the same process to place $n-1$ in the last box of the last column of the shape that remains. Then remove $n-1$ and repeat again to place $n-2$ in the last box of the last column of the shape that remains, and so on.
\end{proof}

Inverting this, we can write
\begin{displaymath}
T = \sigma_{i_k} \sigma_{i_{k-1}} \cdots \sigma_{i_2} \sigma_{i_1} \cdot T_{\smin},
\end{displaymath}
where each adjacent transposition $\sigma_{i_j} = (i \ i+1)$ we apply moves $i$ up and $i+1$ down, meaning that it moves the running tableau $\sigma_{i_{j-1}} \cdots \sigma_{i_1} \cdot T_{\smin}$ up in lexicographic order.

Consider the permutation
\begin{displaymath}
\sigma_T = \sigma_{i_k} \sigma_{i_{k-1}} \cdots \sigma_{i_2} \sigma_{i_1}
\end{displaymath}
such that $T = \sigma_T \cdot T_{\smin}$. Define the {\em column reading word} of a tableau as the arrangement of $1, \ldots , n$ obtained by reading the entries down each column and concatenating the columns, left to right. Thus the column reading word of $T_{\smin}$ is $(1, \ldots , n)$, and the column reading word of $T = \sigma_T \cdot T_{\smin}$ is $(\sigma_T(1), \ldots , \sigma_T(n))$.

Translating to words the fact that $\sigma_{i_j} = (i \ i+1)$ applied to $\sigma_{i_{j-1}} \cdots \sigma_{i_1} \cdot T_{\smin}$ moves $i$ up (and to the right) and $i+1$ down (and to the left), we find that $\sigma_{i_j}$ applied to $\sigma_{i_{j-1}} \cdots \sigma_{i_1} \cdot (1, \ldots , n)$ also moves $i$ to the right and $i+1$ to the left. In other words, $\sigma_{i_j}$ increases the number of inversions of the corresponding column reading word by $1$. Since we can view $\sigma_T = \sigma_{i_k} \sigma_{i_{k-1}} \cdots \sigma_{i_2} \sigma_{i_1}$ as the column reading word of $T$ when we write out $\sigma_T$ as a function from $\{1, \ldots , n\}$ to itself, this proves
\begin{lemma} We have
\begin{enumerate}
\item The number of inversions in $\sigma_T = \sigma_{i_k} \sigma_{i_{k-1}} \cdots \sigma_{i_2} \sigma_{i_1}$ is $k$.
\item This expansion of $\sigma_T$ in terms of adjacent transpositions has minimal length.
\end{enumerate}
\end{lemma}

\begin{proof}
We just proved the first statement. The second follows from the first since any adjacent transposition can only increase the total number of inversions by $1$, and hence any shorter expansion would have fewer inversions.
\end{proof}

Now we are ready to look at the action of $\sigma_T = \sigma_{i_k} \sigma_{i_{k-1}} \cdots \sigma_{i_2} \sigma_{i_1}$ on $v_{T_{\smin}} = f_{T_{\smin}}$. By Proposition \ref{prop_local_action}, we have:
\begin{displaymath}
\sigma_i \cdot v_T = c_1 v_T + c_2 v_{\sigma_i \cdot T},
\end{displaymath}
for constants $c_1$ and $c_2$, where the second term on the right is non-zero if and only if $\sigma_i \cdot T$ is standard.

Iterating this to compute the action of $\sigma_{i_k} \sigma_{i_{k-1}} \cdots \sigma_{i_2} \sigma_{i_1}$, we find that $f_T = \sigma_T \cdot f_{T_{\smin}} = \sigma_T \cdot v_{T_{\smin}}$ is a linear combination of elements of the form
\begin{displaymath}
v_{\sigma_{j_{\kappa}} \sigma_{j_{\kappa-1}} \cdots \sigma_{j_2} \sigma_{j_1} \cdot T_{\smin}},
\end{displaymath}
where $j_{\kappa} j_{\kappa-1} \cdots j_2 j_1$ is a subword of $i_k i_{k-1} \cdots i_2 i_1$. The subword property implies that each product $\sigma_{j_{\kappa}} \sigma_{j_{\kappa-1}} \cdots \sigma_{j_2} \sigma_{j_1}$ lies below $\sigma = \sigma_{i_k} \sigma_{i_{k-1}} \cdots \sigma_{i_2} \sigma_{i_1}$ in the {\em Bruhat order} on $S_n$. This transfers to $\sigma_{j_{\kappa}} \sigma_{j_{\kappa-1}} \cdots \sigma_{j_2} \sigma_{j_1} \cdot T_{\smin}$ lying below $T = \sigma_T \cdot T_{\smin}$ in the induced Bruhat order on standard tableaux.

It is standard that the lexicographic order is a total order on standard tableaux that is a {\em refinement} of the Bruhat order on tableaux induced from the Bruhat order on permutations. This immediately implies
\begin{prop} \label{prop-bruhat-exp}
For any standard tableau $T$, the Specht polynomial $f_T$ is a linear combination of weight vectors $v_T$ and $v_{T_1}, \ldots , v_{T_r}$ where for all $i$, $T_i < T$ in lexicographic order.
\end{prop}

\begin{cor}
If we order the standard tableaux of shape $\lambda$ in lexicographic order, from $T_{\smin}$ to $T_{\smax}$, the transition matrix between the Specht basis $\{f_T\}$ and weight basis $\{v_T\}$ is upper triangular.
\end{cor}

\begin{cor}
For any partition $\lambda$ of $n$, we can write the weight basis $\{v_T\}$ of $V_{\lambda}$ up to scalars as
\begin{equation} \label{eq-weight-irred}
v_T = p_T \cdot \sigma_T \cdot v_{T_{\min}}. 
\end{equation}
\end{cor}

\begin{proof}
Recall that $p_T$ acts in $V_{\lambda}$ as orthogonal projection onto $v_T \in V_T$. By the above, $\sigma_T \cdot v_{T_{\min}}$ is a linear combination of $v_T$ (non-zero coefficient) and a sum of terms $v_{T_i}$ for $T_i < T$. The latter are all killed by $p_T$ since the weight basis is orthogonal, so only the $v_T$ term remains.
\end{proof}

\subsection{Weight Vector Formula in $\mathbb{C}[S_n]$} \label{sec-weight-vector}
We can generalize (\ref{eq-weight-irred}) to an expression for any weight vector $v_{T, T'}$ in $\mathbb{C}[S_n]$, as follows. Given standard tableaux $T$ and $T'$ of the same shape $\lambda$, let $\sigma_{T, T'} \in S_n$ be the permutation such that $\sigma_{T, T'} \cdot T' = T$. (Explicitly, to compute $\sigma_{T, T'}(i)$, we find the box in which $i$ appears in $T'$, and take $\sigma_{T, T'}(i)$ to be the number that appears in that same box in $T$.) Let
\begin{displaymath}
\sigma_{T, T'} = \sigma_{i_k} \sigma_{i_{k-1}} \cdots \sigma_{i_2} \sigma_{i_1}
\end{displaymath}
be a minimal factorization of $\sigma_{T, T'}$ into adjacent transpositions. We can construct this sequence by writing $\sigma_{T, T'}$ as a word
\begin{displaymath}
(\sigma_{T, T'} (1), \sigma_{T, T'} (2), \ldots , \sigma_{T, T'} (n)),
\end{displaymath}
and transforming $(1, 2, \ldots, n)$ into this word by using adjacent transpositions to move $\sigma_{T, T'} (n)$ to the right into the $n$-th entry, then to move $\sigma_{T, T'} (n-1)$ to the right into the $n-1$-st entry, and so on. This sequence is minimal since each adjacent transposition increases the number of inversions by $1$. (Note that these are inversions in the word representation of the permutation, which need not correspond to inversions in tableaux.)

Then, calculating
\begin{displaymath}
\sigma_{T, T'} \cdot v_{T',T'} = \sigma_{i_k} \sigma_{i_{k-1}} \cdots \sigma_{i_2} \sigma_{i_1} \cdot v_{T',T'},
\end{displaymath}
by Proposition \ref{prop_local_action}, we find that $\sigma_{T, T'} \cdot v_{T',T'}$ is a linear combination of weight vectors $v_{T, T'}$ and $v_{T_i, T'}$, where each $T_i$ is different from $T$ by the minimality of the sequence of transpositions we apply. If we act on this by the orthogonal projector $p_T$, we kill every term except $v_{T, T'}$, so we conclude
\begin{prop}
Up to scalars, we have
\begin{displaymath}
v_{T, T'} = p_T \cdot \sigma_{T, T'} \cdot v_{T',T'} = p_T \cdot \sigma_{T, T'} \cdot p_{T'}.
\end{displaymath}
\end{prop}

The minimal factorization of $\sigma_{T,T'}$ into adjacent transpositions also enables us to generalize (\ref{eq-weight-irred}) beyond $T_{\smin}$ to generate any weight vector from any other by applying $p_T \cdot \sigma_{T,T'}$ to $v_{T'}$ for any $T$ and $T'$. $\sigma_{T,T'}$ acting on $v_{T'}$ smears it into a linear combination of $v_T$ and a set of other weight vectors, and $p_T$ acting on this cleans this up by killing every weight vector except $v_T$. Hence we have
\begin{cor} \label{eq-general-wt}
For any weight vector $v_{T'}$ in any irreducible $S_n$-module, we have, up to scalars, 
\begin{displaymath}
p_T \cdot \sigma_{T, T'} \cdot v_{T'} = v_T .
\end{displaymath}
\end{cor}

\subsection{Highest Weight}
Given a partition $\lambda$ of $n$, let us denote by $\lambda^*$ the conjugate partition to $\lambda$, which is obtained by switching rows and columns, or in other words reflecting $\lambda$ through the NW-SE diagonal (the line $y = -x$). For example, if $\lambda = (3,1)$, then $\lambda^* = (2,1,1)$. (This is also frequently denoted by $\lambda'$ in the literature, but we will stay away from this to avoid confusion with the notation we have been using for pairs of tableaux in the weight basis.)

For a tableau $T$ of shape $\lambda$, we will similarly write $T^*$ for the tableau of shape $\lambda^*$ obtained by reflecting $T$ through the same diagonal. Note that the boxes on the diagonal have content $0$, and reflection through it takes a box of content $c$ in $T$ to a box of content $-c$ in $T^*$. Note also that given $T_{\max}$ of shape $\lambda$, $(T_{\max})^*$ is the minimal tableau of shape $\lambda^*$, and given $T_{\min}$ of shape $\lambda$, $(T_{\min})^*$ is the maximal tableau of shape $\lambda^*$.

Conjugation has a straightforward representation theoretic meaning. Consider the map $\eta: \mathbb{C}[S_n] \to \mathbb{C}[S_n]$ defined by
\begin{displaymath}
\eta: \sigma \mapsto (-1)^{\sigma} \sigma.
\end{displaymath}
$\eta$ is an algebra automorphism of $\mathbb{C}[S_n]$ that takes any transposition to its negative.

If $\rho: S_n \to \End(V)$ is any representation of $S_n$, then $\rho \circ \eta: S_n \to \End(V)$ (twisting with $\eta$) is also a representation, which is the tensor product of the representation we started with and the sign representation $V_{(1^{n})}$. Thus we can view $\eta$ as an operation on representations:
\begin{equation} \label{eq-tensor-sign}
\eta: V \to V \otimes V_{(1^{n})}.
\end{equation}
As an automorphism, $\eta$ preserves irreducibility. For irreducible representations, $\eta$ corresponds exactly to conjugating the corresponding partition together with its weight basis:
\begin{prop} \label{prop-flip}
Let $V_{\lambda}$ be the irreducible representation corresponding to $\lambda$, with weight basis $\{v_T\}$. Then we have
\begin{eqnarray*}
\eta(V_{\lambda}) &=& V_{\lambda} \otimes V_{(1^{n})} \cong V_{\lambda^*}, \\
\eta(v_T) &=& v_{T^*}.
\end{eqnarray*}
\end{prop}

\begin{proof}
Since each $X_i$ is a sum of transpositions, $\eta(X_i) = -X_i$ for all $i$. Hence, the action on the weight basis $\{v_T\}$ is given by
\begin{equation} \label{eq-conj-weight}
\eta(X_i) \cdot v_T = - X_i \cdot v_T = - \cont(b(T, i)) \cdot v_T = \cont(b(T^*, i)) \cdot v_T.
\end{equation}
This shows our assertion on the weight vector level. Since an irreducible representation is determined by its weights, the entire statement follows.
\end{proof}
In particular, if $T = T_{\smin}$ for its shape $\lambda$, then $\eta(v_{T_{\smin}}) = v_{T_{\smax}^*}$, where $T_{\smax}^*$ is maximal for its shape $\lambda^*$. We will refer to $v_{T_{\smax}^*}$ as a {\em highest weight vector}.

Proposition \ref{prop-flip} is true for both the left and right actions of $S_n$, which enables us to compute the action of $\eta$ on the regular representation, thereby giving an explicit realization of (\ref{eq-tensor-sign}). Let $v_{T, T'}$ be a weight vector in $\mathbb{C}[S_n]$. Since $\eta$ is an involution as well as an algebra automorphism, when we apply $\eta$ to (\ref{eq-conj-weight}), substituting $v_{T, T'}$ in for $v_T$, we get
\begin{eqnarray*}
X_i \cdot \eta(v_{T, T'}) &=& \cont(b(T^*, i)) \cdot \eta(v_{T, T'}), \\
\eta(v_{T, T'}) \cdot X_i &=& \cont(b((T')^*, i)) \cdot \eta(v_{T, T'})
\end{eqnarray*}
Since weights with respect to the left and right actions determine weight vectors in $\mathbb{C}[S_n]$, it follows that
\begin{displaymath}
\eta(v_{T, T'}) = v_{T^*, (T')^*},
\end{displaymath}
and since weight vectors determine Gelfand-Tsetlin modules in $\mathbb{C}[S_n]$, we conclude that
\begin{eqnarray*}
\eta(V_{\bullet, T}) &=& V_{\bullet, T^*}, \\
\eta(V_{T, \bullet}) &=& V_{T^*, \bullet}.
\end{eqnarray*}

Since $\eta$ exchanges rows with columns and symmetrization with antisymmetrization, we also have
\begin{eqnarray*}
\eta(r_T) &=& c_{T^*} , \\
\eta(c_T) &=& r_{T^*} .
\end{eqnarray*}
Putting all this together, we can use $\eta$ to transform lowest weight vectors into highest weight vectors, and to obtain corresponding identities. Let $T = T_{\smin}$ for its shape $\lambda$. Applying $\eta$ to equation (\ref{eq-lw-sym}), $v_{T,T} = c_T r_T c_T$, we obtain
\begin{eqnarray*}
v_{T^*, T^*} &=& \eta(v_{T,T}) \\
&=& \eta(c_T r_T c_T) \\
&=& r_T^* c_T^* r_T^*.
\end{eqnarray*}
Here $T^* = T_{\smax}$ for the conjugate partition $\lambda^*$. Since every maximal tableau arises as the conjugate of a minimal tableau, we conclude that
\begin{prop}
Let $T = T_{\smax}$ for its shape $\lambda$. Then we have
\begin{displaymath}
v_{T, T} = r_T c_T r_T.
\end{displaymath}
\end{prop}
Similarly, applying $\eta$ to the identities for the Gelfand-Tsetlin modules corresponding to $T = T_{\smin}$:
\begin{eqnarray*}
V_{\bullet, T} &=& \mathbb{C}[S_n] \cdot r_T c_T , \\
V_{T, \bullet} &=& c_T r_T \cdot \mathbb{C}[S_n] ,
\end{eqnarray*}
we obtain corresponding identities for the Gelfand-Tsetlin modules corresponding to $T_{\smax}$:

\begin{prop}
Let $T = T_{\smax}$ for its shape $\lambda$. Then we also have
\begin{eqnarray*}
V_{\bullet, T} &=& \mathbb{C}[S_n] \cdot c_T r_T , \\
V_{T, \bullet} &=& r_T c_T \cdot \mathbb{C}[S_n] .
\end{eqnarray*}
\end{prop}
We will use this in the next section to interpret the Specht polynomial realization of Young’s natural representation as the image of the highest weight Gelfand-Tsetlin module in $\mathbb{C}[S_n]$.

\section{Isomorphism with the Coinvariant Ring} \label{sec_grading_iso}
\subsection{Functional Realization of $\mathbb{C}[S_n]$}
Our goal for the rest of this paper is to connect the weight space decomposition of $\mathbb{C}[S_n]$ with the geometric realizations of the irreducible representations of $S_n$, which typically arise as various function spaces and their quotients. We begin with a kind of delta function\footnote{More properly, the monomials we work with here are {\em limits} of delta functions of permutations, in a way that we will return to at the end of the paper.} realization of $\mathbb{C}[S_n]$.

Consider the subspace $M^{(n-1)}[z_1, \ldots , z_n]$ of the polynomial ring $R[z_1, \ldots, z_n]$ spanned by monomials in $z_1, z_2, \ldots , z_n$ of {\em multidegree} $(0, 1, \ldots, n-1)$. That is, we consider monomials of the form 
\begin{displaymath}
z_1^{i_1} z_2^{i_2} \cdots z_n^{i_n},
\end{displaymath}
where $i_1, i_2, \ldots , i_n$ is a rearrangement of $0, 1, \ldots , n-1$. Any such monomial can be written as 
\begin{displaymath}
\sigma \cdot z_2 z_3^2 \cdots z_n^{n-1}, 
\end{displaymath}
where $\sigma \in S_n$. We will call $z_2 z_3^2 \cdots z_n^{n-1}$, which corresponds to $\sigma = \id$, the {\em base monomial}.

\begin{prop}
For $R = \mathbb{C}$, consider the map $\Phi: \mathbb{C}[S_n] \to M^{(n-1)} [z_1, \ldots , z_n]$ given by
\begin{displaymath}
\Phi: \sigma \mapsto \sigma \cdot z_2 z_3^2 \cdots z_n^{n-1}.
\end{displaymath}
$\Phi$ is an isomorphism of $S_n$-bimodules, with $S_n$ acting on the left in $M^{(n-1)} [z_1, \ldots , z_n]$ by permuting variables in terms of their indices as usual, and with $S_n$ acting on the right by permuting variables in terms of their degrees. Specifically, we have
\begin{eqnarray}
\sigma \cdot z_1^{i_1} z_2^{i_2} \cdots z_n^{i_n} &=& z_{\sigma(1)}^{i_1} z_{\sigma(2)}^{i_2} \cdots z_{\sigma(n)}^{i_n}, \label{eq-left-act-poly} \\
z_{j_1}^0 z_{j_2}^1 \cdots z_{j_n}^{n-1} \cdot \sigma &=& z_{j_1}^{\sigma^{-1}(1)-1} z_{j_2}^{\sigma^{-1}(2)-1} \cdots z_{j_n}^{\sigma^{-1}(n)-1}, \label{eq-right-act-poly}
\end{eqnarray}
where, in particular, $(j \ j+1)$ acts on monomials in $M^{(n-1)} [z_1, \ldots , z_n]$ on the right by interchanging the variable of degree $j-1$ with the variable of degree $j$.
\end{prop}

\begin{proof}
The left action and the expression (\ref{eq-left-act-poly}) for it are standard. For the right action, we compute the action of an adjacent transposition $(j \ j+1)$ on an arbitrary monomial $z_{j_1}^0 z_{j_2}^1 \cdots z_{j_n}^{n-1}$, as follows. Define $\sigma \in S_n$ by $\sigma(i) = j_i$. Then 
\begin{displaymath}
z_{j_1}^0 z_{j_2}^1 \cdots z_{j_n}^{n-1} = z_{\sigma(1)}^0 z_{\sigma(2)}^1 \cdots z_{\sigma(n)}^{n-1} = \sigma \cdot z_1^0 z_2^1 \cdots z_n^{n-1} = \sigma \cdot \Phi(\id).
\end{displaymath}
We have $\id \cdot (j \ j+1) = (j \ j+1) \cdot \id$ in $\mathbb{C}[S_n]$, i.e. the left action on the identity is equal to the right action. Hence the action induced by $\Phi$ on $M^{(n-1)} [z_1, \ldots , z_n]$ must satisfy 
\begin{displaymath}
\sigma \cdot \Phi(\id) \cdot (j \ j+1) = \sigma \cdot (j \ j+1) \cdot \Phi(\id).
\end{displaymath}
Thus we can compute the right action of $(j \ j+1)$ in $M^{(n-1)} [z_1, \ldots , z_n]$ by
\begin{eqnarray*}
z_{\sigma(1)}^0 z_{\sigma(2)}^1 \cdots z_{\sigma(n)}^{n-1} \cdot (j \ j+1) &=& \sigma \cdot (j \ j+1) \cdot \Phi(\id) \\
&=& \sigma \cdot z_1^0 z_2^1 \cdots z_j^j z_{j+1}^{j-1} \cdots z_n^{n-1} \\
&=& z_{\sigma(1)}^0 z_{\sigma(2)}^1 \cdots z_{\sigma(j)}^j z_{\sigma(j+1)}^{j-1} \cdots z_{\sigma(n)}^{n-1}.
\end{eqnarray*}
This shows that $(j \ j+1)$ acting on the right interchanges the variables of degrees $j-1$ and $j$, as we have asserted. The general equation (\ref{eq-right-act-poly}) follows from this by iterating the action of the adjacent transpositions.
\end{proof}

The fact that the identity element $\id \in S_n$ maps to the base monomial under $\Phi$ makes this functional realization especially convenient for modeling the inductive chain $S_1 \subseteq S_2 \subseteq \cdots S_{n-1} \subseteq S_n$. This is because $\id \in S_{n-1}$ corresponds to $z_2 z_3^2 \cdots z_{n-1}^{n-2} \in M^{(n-2)} [z_1, \ldots , z_{n-1}]$, and so the embedding 
\begin{eqnarray*}
M^{(n-2)} [z_1, \ldots , z_{n-1}] & \to & M^{(n-1)} [z_1, \ldots , z_n] \\
f(z_1, \ldots , z_{n-1}) & \mapsto & f(z_1, \ldots , z_{n-1}) \cdot z_n^{n-1}
\end{eqnarray*}
corresponds to the embedding $\mathbb{C}[S_{n-1}] \subseteq \mathbb{C}[S_n]$ (considered as a representation of $S_{n-1}$ induced to $S_n$). We can think of iterating this map as building up the appropriate base monomial inductively, one variable at a time.

Now, given a partition $\lambda$ of $n$, let $T = T_{\smax}$ be the corresponding maximal tableau, and let $V_{\bullet, T} = \mathbb{C}[S_n] \cdot p_T = \mathbb{C}[S_n] \cdot c_T r_T$ be the corresponding Gelfand-Tsetlin module for the left action (equivalently, the highest weight space for the right action). Under $\Phi$, $c_T r_T = c_T r_T \cdot \id$ maps to
\begin{displaymath}
c_T r_T \cdot \Phi(\id) = c_T r_T \cdot z_2 z_3^2 \cdots z_n^{n-1},
\end{displaymath}
so this polynomial generates $\Phi(V_{\bullet, T})$ (which from now on we will just identify with $V_{\bullet, T}$) as a left $S_n$-module inside $M^{(n-1)} [z_1, \ldots , z_n] \cong \mathbb{C}[S_n]$.

In the case of $T = T_{\smax}$, we have already seen (Remark \ref{rem-specht-proj}) that the corresponding Gelfand-Tsetlin module $V_{\bullet, T} = \mathbb{C}[S_n] \cdot c_T r_T$ projects onto the realization of $V_{\lambda}$ in terms of Specht polynomials, by mapping $c_T r_T$ to the Specht polynomial $f_T(z_1, \ldots , z_n)$ for $T = T_{\smax}$ and extending by the action of $S_n$. Let us write $f_T$ more compactly as
\begin{equation} \label{eq-specht-degree}
f_T(z_1, \ldots , z_n) = c_T r_T \cdot \prod_{i=1}^n z_i^{r_i-1},
\end{equation}
where we write $r_i$ to denote the row in which $i$ appears in the maximal tableau $T_{\smax}$. Thus, if $\lambda = (\lambda_1, \lambda_2, \ldots , \lambda_l)$ as usual, then 
\begin{eqnarray*}
1 &=& r_1 = r_2 = \cdots = r_{\lambda_1}, \\
2 &=& r_{\lambda_1+1} = r_{\lambda_1+2} = \cdots = r_{\lambda_1+\lambda_2},
\end{eqnarray*}
and so on.

After identifying $c_T r_T \in \mathbb{C}[S_n]$ with $c_T r_T \cdot z_2 z_3^2 \cdots z_n^{n-1}$ in $M^{(n-1)} [z_1, \ldots , z_n]$, this projection becomes nothing more than an operation on the degree of each monomial:
\begin{equation} \label{eq-fdr-specht}
\prod_{i=1}^n z_i^{i-1} \mapsto \prod_{i=1}^n z_i^{r_i-1},
\end{equation}
where (\ref{eq-fdr-specht}) applies specifically to the base monomial, and the left action of $S_n$ extends (\ref{eq-fdr-specht}) to all monomials. This just means that, for any monomial in $M^{(n-1)} [z_1, \ldots , z_n]$, our projection reduces the variable of degree $i-1$ down to degree $r_i-1$. Since it commutes with the left action of $S_n$, this projection takes weight vectors in $V_{\bullet, T}$ to weight vectors in the Specht polynomial representation of $V_{\lambda}$.

\subsection{Coinvariant Ring and its Decomposition} \label{sec-charge}
We would like to use the weight basis decomposition to generalize the projection we just defined for the $T = T_{\smax}$ case to every standard tableau $T$. In other words, we want to extend the projection to {\em every} left Gelfand-Tsetlin module $V_{\bullet, T}$; taking their direct sum, this means extending it to the entire regular representation.

To do this, we need a larger space in which we can generalize the Specht polynomial realization of $V_{\lambda}$. This turns out to be the well-known {\em ring of coinvariants}, which we now describe.

Letting $R = \mathbb{C}$ still, consider the full polynomial ring $R[z_1, \ldots , z_n]$ as a left $S_n$-module via the usual action (\ref{eq-left-act-poly}). The subring $\Lambda_n$ of {\em symmetric functions} in $z_1, \ldots , z_n$, which are fixed by this action, is the isotypic component of the trivial representation. As we will be taking quotients, we need to take out the degree $0$ piece ($R$ itself) to keep from killing everything when we pass to the quotient. So, let $\Lambda_n^+$ denote the subring of symmetric functions of positive degree; this is also a (trivial) subrepresentation.

A larger, non-trivial subrepresentation is the {\em ideal} $I_n^+ \subset R[z_1, \ldots , z_n]$ generated by the symmetric functions of positive degree. Both $\Lambda_n^+$ and $I_n^+$ have several well-known sets of generators, including the elementary symmetric functions
\begin{displaymath}
e_i(z_1, \ldots , z_n) = \sum_{1 \leq j_1 < j_2 < \cdots < j_i \leq n} z_{j_1} z_{j_2} \cdots z_{j_i},
\end{displaymath}
and the power sum symmetric functions
\begin{displaymath}
p_i(z_1, \ldots , z_n) = z_1^i + \cdots + z_n^i.
\end{displaymath}
Here the degree $i=1, \ldots , n$, and both the $e_i$ and $p_i$ generate $\Lambda_n^+$ as a subring (without identity) and $I_n^+$ as an ideal in $R[z_1, \ldots , z_n]$.

Finally, the {\em ring of coinvariants} is defined as the quotient ring
\begin{displaymath}
R[z_1, \ldots , z_n]_{S_n} = R[z_1, \ldots , z_n] / I_n^+.
\end{displaymath}
$R[z_1, \ldots , z_n]$ is graded by the usual degree on polynomials, and this grading passes to the ideal $I_n^+$. The generators $e_i$ (or $p_i$) are homogeneous polynomials under this grading, i.e. $I_n^+$ is a {\em homogeneous ideal}. This makes the quotient $R[z_1, \ldots , z_n]_{S_n}$ a graded ring as well, meaning that we can decompose it by degree in a well-defined way into subspaces.

Since $S_n$ preserves $I_n^+$, the action of $S_n$ passes to the quotient $R[z_1, \ldots , z_n]_{S_n}$. Moreover, since $S_n$ preserves the degree $k$ component of $R[z_1, \ldots , z_n]$ for every $k$ (permuting individual variables doesn’t change the total degree of a monomial), and $I_n^+$ is homogeneous, the action of $S_n$ on $R[z_1, \ldots , z_n]_{S_n}$ preserves the degree $k$ component for every $k$ as well.

The degree decompositions of $R[z_1, \ldots , z_n]$, $I_n^+$, and $R[z_1, \ldots , z_n]_{S_n}$ enable us to define {\em graded characters} for each one as a left $S_n$-module. Computing the graded character of $R[z_1, \ldots , z_n]_{S_n}$ gives rise to a formula for the multiplicity of each irreducible representation of $S_n$ in the degree $k$ component of $R[z_1, \ldots , z_n]_{S_n}$. To state it, we will need some combinatorial notation.

Let $T$ be a standard tableau of shape $\lambda$. Define the {\em charge tableau} of $T$, denoted by $\ch(T)$, as the tableau of shape $\lambda$ obtained by replacing each entry $i$ in $T$ with $w_i$, where $w_1 = 0$, and, for $i>1$, $w_i$ is defined inductively as follows:
\begin{equation} \label{eq-charge-def}
w_{i+1} = \left\{ \begin{array}{cll} w_i &\mbox{\rm if}& \cont(b(T,i+1)) > \cont(b(T,i)) \\ w_i + 1 &\mbox{\rm if}& \cont(b(T,i+1)) < \cont(b(T,i)). \end{array} \right.
\end{equation}
To unpack this, recall that $b(T,i)$ and $b(T,i+1)$ are the boxes in $T$ containing $i$ and $i+1$. $\cont(b(T,i+1)) > \cont(b(T,i))$ means that the diagonal containing $b(T,i+1)$ lies above that containing $b(T,i)$, or, equivalently, that $b(T,i+1)$ appears northeast of $b(T,i)$. $\cont(b(T,i+1)) < \cont(b(T,i))$ means the reverse, i.e., $b(T,i+1)$ appears southwest of $b(T,i)$. This really means $i+1$ appears below $i$, since it can’t be immediately west if the tableau is standard. In this case, we call $i$ a {\em descent} of $T$.

To give an example, below we list the standard tableaux of shape $(3,2)$ in lexicographic order, from smallest to largest, with the charge tableau of each one shown directly below it:
\begin{eqnarray*}
\ytableaushort{135,24} 
& 
\ytableaushort{134,25} 
& 
\ytableaushort{125,34} 
\hspace{5mm} 
\ytableaushort{124,35} 
\hspace{5mm} 
\ytableaushort{123,45} 
\\
\ytableaushort{012,12} 
& 
\ytableaushort{011,12} 
& 
\ytableaushort{001,11} 
\hspace{5mm} 
\ytableaushort{001,12} 
\hspace{5mm} 
\ytableaushort{000,11} 
\end{eqnarray*}
Finally, define the {\em charge} of a standard tableau $T$, written $|\ch(T)|$, to be the sum of the entries of its charge tableau: $|\ch(T)| = w_1 + \cdots + w_n$.

There is a parallel notion of {\em cocharge}. Define $\tilde{w}_1 = 0$ and
\begin{displaymath}
\tilde{w}_{i+1} = \left\{ \begin{array}{cll} \tilde{w}_i &\mbox{\rm if}& \cont(b(T,i+1)) < \cont(b(T,i)) \\ \tilde{w}_i + 1 &\mbox{\rm if}& \cont(b(T,i+1)) > \cont(b(T,i)). \end{array} \right.
\end{displaymath}
The {\em cocharge tableau} $\coch(T)$ of $T$ has $\tilde{w}_i$ written in each box in place of $i$. The {\em cocharge} of $T$ is the sum $|\coch(T)| = \tilde{w}_1 + \cdots + \tilde{w}_n$.\footnote{Lascoux and Schutzenberger\cite{bib_LS} originally defined charge and cocharge in a more general way for {\em semistandard tableaux} \cite{bib_LS}. We have interchanged charge and cocharge to make the notation more intuitive for our setting, but apart from this, our definition is consistent with \cite{bib_LS} for the case of standard tableaux.}

If we look at $T$ together with $T^*$, it is clear for any $i = 1, \ldots , n$ that $i$ is a descent of either $T$ or $T^*$, but not both. Thus, by an easy induction on $i$, $w_i$ for $T$ equals $\tilde{w}_i$ for $T^*$ for all $i$, and vice versa. This means conjugation exchanges the charge and cocharge tableaux, and hence charge and cocharge:
\begin{eqnarray*}
\ch(T^*) &=& (\coch(T))^*, \\
|\ch(T^*)| &=& |\coch(T)|.
\end{eqnarray*}
Further, for all $i = 0, \ldots , n-1$, we have
\begin{displaymath}
w_{i+1} + \tilde{w}_{i+1} = w_i + \tilde{w}_i + 1 = i.
\end{displaymath}
Hence, for any standard tableau $T$,
\begin{displaymath}
|\ch(T)| + |\coch(T)| = \sum_{i=1}^n (w_i + \tilde{w}_i) = \sum_{i=1}^n (i-1) = \binom{n}{2}.
\end{displaymath}
Hence minimizing the charge corresponds to maximizing the cocharge, and vice versa. Across all partitions of $n$, charge and cocharge are maximal (and their counterparts minimal) when $\lambda = (1^n)$ and $\lambda = (n)$, respectively.

To minimize the charge for a general shape, we have:
\begin{prop} \label{prop-min-charge}
Let $\lambda = (\lambda_1 , \ldots , \lambda_l)$, where $\lambda_l > 0$. Then:
\begin{enumerate}
\item Any tableau $T$ of shape $\lambda$ has at least $l-1$ descents.
\item $T_{\smax}$ has minimal charge among all tableaux of shape $\lambda$.
\end{enumerate}
\end{prop}

\begin{proof}
By induction on the number of boxes in $T$. The only non-trivial case of the first part is if $T$ has more rows than $T^{(n-1)}$, which means that $\lambda_l = 1$, the single box in the last row contains $n$, and hence adding this box to $T^{(n-1)}$ makes $n-1$ into an additional descent.

For the second statement, the first implies that $w_n \geq l-1$. We can compute
\begin{equation} \label{eq-min-charge}
|\ch(T_{\smax})| = \lambda_2 + 2 \lambda_3 + \cdots + (l-1) \cdot \lambda_l.
\end{equation}
Assume by induction that (\ref{eq-min-charge}) is the minimal possible charge for tableaux with $n-1$ boxes. Since the $n$-th box must contribute at least $l-1$ to the total charge, (\ref{eq-min-charge}) must be the minimal possible charge for tableaux with $n$ boxes as well. 
\end{proof}
The maximal charge for $T$ corresponds to the minimal cocharge, which in turn is the minimal charge for $T^*$. By Proposition \ref{prop-min-charge}, this is given by the maximal tableau of shape $\lambda^*$ in lexicographic order, the conjugate of which is $T_{\smin}$. Thus we have
\begin{cor}
The maximal charge for any shape $\lambda$ is realized at $T = T_{\smin}$, and is given by
\begin{displaymath}
|\ch(T_{\smin})| = \binom{n}{2} - (\lambda_2' + 2 \lambda_3' + \cdots + (l' - 1) \cdot \lambda_{l'}'),
\end{displaymath}
where $(\lambda_1', \ldots , \lambda_{l'}')$ are the column lengths of $\lambda$, i.e., the parts of $\lambda^*$.
\end{cor}

The significance of the charge statistic for the analysis of $R[z_1, \ldots , z_n]_{S_n}$ as an $S_n$-module was first realized via the following theorem \cite{bib_S}:
\begin{thm} \label{thm-cov-decomp}
Denote by $(R[z_1, \ldots , z_n]_{S_n})_d$ the degree $d$ component of the graded left $S_n$-module $R[z_1, \ldots , z_n]_{S_n}$. If $\lambda$ is a partition of $n$, the multiplicity of the irreducible module $V_{\lambda}$ in $(R[z_1, \ldots , z_n]_{S_n})_d$ is equal to the number of standard tableaux $T$ with shape $\lambda$ and charge $d$.
\end{thm}
This is sometimes stated as saying that the {\em graded character} of $R[z_1, \ldots , z_n]_{S_n}$ as an $S_n$-module is given by the Kostka polynomial $K_{\lambda \mu}(q)$, where $\mu = (1^n)$. In this formulation, $K_{\lambda \mu}(q)$ is the generating function for (semistandard) tableaux of shape $\lambda$ and weight (number of $1$’s, $2$’s, etc.) $\mu$, grouped by charge, so setting $\mu = (1^n)$ we get the number of standard tableaux of shape $\lambda$ and charge $d$ as the coefficient of $q^d$.

As a corollary, the multiplicity of $V_{\lambda}$ in $R[z_1, \ldots , z_n]_{S_n}$ is equal to the number of standard tableaux of shape $\lambda$. This is the multiplicity of $V_{\lambda}$ in the regular representation, so {\em $R[z_1, \ldots , z_n]_{S_n}$ is isomorphic to $\mathbb{C}[S_n]$ as a left $S_n$-module}.

The minimal degree copy of $V_{\lambda}$ corresponds to the minimal charge tableau $T_{\smax}$, and has degree $\lambda_2 + 2 \lambda_3 + \cdots + (l-1) \cdot \lambda_l$. By (\ref{eq-specht-degree}), this is exactly the degree of the Specht polynomial $f_T(z_1, \ldots , z_n)$. This enables us to reinterpret the projection (\ref{eq-fdr-specht}) of $V_{\bullet, T}$ into the Specht module as a map into $R[z_1, \ldots , z_n]_{S_n}$, whose image is the lowest degree copy of $V_{\lambda}$, realized in terms of Specht polynomials. The next step is to generalize this to all of $\mathbb{C}[S_n]$ and to the entire ring $R[z_1, \ldots , z_n]_{S_n}$.

\subsection{An Example Beyond Specht Modules} \label{sec-example}
While $R[z_1, \ldots , z_n]_{S_n}$ is isomorphic to $\mathbb{C}[S_n]$ as a left $S_n$-module, it does not come pre-equipped with a right $S_n$-module structure, and thus does not have a decomposition into unique irreducible $S_n$-bimodules. However, we can use the charge statistic to map the weight basis of $\mathbb{C}[S_n]$ to a weight basis inside $R[z_1, \ldots , z_n]_{S_n}$, and thereby to recover much of the structure of $\mathbb{C}[S_n]$.

The key is the second point in the discussion of Proposition \ref{prop_local_action}. In this subsection, we will work out a detailed example ($n=3$), which already captures the main idea. In the two subsections that follow, we will define the charge map in full generality and prove that it is a (graded) $S_n$-module isomorphism from $\mathbb{C}[S_n]$ to $R[z_1, \ldots , z_n]_{S_n}$.

When $n=3$, we have three shapes: $(3)$, $(2,1)$, and $(1,1,1)$. We already understand the first and last, which are the trivial and sign representations. They are $1$-dimensional, and have multiplicity $1$ in $\mathbb{C}[S_n]$, meaning that the Specht module realization of each in $R[z_1, z_2 , z_n]_{S_3}$ is the only one time that isomorphism type appears. For $\lambda = (3)$, the corresponding $S_3$-module projection (\ref{eq-fdr-specht}) $\mathbb{C}[S_3] \to R[z_1, z_2 , z_n]_{S_3}$ takes $z_2 z_3^3$ and every other monomial to $1$ because all $r_i = 1$. For $\lambda = (1,1,1)$, (\ref{eq-fdr-specht}) is the identity, as we have $(r_1, r_2, r_3) = (1,2,3)$. The intuition is that the weight vector corresponding to the sign representation is given by the full Vandermonde determinant, and any reduction of exponents in its expansion would zero it out.

For $\lambda = (2,1)$, we have $T = T_{\smax} = \ytableaushort{12,3}$ and $T' = T_{\smin} = \ytableaushort{13,2}$. The corresponding diagonal weight vectors in the functional realization are
\begin{eqnarray*}
v_{T, T} &=& -\frac{1}{6} (X_2 + 1) (X_3 - 2) \cdot z_2 z_3^2 \\
&=& \frac{1}{6} \left( 2(z_1+z_2) z_3^2 - (z_2+z_3) z_1^2 - (z_1+z_3) z_2^3 \right), \\
v_{T',T'} &=& -\frac{1}{6} (X_2 -1 ) (X_3 + 2) \cdot z_2 z_3^2 \\
&=& \frac{1}{6} \left( 2(z_2-z_1) z_3^2 + (z_2-z_3) z_1^2 + (z_3-z_1) z_2^2 \right).
\end{eqnarray*}
The projection (\ref{eq-fdr-specht}) to the Specht module realization, corresponding to $T = T_{\smax}$, is $z_i z_j^2 \mapsto z_j$, which takes $v_{T,T}$ to the highest weight vector $\frac{1}{3}(2z_3 - z_1 - z_2)$ and {\em kills $v_{T',T'}$}. To see why one survives and the other doesn’t, let us go back to the equation
\begin{equation} \label{eq_sign_charge}
\sigma_i \cdot v_T = \frac{1}{b-a} \left( v_T - v_{\sigma_i T} \right),
\end{equation}
where $a = \cont(b(T, i))$, $b = \cont(b(T, i+1))$, and by convention $v_{\sigma_i T} = 0$ if $\sigma_i \cdot T$ is not standard (that is, if $i$ and $i+1$ are immediately adjacent). The projection (\ref{eq-fdr-specht}) takes $z_2 z_3^2$ and $(1 \ 2) \cdot z_2 z_3^2 = z_1 z_3^2$ to the same element of $R[z_1, z_2 , z_3]_{S_3}$, which is just $z_3$. $v_{T,T}$ has $b-a = 1$ ($i+1=2$ lies northeast of $i=1$), and so by (\ref{eq_sign_charge}), $v_{T,T}$ and $(1 \ 2) \cdot v_{T,T}$ have the same sign (positive) for the coefficient of $z_2 z_3^2$. This in turn means that $v_{T,T}$ has a positive coefficient for both $z_2 z_3^2$ and $(1 \ 2) \cdot z_2 z_3^2 = z_1 z_3^2$, by the $S_3$-invariance of the inner product that picks out these coefficients:
\begin{displaymath}
\langle v_{T,T} \ , \ z_1 z_3^2 \rangle = \langle v_{T,T} \ , \ (1 \ 2) \cdot z_2 z_3^2 \rangle = \langle (1 \ 2) \cdot v_{T,T} \ , \ z_2 z_3^2 \rangle .
\end{displaymath}
Since both monomials that project down to $z_3$ have positive coefficients in $v_{T,T}$, the coefficient of $z_3$ in the image of $v_{T,T}$ must be positive, making the projection non-zero.

On the other hand, $v_{T',T'}$ has $i+1=2$ southwest of $i=1$, so $b-a = -1$. By (\ref{eq_sign_charge}), $(1 \ 2) \cdot v_{T',T'}$ has a negative coefficient for $z_2 z_3^2$, and correspondingly $z_2 z_3^2$ and $z_1 z_3^2$ have {\em opposite} coefficients in $v_{T',T'}$, which cancel out under the projection.

What saves the day for $v_{T',T'}$, and enables it to find a home elsewhere in $R[z_1, z_2 , z_3]_{S_3}$, is that $3$ lies northeast of $2$ in $T'$. Applying (\ref{eq_sign_charge}) and the logic above to the action of $(2 \ 3)$ on $v_{T',T'}$, we find that the coefficients of $z_2 z_3^2$ and $(2 \ 3) \cdot z_2 z_3^2 = z_2^2 z_3$ in $v_{T',T'}$ should have the same (positive) sign, which is indeed what we see in the explicit expression for it above. We can define a different projection $V_{\bullet, T'}  \to R[z_1, z_2, z_3]_{S_3}$ that sends these two monomials to the same term by taking $z_i z_j^2 \mapsto z_i z_j$. This maps $v_{T',T'}$ to 
\begin{displaymath}
\frac{1}{6} \left( 2(z_2-z_1) z_3 + (z_2-z_3) z_1 + (z_3-z_1) z_2 \right) = \frac{1}{2} (z_2 - z_1) z_3,
\end{displaymath} 
which is exactly the lowest weight vector in the degree 2 isotypic component of $V_{(2,1)}$ inside $R[z_1, z_2, z_3]_{S_3}$ predicted by Theorem \ref{thm-cov-decomp}.

The projection $z_i z_j^2 \mapsto z_i z_j$ also identifies pairs of terms in the expansion of $v_{T',T'}$: $z_2 z^3$ and $z_2^2 z_3$ both have positive sign and map to $z_2 z_3$, while $z_1 z_3^2$ and $z_1^2 z_3$ both have negative sign and map to $z_1 z_3$. Finally, $z_1^2 z_2$ and $z_1 z_2^2$ have {\em opposite} signs and cancel out under the projection.

To explain this using group theory, let us look at the elements of $S_3$ that generate each of the above monomials from the base monomial, in pairs corresponding to fibers of our projection. We analyzed the fiber of $z_2 z_3$ above in terms of $(2 \ 3)$. The fiber of $z_1 z_3$ consists of $z_1 z_3^2 = (1 \ 2) z_2 z_3^2$ and $z_1^2 z_3 = (1 \ 2)(2 \ 3) z_2 z_3^2 = (1 \ 2) z_2^2 z_3$. Both $z_2 z_3^2$ and $z_2^2 z_3$ have positive sign by the above, and the action of $(1 \ 2)$ flips both to negative, because $2$ lies below $1$ in $T'$. Finally, consider the (hypothetical) fiber of $z_1 z_2$, consisting of $z_1 z_2^2$ and $z_1^2 z_2$. We have $z_1 z_2^2 = (2 \ 3)(1 \ 2) z_2 z_3^2$ and $z_1^2 z_2 = (1 \ 2)(2 \ 3)(1 \ 2) z_2 z_3^2$.  The action of $(1 \ 2)$ on $(2 \ 3)(1 \ 2) z_2 z_3^2$ flips the sign, enabling these two monomials to cancel each other out when projected to $R[z_1,z_2,z_3]_{S_3}$.

The punch line of our example is that the appropriate projection for each $T$, i.e., for each $V_{\bullet, T}$, is determined by where each entry of $T$ is positioned relative to the preceding one (whether $i+1$ lies northeast or southwest of $i$) — in other words, exactly by the charge tableau $\ch(T)$. The next step is to generalize this for all $n$ and $T$.

\subsection{The Charge Map}
Given a standard tableau $T$, let $(w_1 , \ldots , w_n)$ be the entries of the charge tableau $\ch(T)$, as in Section \ref{sec-charge}. Define the {\em charge monomial} $z^{\minich(T)}$ corresponding to $T$ to be the product $\prod_{i=1}^n z_i^{w_i}$. Thus, for the last example from the previous subsection, we have $T' = \ytableaushort{13,2}$, $\ch(T') = \ytableaushort{01,1}$, $(w_1, w_2, w_3) = (0,1,1)$, and $z^{\minich(T')} = z_2 z_3$. Of course, the charge monomial $z^{\minich(T')} = z_2 z_3$ is exactly where we needed to map base monomial $z_2 z_3^2$ to project $V_{\bullet , T'}$ correctly into $R[z_1, z_2, z_3]_{S_3}$.

Define the {\em charge map} of $T$, written $\Ch_T: M^{(n-1)} [z_1, \ldots , z_n] \to R[z_1, \ldots , z_n] $ by first mapping the base monomial to the charge monomial of $T$, i.e., defining
\begin{equation} \label{eq-fdr}
\Ch_T: \prod_{i=1}^n z_i^{i-1} \mapsto \prod_{i=1}^n z_i^{w_i} = z^{\minich(T)},
\end{equation}
then extending $\Ch_T$ to all monomials in $M^{(n-1)} [z_1, \ldots , z_n]$ using the action of $S_n$. This is a generalization of (\ref{eq-fdr-specht}), because when $T = T_{\smax}$, we have $w_i = r_i - 1$. As in that earlier case, $\Ch_T$ commutes with the left action of $S_n$.
 
We have a kind of duality between $z^{\minich(T)}$ and $z^{\minich(T^*)}$, and hence between $\Ch_T$ and $\Ch_{T^*}$, because
\begin{equation} \label{eq-dual-monom}
z^{\minich(T)} \cdot z^{\minich(T^*)} = \prod_{i=1}^n z_i^{w_i} \cdot \prod_{i=1}^n z_i^{\tilde{w}_i} = \prod_{i=1}^n z_i^{w_i + \tilde{w}_i} = \prod_{i=1}^n z_i^{i-1}.
\end{equation}

Defining $\overline{v}_{T,T} = \Ch_T (v_{T, T}) $, we find 
\begin{eqnarray*}
\overline{v}_{T,T} &=& \Ch_T \left( p_T \cdot \prod_{i=1}^n z_i^{i-1} \right) \\
&=& p_T \cdot \Ch_T \left( \prod_{i=1}^n z_i^{i-1} \right) \\
&=& p_T \cdot \prod_{i=1}^n z_i^{w_i} = p_T \cdot z^{\minich(T)} .
\end{eqnarray*}
In this subsection we will prove 
\begin{prop} \label{prop-nonzero-proj}
$\overline{v}_{T,T} \neq 0$ in $R[z_1, \ldots , z_n]$.
\end{prop}
Subsequently, we will prove that $\overline{v}_{T,T} \neq 0$ in $R[z_1, \ldots , z_n]_{S_n}$, as well.

If we define $S_T \subseteq S_n$ to be the stabilizer of the charge vector $w = (w_1, \ldots , w_n)$ of $T$, or, equivalently, of the charge monomial $z^{\minich(T)}$ (here $S_n$ permutes either the entries of $w$ or the variables in $z^{\minich(T)}$ as usual), then $\Ch_T$ is nothing more than taking the quotient by the action of $S_T$:
\begin{lemma} \label{lemma-stab}
Let $\tau \in S_n$, and let $\tau \cdot z^{\minich(T)}$ be a monomial. Then $\Ch_T^{-1} (\tau \cdot z^{\minich(T)}) = \mathbb{C}[\tau \cdot S_T]$ (as usual, $\tau \cdot S_T$ is the left coset of $S_T$ by $\tau$).
\end{lemma}

\begin{proof}
If $\sigma \in S_T$, then 
\begin{displaymath}
\Ch_T \left (\tau \cdot \sigma \cdot \prod_{i=1}^n z_i^{i-1} \right) = \tau \cdot \sigma \cdot \Ch_T \left( \prod_{i=1}^n z_i^{i-1} \right) = \tau \cdot \sigma \cdot z^{\minich(T)} = \tau \cdot z^{\minich(T)}.
\end{displaymath}
Conversely, given any monomial $\tau' \cdot z^{\minich(T)} \in \Ch_T^{-1}(\tau \cdot z^{\minich(T)})$, we have
\begin{displaymath}
\tau \cdot z^{\minich(T)} = \Ch_T \left( \tau' \cdot \prod_{i=1}^n z_i^{i-1} \right) = \tau' \cdot z^{\minich(T)} ,
\end{displaymath}
and thus $\sigma = \tau^{-1} \cdot \tau' \in S_T$, and $\tau' = \tau \cdot \sigma$ with $\sigma \in S_T$.
\end{proof}
Lemma \ref{lemma-stab} enables us to work with the image of $\Ch_T$ in terms of the functional realization of $\mathbb{C}[S_n]$ in a straightforward way.

If we let $i_1, i_2, \dots , i_k$ be the full set of descents of $T$, written in increasing order, then we have
$w_{i_j+1} = w_{i_j+2} = \cdots = w_{i_{j+1}} = j$ for all $j$, and we can write
\begin{displaymath}
z^{\minich(T)} = (z_{i_1+1} \cdots z_{i_2}) \cdot (z_{i_2+1} \cdots z_{i_3})^2 \cdots (z_{i_k+1} \cdots z_n)^k.
\end{displaymath}
We can then write
\begin{equation} \label{eq-stab-decomp}
S_T = S_{ \{1, \ldots , i_1\} } \times S_{ \{i_1+1, \ldots ,i_2\} } \times \cdots \times S_{ \{i_k+1, \ldots , n \} },
\end{equation}
where $S_{\{ a_1, \ldots , a_j \} } \subseteq S_n$ consists of permutations that only move $\{ a_1, \ldots , a_j \} \subseteq \{1, \ldots , n\}$ and fix everything else. Alternatively, we can view $S_T$ as the subgroup of $S_n$ generated by the adjacent transpositions $\sigma_i$ where $i$ {\em is not} a descent of $T$. Since every $i$ is a descent of either $T$ or $T^*$ but not both, the complementary adjacent transpositions $\sigma_i$, where $i$ {\em is} a descent of $T$, generate $S_{T^*}$.

The following lemma expresses a useful complementarity between $S_T$ and $S_{T^*}$:
\begin{lemma} \label{lemma-group-complement}
Let $\sigma \in S_T$, $\tau \in S_{T^*}$. Write minimal length expansions of $\sigma$ and $\tau$ into adjacent transpositions as in Section \ref{sec-weight-vector}:
\begin{eqnarray*}
\sigma &=& \sigma_{i_k} \sigma_{i_{k-1}} \cdots \sigma_{i_2} \sigma_{i_1}, \\
\tau &=& \sigma_{j_{\kappa}} \sigma_{j_{\kappa-1}} \cdots \sigma_{j_2} \sigma_{j_1}.
\end{eqnarray*}
Then $\sigma_{i_k} \sigma_{i_{k-1}} \cdots \sigma_{i_2} \sigma_{i_1} \cdot \sigma_{j_{\kappa}} \sigma_{j_{\kappa-1}} \cdots \sigma_{j_2} \sigma_{j_1}$ is a minimal length expansion of $\sigma \cdot \tau$.
\end{lemma}

\begin{proof}
Let $m_1, m_2, \ldots , m_d$ be the descents of $T$. Partition $1, \ldots , n$ into subsets as
\begin{displaymath}
\{ 1, \ldots , m_1 \}, \{m_1+1, \ldots , m_2\}, \ldots , \{m_d+1, \ldots , n \},
\end{displaymath}
where some of these sets might have just a single element (if $m_i$ and $m_{i+1}$ are both descents of $T$). The adjacent transpositions that lie in $S_T$ act within each subset, while the adjacent transpositions that lie in $S_{T*}$ interchange the endpoints of adjacent subsets. For example, if $m_1$ is the first descent of $T$, then it is the first non-descent of $T^*$, and thus $\sigma_{m_1}$ is the first adjacent transposition in $S_{T^*}$, interchanging $m_1$ and $m_1+1$.

The key point is that since $\tau$ acts only on the endpoints of adjacent subsets, it preserves the order of each subset. (Given a subset $\{m_l+1, \ldots , m_{l+1}\}$, $S_{T^*}$ can only switch $m_l+1$ with $m_l$ and lower values, can only switch $m_{l+1}$ with $m_{l+1}+1$ and higher values, and cannot touch any elements in between.) As the expansion of $\tau$ into adjacent transpositions is minimal, applying each one in turn increases the number of inversions by $1$. Consequently, when $\sigma$ follows $\tau$, it acts on each subset separately by (\ref{eq-stab-decomp}). Since each subset starts out in order, and since the expansion of $\sigma$ into adjacent transpositions is minimal, applying each one in turn also increases the number of inversions by $1$. Hence the number of inversions increases at each step in the expansion of $\sigma \cdot \tau$, which implies that expansion must be minimal.
\end{proof}

The main idea underlying Proposition \ref{prop-nonzero-proj} is the first part of the following lemma:
\begin{lemma} \label{lemma-sign} Sign Lemma:

Let $\sigma \in S_T$, $\tau \in S_{T^*}$. Then the following hold for the standard $S_n$-invariant inner product on $\mathbb{C}[S_n]$:
\begin{enumerate}
\item $\langle \sigma \cdot v_{T,T} \ , \ \id \rangle$ has the same sign as $\langle v_{T,T} \ , \ \id \rangle$.
\item $\langle \sigma \cdot \tau \cdot v_{T,T} \ , \ \id \rangle$ and $\langle \tau \cdot \sigma \cdot v_{T,T} \ , \ \id \rangle$ have the same sign as $\langle \tau \cdot v_{T,T} \ , \ \id \rangle$.
\end{enumerate}
\end{lemma}
The action of $S_{T^*}$ will come into play when we extend $\Ch_T$ to a projection onto $R[z_1, \ldots , z_n]_{S_n}$ in the next subsection.

\begin{proof}
For the first part, by (\ref{eq-stab-decomp}), we can assume that $\sigma \in S_{ \{ j, j+1, \ldots , j+j' \} }$, where $w_j = w_{j+1} = \cdots = w_{j+j'}$ and $j+m+1$ lies northeast of $j+m$ in $T$ for $m=0,1, \ldots, j'-1$.

As usual, let $\sigma = \sigma_{i_k} \sigma_{i_{k-1}} \cdots \sigma_{i_2} \sigma_{i_1}$ be a minimal factorization into adjacent transpositions, each of which increases the number of inversions (of $j, j+1, \ldots , j+j'$) by $1$. Equivalently, when each $\sigma_{i_m}$ in this factorization acts on $T$, the lower number is southwest of the higher one. By Proposition \ref{prop_local_action}, this means that {\em when $\sigma_{i_m}$ acts on $v_{T,T}$, the coefficient of $v_{T,T}$ is positive. As we will see, this is the essence of how weight vectors come to know their charge.}

Now consider the expansion of $\sigma \cdot v_{T,T} = \sigma_{i_k} \sigma_{i_{k-1}} \cdots \sigma_{i_2} \sigma_{i_1} \cdot v_{T,T}$, with each $\sigma_{i_m}$ term acting in turn using Proposition \ref{prop_local_action}. As in the proof of Proposition \ref{prop-bruhat-exp}, we find that this is a sum of a term $a \cdot v_{T,T}$, where $a > 0$ by the above, along with terms of the form
\begin{displaymath}
a' \cdot v_{\sigma_{j_{\kappa}} \sigma_{j_{\kappa-1}} \cdots \sigma_{j_2} \sigma_{j_1} \cdot T \ , \ T} \ ,
\end{displaymath}
where $j_{\kappa} j_{\kappa-1} \cdots j_2 j_1$ is a non-empty subword of $i_k i_{k-1} \cdots i_2 i_1$ (the empty subword corresponds to the $a \cdot v_{T,T}$ term). As increasing the length increases the number of inversions, all these words must have at least one inversion, and so $\sigma_{j_{\kappa}} \sigma_{j_{\kappa-1}} \cdots \sigma_{j_2} \sigma_{j_1} \cdot T \neq T$ for any of these terms. Thus we have
\begin{eqnarray*}
\langle \sigma \cdot v_{T,T}, \id \rangle &=& \left\langle a \cdot v_{T,T} + \sum a_{T'} v_{T',T} , \id \right\rangle \\
&=& a \cdot \langle v_{T,T}, \id \rangle + \sum a_{T'} \cdot \langle v_{T',T} , \id \rangle
\end{eqnarray*}
for some collection of $T' \neq T$. By Corollary \ref{cor-orthog}, all terms but the first one are zero, and now $a > 0$ yields the first part of the lemma.

For the second part, we revisit the argument for the first part in light of Lemma \ref{lemma-group-complement}. Writing out minimal expansions $\sigma = \sigma_{i_k} \sigma_{i_{k-1}} \cdots \sigma_{i_2} \sigma_{i_1}$ and $\tau = \sigma_{j_{\kappa}} \sigma_{j_{\kappa-1}} \cdots \sigma_{j_2} \sigma_{j_1}$ into adjacent transpositions, we know that the product of the expansions of $\sigma$ and $\tau$ is a minimal expansion for $\sigma \tau$. Switching which tableau we call $T$ and which we call $T^*$, we see that the product of the expansions of $\tau$ and $\sigma$ is a minimal expansion of $\tau \sigma$ as well. Hence, when $\sigma \tau$ and $\tau \sigma$ act on $v_{T,T}$, we still only need to follow the path of $v_{T,T}$. Each adjacent transposition in the expansion of $\tau$ flips the sign of $v_{T,T}$ and each adjacent transposition in the expansion of $\sigma$ preserves it. Letting $\tau$ act first by applying each adjacent transposition in turn, we can either stop at $\sigma_{j_{\kappa}}$ (computing the action of $\tau$ only) or continue all the way to $\sigma_{i_k}$ (computing the action of $\sigma \tau$). Either way, Corollary \ref{cor-orthog} knocks out every term but $v_{T,T}$. Its sign must stay fixed at whatever it was after the action of $\tau$. Similarly, if we let $\sigma$ act first, it keeps the sign of $v_{T,T}$ unchanged, so if we apply the adjacent transposition expansion of $\tau$, the sign of $v_{T,T}$ is the same whether we start from $v_{T,T}$ or from $\sigma \cdot v_{T,T}$, and Corollary \ref{cor-orthog} knocks out every other term again. This proves the second part of the lemma.
\end{proof}

Proposition \ref{prop-nonzero-proj} follows from the Sign Lemma:
\begin{proof}
We show that the coefficient of the reduced base monomial $z^{\minich(T)}$ is non-zero. By Lemma \ref{lemma-stab}, the pre-image of $z^{\minich(T)}$ under $\Ch_T$ is the $\mathbb{c}[S_T]$-orbit of the ordinary base monomial $z_2 z_3^2 \cdots z_n^{n-1}$. Generalizing Section \ref{sec-example}, it is enough to show that in the expansion of $v_{T,T} = p_T \cdot z_2 z_3^2 \cdots z_n^{n-1}$ as a linear combination of monomials, each term in the $S_T$-orbit of $z_2 z_3^2 \cdots z_n^{n-1}$ has the same sign.

The standard $S_n$-invariant inner product on $\mathbb{C}[S_n]$ given by $\langle g, g' \rangle = \delta_{g, g'}$ becomes an $S_n$-invariant inner product on the functional realization $M^{(n-1)} [z_1, \ldots , z_n]$, with the monomials making up an orthonormal basis. Given $\sigma \in S_T$, the coefficient of $\sigma \cdot z_2 z_3^2 \cdots z_n^{n-1}$ in $v_{T,T}$ is
\begin{displaymath}
\langle v_{T,T}, \sigma \cdot z_2 z_3^2 \cdots z_n^{n-1} \rangle = \langle \sigma^{-1} \cdot v_{T,T}, z_2 z_3^2 \cdots z_n^{n-1} \rangle = \langle \sigma^{-1} \cdot v_{T,T}, \id \rangle,
\end{displaymath}
where the last inner product is taken directly in $\mathbb{C}[S_n]$. By (the first part of) the lemma, these all have the same sign for any $\sigma \in S_T$ as for $\sigma = \id$, and we are done.
\end{proof}
As an extension, let $\tau \in S_{T^*}$, and consider the preimage under $\Ch_T$ of $\tau \cdot z^{\minich(T)}$. By Lemma \ref{lemma-stab}, this is just $\mathbb{C}[\tau \cdot S_T] \cdot z_2 z_3^2 \cdots z_n^{n-1}$. The coefficient of $\tau \sigma \cdot z_2 z_3^2 \cdots z_n^{n-1}$ in $v_{T,T}$ is 
\begin{displaymath}
\langle v_{T, T} , \tau \sigma \cdot z_2 z_3^2 \cdots z_n^{n-1} \rangle = \langle \sigma^{-1} \tau^{-1} \cdot v_{T, T} , z_2 z_3^2 \cdots z_n^{n-1} \rangle .
\end{displaymath}
By the second part of Lemma \ref{lemma-sign}, for any $\sigma \in S_T$, $\langle \sigma^{-1} \tau^{-1} \cdot v_{T, T}, z_2 z_3^2 \cdots z_n^{n-1} \rangle$ will have the same sign as it does in the $\sigma = \id$ case. We claim that when $\sigma = \id$, $\langle \tau^{-1} \cdot v_{T, T} , z_2 z_3^2 \cdots z_n^{n-1} \rangle$ has sign $(-1)^{\tau}$. To see this, we write a minimal expansion $\tau^{-1} = \sigma_{i_{\kappa}} \sigma_{i_{\kappa-1}} \cdots \sigma_{i_2} \sigma_{i_1}$ into adjacent transpositions, where each $\sigma_{i_j}$ fixes $z^{T^*}$ and increases the number of inversions. As in the proof of the sign lemma, to compute how this sequence acts on $v_{T,T}$, we only need to track the trajectory of the $v_{T,T}$ term as we apply each adjacent transposition, because of Corollary \ref{cor-orthog} again. Unlike the proof of the sign lemma, in this case $\sigma_{i_j} \in S_{T^*}$ means that $i_j+1$ lies northeast of $i_j$ in $T^*$ and southwest of $i_j$ in $T$, so by Proposition \ref{prop_local_action}, each $\sigma_{i_j}$ {\em flips} the sign of $v_{T,T}$. This means that the aggregate sign applied is $(-1)^{\tau^{-1}} = (-1)^{\tau}$, as asserted.

Projecting every monomial in the fiber of $\tau \cdot z^{\minich(T)}$ down via $\Ch_T$, we conclude
\begin{prop} \label{prop-neg-sign}
Let $\tau \in S_{T^*}$. Then the coefficient of $\tau \cdot v^T$ in $\overline{v}_{T,T} = \Ch_T(v_{T,T})$ is non-zero and has sign $(-1)^{\tau}$ in $R[z_1, \ldots , z_n]$.
\end{prop}
This will be useful in the next subsection.

Replacing $T$ by $T'$ for consistency of notation, consider the restriction of $\Ch_{T'}$ to $V_{\bullet, T'}$. We have just shown it is non-zero, so by Schur’s lemma it must be an isomorphism onto its image  $\Ch_{T'}(V_{\bullet, T'})$. Thus we can define the {\em induced Gelfand-Tsetlin modules} in $R[z_1, \ldots , z_n]$ as
\begin{displaymath}
\overline{V}_{\bullet , T'} = \Ch_{T'}(V_{\bullet, T'}).
\end{displaymath}
The image of every weight vector in $V_{\bullet, T'}$ is a non-zero weight vector in $\overline{V}_{\bullet , T'}$, and we can write the generic weight vector as 
\begin{eqnarray*}
\overline{v}_{T, T'} &=& \Ch_{T'}(v_{T, T'}) \\
&=& p_T \cdot \sigma_{T, T'} \cdot \Ch_{T'}(v_{T', T'}) \\
&=& p_T \cdot \sigma_{T, T'} \cdot \overline{v}_{T', T'} \\
&=& p_T \cdot \sigma_{T, T'} \cdot p_{T'} \cdot z^{T'} .
\end{eqnarray*}

Now that we know that each $\Ch_{T}$ is non-zero on the Gelfand-Tsetlin module $V_{\bullet, T}$, we can glue them together into intertwiner on all of $\mathbb{C}[S_n] = \bigoplus_T V_{\bullet, T}$. Define the full {\em charge map}
\begin{displaymath}
\Ch = \bigoplus_T \Ch_{T} ,
\end{displaymath}
where we restrict each $\Ch_{T}$ to acting on $V_{\bullet, T}$. We know that each $\Ch_{T}$ restricted to $V_{\bullet, T}$ is injective, but not yet that their direct sum $\Ch$ is injective, i.e. that the sum of the $\overline{V}_{\bullet , T'}$ inside $R[z_1, \ldots , z_n]$ is direct. We turn to this next.

The idea is just to project the sequence of induced representations $\Ind_{S_k}^{S_{k+1}} V$ built along the chain $S_1 \subset S_2 \subset S_3 \subset \cdots $ down by $\Ch$, and to use the fact that each step is multiplicity-free, together with Schur’s lemma on each irreducible component, to conclude that $\Ch$ must be an isomorphism onto its (reducible) image. However, in the image of $\Ch$, given a tableau $T$, some steps $S_k \hookrightarrow S_{k+1}$ can be realized as true induced representations and others cannot, so, as we will see, to make the idea work, we will need to go back and forth between conjugate partitions.
\begin{prop} \label{prop-reduced-basis}
The weight vectors $\overline{v}_{T, T'}$ are independent, making $\Ch$ an isomorphism onto its image in $R[z_1, \ldots , z_n]$.
\end{prop}

\begin{proof}
Two reductions of the problem will make our lives significantly easier. First, different polynomials can only cancel each other out in a dependence relation if their degrees match by variable. Thus, it is enough to consider dependence relations among weight vectors $\overline{v}_{T, T'}$ with the same charge vector $w = (w_1, \ldots w_n)$ corresponding to the induction tableau $T'$, and hence the same charge map $\Ch_{T'}$.

Second, if we write $\overline{v}_{T,T'} = p_T \sigma_{T,T'} p_{T'} \cdot z^{\minich(T)}$, we can impose some extra structure via the automorphism $\eta$ of $\mathbb{C}[S_n]$ given by $\eta(\sigma) = (-1)^{\sigma} \sigma$. It is straightforward that $\eta$ takes $p_T$ to the conjugate projector $p_{T^*}$, and hence acts on the sum of induced Gelfand-Tsetlin modules $\overline{V}_{\bullet, T}$ by $\eta(\overline{v}_{T,T'}) = \overline{v}_{T^*,(T')^*}$. Since $\eta$ is an automorphism of $\mathbb{C}[S_n]$, it will preserve dependency relations among weight vectors, which implies that {\em a dependency relation among a set of weight vectors with weight $w = (w_1, \ldots , w_n)$ generates an equivalent dependency relation among the conjugate set of weight vectors}, with the complementary charge vector $(\tilde{w}_1, \ldots , \tilde{w}_n)$. By {\em complementary}, we mean that the sum $(w_1 + \tilde{w}_1, \ldots , w_n+\tilde{w}_n)$ equals $(0, 1, \ldots , n-1)$, which is the relationship between the charge vectors of a tableau and its conjugate tableau. This implies that the sum of the entries of two complementary charge vectors will always equal $\binom{n}{2}$.

Now let $T$ be a standard tableau with $k$ boxes and charge vector $w = (w_1, \ldots , w_k)$, and consider the spaces 
\begin{eqnarray*}
\overline{V}_{T , w_k} &=& \mathbb{C}[S_{k+1}] \cdot \overline{V}_{\bullet, T} \cdot z_{k+1}^{w_k} , \\
\overline{V}_{T , w_k+1} &=& \mathbb{C}[S_{k+1}] \cdot \overline{V}_{\bullet, T} \cdot z_{k+1}^{w_k+1} 
\end{eqnarray*}
inside $R[z_1, \ldots , z_n]$. We can project 
\begin{displaymath}
\Ind_{S_k}^{S_{k+1}} V_{\bullet, T} = \mathbb{C}[S_{k+1}] \cdot V_{\bullet, T} \cdot z_{k+1}^k \subseteq M^{(k)}[z_1, \ldots , z_{k+1}]
\end{displaymath}
onto both spaces if we extend $\Ch_T$ by changing the degree of the new variable $z_{k+1}$ from $k$ to $w_k$ to project onto $\overline{V}_{T , w_k}$, and from $k$ to $w_k+1$ to project onto $\overline{V}_{T , w_k+1}$. Thus, both $\overline{V}_{T , w_k}$ and $\overline{V}_{T , w_k+1}$ are quotients of $\Ind_{S_k}^{S_{k+1}} V_{\bullet, T}$, and hence multiplicity-free by Young’s rule. (Since the projections are isomorphisms on the irreducibles that they do not kill, we can identify them with subspaces of $\Ind_{S_k}^{S_{k+1}} V_{\bullet, T}$, as well.) If $T'$ and $T''$ are $k+1$-box tableaux containing $T$, where $k+1$ is northeast of $k$ in $T'$ and southwest of $k$ in $T''$, then $\Ch$ restricted to $V_{\bullet, T'}$ is equal to the former projection and $\Ch$ restricted to $V_{\bullet, T''}$ is equal to the latter one.

Importantly, when $k+1$ is southwest of $k$, we have
\begin{displaymath}
\mathbb{C}[S_{k+1}] \cdot \overline{V}_{\bullet, T} \cdot z_{k+1}^{w_k+1} \cong \Ind_{S_k}^{S_{k+1}} \overline{V}_{\bullet, T},
\end{displaymath} 
because the spaces $(i \ k+1) \cdot \overline{V}_{\bullet, T} \cdot z_{k+1}^{w_k+1}$ are all independent of each other, as they all have a different variable (just one!) of degree $w_k+1$. This corresponds to the fact that the map that changes the degree of $z_{k+1}$ to $w_k$ extends the one that changes the degree to $w_k+1$, so all irreducible components of $\mathbb{C}[S_{k+1}] \cdot V_{\bullet, T} \cdot z_{k+1}^k = \Ind_{S_k}^{S_{k+1}} V_{\bullet, T}$ survive the projection to $\mathbb{C}[S_{k+1}] \cdot \overline{V}_{\bullet, T} \cdot z_{k+1}^{w_k+1}$, regardless of where $k+1$ is positioned relative to $k$ in the tableau that extends $T$.

This enables us to prove Proposition \ref{prop-reduced-basis} by induction, as follows. Suppose we know all weight vectors in the image of $\Ch$ corresponding to tableaux with $k$ boxes are independent. Let $(w_1, \ldots , w_k, w_{k+1})$ be an admissible degree vector, i.e., $w_1 = 0$, $w_{i+1} = w_i$ or $w_{i+1} = w_i + 1$. Assuming $w_{k+1} = w_k + 1$, look at the spaces $\overline{V}_{\bullet , T}$ corresponding to all tableaux $T$ having charge vector $(w_1, \ldots , w_k)$. They are all independent, and so are their weight vectors, by the inductive hypothesis. Any tableau $T'$ having charge vector $(w_1, \ldots , w_k, w_{k+1})$ has one of these $T$ as the subtableau obtained by removing the box containing $k+1$, and this box lies southwest of the box containing $k$ because we assumed $w_{k+1} = w_k + 1$. Thus we have
\begin{displaymath}
\overline{V}_{\bullet, T'} \subseteq \mathbb{C}[S_{k+1}] \cdot \overline{V}_{\bullet, T} \cdot z_{k+1}^{w_k+1} \cong \Ind_{S_k}^{S_{k+1}} \overline{V}_{\bullet, T},
\end{displaymath}
and the independence of the $\overline{V}_{\bullet, T}$ implies the independence of their images under the action of $\mathbb{C}[S_{k+1}]$ by the above. Looking inside $\mathbb{C}[S_{k+1}] \cdot \overline{V}_{\bullet, T} \cdot z_{k+1}^{w_k+1}$, its decomposition into irreducibles is multiplicity-free, so different $\overline{V}_{\bullet, T'}$ for $T'$ that extend $T$ are independent from each other, and the weight bases inside are independent by construction.

The final trick is to handle the $w_{k+1} = w_k$ case by flipping it. That is, define the complementary degrees $\tilde{w}_i = i-1 - w_i$ and apply the above argument to the degree vector $(\tilde{w}_1, \ldots , \tilde{w}_{k+1})$. This conjugates the corresponding tableaux, and now our earlier discussion implies we are done: the absence of a dependence relation for a given shape and degree vector $\tilde{w}$ implies the absence of such a relation for its conjugate $w$.
\end{proof}

\subsection{Weight Basis in the Coinvariant Ring}
While we have viewed the charge map $\Ch$ as a map into $R[z_1, \ldots , z_n]$ so far, we can compose it with the projection $R[z_1, \ldots , z_n] \to R[z_1, \ldots , z_n]_{S_n}$ and thereby view it as a left $S_n$-intertwiner 
\[
\Ch: \mathbb{C}[S_n]  \to R[z_1, \ldots , z_n]_{S_n}.
\] 
As we know that these two spaces are isomorphic as left $S_n$-modules, we would like to see $\Ch$ turn out to be an isomorphism, and indeed it does:
\begin{thm} \label{thm-iso}
$\Ch: \mathbb{C}[S_n]  \to R[z_1, \ldots , z_n]_{S_n}$
is an isomorphism of left $S_n$-modules, and the weight vectors 
\begin{displaymath}
\overline{v}_{T, T'} = \Ch(v_{T, T'}) = p_T \cdot \sigma_{T, T'} \cdot p_{T'} \cdot z^{T'}
\end{displaymath}
form a basis for $R[z_1, \ldots , z_n]_{S_n}$.
\end{thm}

\begin{proof}
We have done most of the work behind this theorem already. First, let us extend the argument of the preceding subsection to prove that the weight vector $\overline{v}_{T,T}$ in $R[z_1, \ldots , z_n]$ remains non-zero in $R[z_1, \ldots , z_n]_{S_n}$.

The following technical tool, adapted from \cite{bib_ATY}, will help us work modulo $I^+$. Given an admissible degree type $w = (w_1, \ldots , w_n)$, let $\tilde{w}$ be the complementary degree type, let $c_n = \sum (-1)^{\sigma} \sigma$ be the full antisymmetrizer in $S_n$, and define, for any homogeneous polynomial $f(z_1, \ldots , z_n)$ of multi-degree $w$ (this means we allow the degrees of a monomial to be any permutation of $(w_1, \ldots , w_n)$), the operator
\begin{displaymath}
\mathfrak{c}_w : f(z_1, \ldots , z_n) \mapsto c_n \cdot (f(z_1, \ldots , z_n) \cdot z^{\tilde{w}}).
\end{displaymath}
The idea here is that if $f \in I^+$, then $f \cdot z^{\tilde{w}}$, which is homogeneous of total degree $\binom{n}{2}$ by complementarity, will be a linear combination of terms of the form $h(z_1, \ldots , z_n) \tilde{f}(z_1, \ldots , z_n)$, where $h$ is symmetric of positive degree and $\tilde{f}$ consequently has degree less than $\binom{n}{2}$. We have $c_n \cdot (h \cdot \tilde{f} ) = h \cdot (c_n \cdot \tilde{f})$ by symmetry of $h$, and since $\tilde{f}$ has degree less than $\binom{n}{2}$, we conclude $c_n \cdot \tilde{f} = 0$ and so $c_n \cdot f = 0$. Thus, to prove $\overline{v}_{T,T}$ is non-zero when we divide out $I^+$, it suffices to prove it is not killed by $\mathfrak{c}_w$, where $w$ is the charge vector of $T$.

We have $c_n \cdot z_1^{d_1} \ldots z_n^{d_n} = 0$ unless $(d_1, \ldots , d_n)$ is a rearrangement of $(0, \ldots , n-1)$, in which case it is equal to $z_1^{d_1} \ldots z_n^{d_n}$ multiplied by the sign of the element of $S_n$ that takes one to the other. Let us call two monomials, or two (not necessarily increasing) vectors, {\em semi-complementary} if the sum of the degrees of the monomials, or the entries of the vectors, is such a rearrangement. To compute 
\begin{displaymath}
\mathfrak{c}_w(\overline{v}_{T, T}) = \mathfrak{c}_w(\Ch_T(v_{T,T})) = \mathfrak{c}_w(\Ch_T(p_T \cdot z_2 z_3^2 \cdots z_n^{n-1})),
\end{displaymath}
we just need to keep track of which monomials in $\overline{v}_{T,T}$ are semi-complementary to $z^{T^*}$, and which monomials in $v_{T,T}$ project to these under $\Ch_T$.

Starting with
\begin{displaymath}
w + \tilde{w} = (w_1, \ldots , w_n) + (\tilde{w}_1, \ldots , \tilde{w}_n) = (d_1, \ldots , d_n) = (0, 1, \ldots , n-1),
\end{displaymath}
observe that if we act on $z^{T}$, or on $w$, by $\sigma \in S_{T^*} = \Stab(\tilde{w})$, we simply rearrange the sums $d_i$, so the resulting monomial $\sigma \cdot z^{\minich(T)}$ remains semi-complementary to $z^{T^*}$. Conversely, for any monomial that is semi-complementary to $z^{T^*}$, we can always find such a $\sigma$:
\begin{lemma} \label{lemma-complement}
Let $u = (u_1, \ldots , u_n)$ be a charge vector, i.e., $u_1 = 0$, and $u_{i+1}$ equals either $u_i$ or $u_i + 1$. Let $\tilde{u}$ be the vector of complementary degrees, i.e., $\tilde{u}_i = i-1-u_i$. Let $x = (x_1, \ldots, x_n)$ be any vector semi-complementary to $u$. Then we can find $\tau \in \Stab(u) \subseteq S_n$ such that $x = \tau \cdot \tilde{u}$.
\end{lemma}

\begin{proof}
Let $i_1+1$ be the location of the first non-zero entry of $u$ (intuitively, $i_1$ is the first descent of the tableau from which $u$ originates), so that $u_1 = \ldots = u_{i_1} = 0$, $1 = u_{i_1+1} \leq u_{i_1+2} \leq \cdots $. Since $x$ and $u$ are semi-complementary, one of the sums $d_i = x_i + u_i$ is equal to $0$. This forces $x_i = u_i = 0$, so $i \leq i_1$. If we remove this entry from $x$ and $u$, and drop all the entries of $x$ by $1$, the truncated vectors $x'$ and $u'$ will still be semi-complementary. Repeating this, we find an entry of $x'$ that equals $0$ and must be paired with one of the $0$ entries of $u'$, corresponding to an entry of $x$ that equals $1$ and must be paired with one of the $0$ entries of $u$. Continuing, this shows that the full set of entries of $x$ paired with the $i_1$ $0$’s in $u$ must be $0,1, \ldots , i_1-1$. We can put them in order with a permutation that acts on the first $i_1$ entries alone, and thus fixes $u$. If we now remove the first $i_1$ entries for good, drop the entries of $u$ by $1$, and drop the entries of $x$ by $i_1$, then the resulting truncated vectors are still semi-complementary. We can thus repeat this process for each level set of $u$. The resulting permutation $\sigma$ fixes $u$ and puts the entries of $x$ in order across from each level set of $u$, ending up with a vector whose entries increase by $1$ along each level set and stay fixed as we move from one level set to another. This exactly characterizes $\tilde{u}$, so $\sigma \cdot x = \tilde{u}$ and $\tau = \sigma^{-1}$ satisfies the lemma.
\end{proof}
Now let us expand $\overline{v}_{T, T}$ into monomials. The ones that contribute to $\mathfrak{c}_w (\overline{v}_{T, T} )$ are those that are semi-complementary to $z^{T^*}$. By Lemma \ref{lemma-complement}, we can write any such monomial in the form $\tau \cdot z^{\minich(T)}$ for some $\tau \in S_{T^*}$. This puts us exactly in the setting of Proposition \ref{prop-neg-sign}, and we conclude that for fixed $\tau$, the sign of $\tau \cdot z^{\minich(T)}$ in $\overline{v}_{T, T}$ is equal to $(-1)^{\tau}$. This implies that $\tau \cdot z^{\minich(T)}$ contributes {\em positively} to $\mathfrak{c}_w (\overline{v}_{T, T} )$, because 
\begin{eqnarray*}
\mathfrak{c}_w ( (-1)^{\tau} \cdot \tau \cdot z^{\minich(T)} ) &=& c_n ( (-1)^{\tau} \cdot ( \tau \cdot z^{\minich(T)} ) \cdot z^{T^*}) \\
& = & (-1)^{\tau} \cdot c_n ( \tau \cdot z_2 z_3^2 \cdots z_n^{n-1} ) \\
& = & (-1)^{\tau} \cdot (-1)^{\tau} \Delta_n = \Delta_n .
\end{eqnarray*}
Summing up over $\tau \in S_{T^*}$, we conclude that $\mathfrak{c}_w (\overline{v}_{T, T} )$ is a sum of positive terms, hence positive.

Once we know $\overline{v}_{T,T}$ is non-zero modulo $I^+$, we can immediately assert that the entire induced Gelfand-Tsetlin module $\overline{V}_{\bullet, T}$ must be non-zero modulo $I_n^+$ as well, for otherwise the quotient would be an $S_n$-invariant subspace and also a non-trivial quotient of an irreducible $S_n$-module, hence identically $0$. Next, as the weight vectors $\overline{v}_{T, T'}$ for fixed $T'$ span $\overline{V}_{\bullet, T'}$ over $\mathbb{C}$, they certainly span it over $I_n^+$, so they form a basis for $\overline{V}_{\bullet, T'}$ over $I_n^+$ by comparing dimensions, and hence must be independent.

The final step is to prove independence across induced Gelfand-Tsetlin modules, from which the fact that the $\overline{v}_{T,T'}$ make up a basis for $R[z_1, \ldots , z_n]_{S_n}$ follows by a dimension count. It is enough to check that the inductive step in the proof of Proposition \ref{prop-reduced-basis} preserves independence over $I^+$, not just over $\mathbb{C}$. So, assume as part of the inductive hypothesis in that argument, that for $S_k$, the weight vectors $\overline{v}_{T, T'}$ are independent over the ideal $I_k^+$ of positive degree symmetric functions in $k$ variables. The key point is the embedding
\begin{equation} \label{eq-coinv-induct}
\overline{V}_{\bullet, T} \hookrightarrow \mathbb{C}[S_{k+1}] \cdot \overline{V}_{\bullet, T} \cdot z_{k+1}^{w_k+1}.
\end{equation}
Here $w_k$ is fixed, and $T$ is any standard tableau with $k$ boxes whose charge vector $w$ has $w_k$ as its last entry, i.e., $w = (w_1, \ldots , w_k)$. We know this embedding realizes the induced representation $\Ind_{S_k}^{S_{k+1}} \overline{V}_{\bullet, T}$, and hence in particular contains $\overline{V}_{\bullet, T'}$ for every $T'$ obtained from $T$ by adding one box southwest of the one containing $k$. Take the direct sum of (\ref{eq-coinv-induct}) over all $T$ corresponding to a fixed $w_k$. What we need to show is that given $f_1(z_1, \ldots , z_k)$, $f_2(z_1, \ldots z_k)$, etc., where each $f_i$ comes from $\overline{V}_{\bullet, T}$ for a different $T$, (hence the $f_i$ are independent over $I_k^+$ by the inductive hypothesis), the images of the embedding, i.e., $f_1(z_1, \ldots , z_k) \cdot z_{k+1}^{w_k+1}$, $f_2(z_1, \ldots , z_k) \cdot z_{k+1}^{w_k+1}$, etc., are independent over $I_{k+1}^+$.

So, consider a potential non-trivial dependence relation over $I_{k+1}^+$:
\begin{equation} \label{eq-r-dependence}
f_1 z_{k+1}^{w_k + 1} \cdot g_1(e_1, \ldots , e_{k+1}) + \cdots + f_m z_{k+1}^{w_k + 1} \cdot g_m(e_1, \ldots , e_{k+1}) = 0.
\end{equation}
Here $e_i = e_i(z_1, \ldots , z_{k+1})$ is the $i$-th elementary symmetric function, so $g_j(e_1, \ldots , e_{k+1})$ is an arbitrary element of $I_{k+1}^+$. The elementary symmetric functions satisfy the relation
\begin{equation} \label{eq-esym-induct}
e_i(z_1, \ldots , z_{k+1}) = e_i(z_1, \ldots , z_k) + z_{k+1} \cdot e_{i-1}(z_1, \ldots , z_k).
\end{equation}
Now, substitute the right hand side of (\ref{eq-esym-induct}) into (\ref{eq-r-dependence}) and expand (\ref{eq-r-dependence}) out as a polynomial in $z_{k+1}$. Assuming that each $f_i$ is non-trivial, it will have a term that has degree $0$ in $z_{k+1}$ by (\ref{eq-esym-induct}). Thus each term in (\ref{eq-r-dependence}) has ${w_k+1}$ as its lowest degree in $z_{k+1}$, and when we divide (\ref{eq-r-dependence}) by $z_{k+1}^{w_k+1}$ and then set $z_{k+1} = 0$, we obtain a non-trivial dependence relation among the $f_i$ whose coefficients are combinations of the $e_i$ considered now as functions of $z_1, \ldots , z_k$. In other words, this is a dependence relation for the $f_i$ over $I_k^+$, which cannot exist by our inductive hypothesis.

Finally, fixing $T$, since the irreducible components of $\mathbb{C}[S_{k+1}] \cdot \overline{V}_{\bullet, T} \cdot z_{k+1}^{w_k+1}$ are all distinct, and each is the image of a projector $p_{T',T}(X_{k+1})$, we conclude that all the irreducible components $V_{\bullet, T'}$, where $T'$ extends some $T$ corresponding to $w_k$, are independent of each other. This completes the proof of the theorem.
\end{proof}

\subsection{Identity Element in $S_n$ and Coinvariant Ring as Induced Module}
Now that we have an isomorphism $\mathbb{C}[S_n] \to R[z_1, \ldots , z_n]_{S_n}$ and a clear inductive structure on $\mathbb{C}[S_n]$, it is natural to ask about the inductive structure on $R[z_1, \ldots , z_n]_{S_n}$. Our proofs of Proposition \ref{prop-reduced-basis} and Theorem \ref{thm-iso} make implicit use of this already. Beyond this, since the identity element $\id \in S_n$ is the same for every $n$, in the sense that it stays fixed along the inductive chain $S_n \hookrightarrow S_{n+1}$, we can use its image in $R[z_1, \ldots , z_n]_{S_n}$ to anchor the same chain. The weight basis enables us to write this concretely, as we can use our realization (\ref{eq-part-unity}) of the identity in terms of the weight basis in $\mathbb{C}[S_n]$ to project it down to $R[z_1, \ldots , z_n]_{S_n}$:
\begin{displaymath}
\Ch: \id = \sum_{|T|=n} p_T \mapsto \sum_{|T|=n} p_T z^{\minich(T)} .
\end{displaymath}
Then the inductive chain inside $R[z_1, \ldots , z_n]_{S_n}$ becomes
\begin{displaymath}
\mathbb{C}[S_1] \cdot  \sum_{|T|=n} p_T z^{\minich(T)} \subset \mathbb{C}[S_2] \cdot  \sum_{|T|=n} p_T z^{\minich(T)} \subset \mathbb{C}[S_2] \cdot  \sum_{|T|=n} p_T z^{\minich(T)} \subset \cdots .
\end{displaymath}

\section{Further Steps}
We conclude by listing what we see as natural directions to extend this work.
\begin{itemize}
\item {\em RSK correspondence.} The appearance of a natural basis (\ref{left_right_tableaux}) for $\mathbb{C}[S_n]$ indexed by pairs of standard tableaux is a strong indication that the {\em Robinson-Schensted-Knuth correspondence} should be lurking nearby. RSK is a bijection between permutations in $S_n$ and pairs $(T, T')$ of standard tableaux of the same shape with $n$ boxes, inspired by comparing dimensions in the decomposition of the regular representation into isotypic components. In a sequel \cite{bib_ES} to this paper, we realize RSK as a change of basis between two weight bases, one labeled by pairs of tableaux constructed here, and another, labeled by permutations, coming from an external degenerate affine Hecke algebra acting on $\mathbb{C}[S_n]$. The transformation of weight bases corresponds to the quotient map $H_n \to \mathbb{C}[S_n]$, which takes $Y_i \mapsto X_i$ (external translations in $H_n$ to JM-elements in $\mathbb{C}[S_n]$) as mentioned in Section \ref{sec_dAHA}.

\item{\em Geometry.} The following geometric construction of an $S_n$-module isomorphism $\mathbb{C}[S_n] \to R[z_1, \ldots , z_n]_{S_n}$ is parallel to the charge map. Choose distinct points $a_1, \ldots , a_n \in \mathbb{C}$, and identify $S_n$ with the $n!$ possible rearrangements of these points, or, equivalently, with the orbit $\mathcal{O} = S_n \cdot (a_1, \ldots , a_n) \subset \mathbb{C}^n$. Viewing $\mathbb{C}[S_n]$ (more precisely, its dual) as functions on $S_n$, we can identify it with the coordinate ring $R(\mathcal{O}) = R[z_1, \ldots , z_n] / I(\mathcal{O})$, where $I(\mathcal{O})$ is the ideal of polynomials vanishing on $\mathcal{O}$. It is not hard to show that $I(\mathcal{O})$ is generated by polynomials of the form $e_i(z_1, \ldots , z_n) - e_i(a_1, \ldots , a_n)$, where the $e_i$ are the elementary symmetric polynomials.

$I(\mathcal{O})$ is not a homogeneous ideal because of the constant term, however $R(\mathcal{O})$ is still filtered by degree. Constructing the associated graded space $\gr (R(\mathcal{O}))$ amounts to setting the constant terms equal to $0$, i.e., setting all $a_i = 0$, which turns $I(\mathcal{O})$ into $I^+$ (in the notation of this paper), showing that $\gr (R(\mathcal{O})) = R[z_1, \ldots , z_n]_{S_n}$. Thus, the projection $R(\mathcal{O}) \to \gr (R(\mathcal{O}))$ is a geometric version of our charge map.\footnote{The functional realization of $\mathbb{C}[S_n]$ used in this paper sits somewhere in the middle of this picture, as our monomials $z_1^{i_1} z_2^{i_2} \cdots z_n^{i_n}$ can be obtained by taking delta functions in $R(\mathcal{O})$ and letting the $a_i \to 0$.}

We can begin to relate the algebraic and geometric constructions as follows. Our algebraic construction associates a grading of $\mathbb{C}[S_n]$ to an ascending chain of subgroups of $S_n$, while the geometric construction associates a grading of $\mathbb{C}[S_n]$ to a set of points $a_1, \ldots , a_n$. We can put these on a common footing by extending the $S_n$ action to a degenerate affine Hecke algebra $H_n$, in which the external translations $Y_i$ act in terms of the parameters $a_i$ (this is already related to the realization of RSK mentioned above). Then we can try to write the weight basis in $R(\mathcal{O})$ explicitly in terms of this $H_n$ action, and take appropriate limits to obtain the algebraic weight basis constructed here.

In terms of modern algebraic geometry, the limit $a_i \to 0$ corresponds to replacing the orbit $\mathcal{O}$ by a thickened point, which maintains a memory of how the $a_i$ come together. $R[z_1, \ldots , z_n]_{S_n}$ can then be viewed as the coordinate ring of this thickened point in a scheme-theoretic sense. At the same time, $R[z_1, \ldots , z_n]_{S_n}$ is also the cohomology ring of the variety $\mathcal{F}_n$ of full flags in $\mathbb{C}^n$. The flag picture sits above the thickened point picture if we view $(0,\ldots, 0)$ as the set of eigenvalues of the (maximally nilpotent) $n \times n$ $0$-matrix $\mathbf{0}$. $\mathbf{0}$ is a singular point in the variety $\mathcal{N}$ of all nilpotent matrices, and the full flag variety $\mathcal{F}_n$ is the fiber over $\mathbf{0}$ (Springer fiber) in the resolution of singularities (Springer resolution) over $\mathcal{N}$. From the point of view of the thickened point, this is a mechanism for reversing the $a_i \to 0$ limit and recovering (abstractly) the $a_i$. This gives more geometric meaning to the projection $R(\mathcal{O}) \to \gr (R(\mathcal{O})) \cong H^*(\mathcal{F})$, and makes it a natural setting in which to try to interpret the weight basis for $\mathbb{C}[S_n]$, as well as the charge statistic, geometrically, on both the flag and configuration levels.

\item {\em Young modules.} For any partition $\lambda$ of $k$, the representation $\Ind_{S_k \times S_l}^{S_{k+l}} V_{\lambda} \otimes V_{\triv}$, where $V_{\triv}$ is the trivial representation of $S_l$, is also multiplicity-free (Pieri’s rule). It decomposes as the direct sum of $V_{\lambda'}$, where $\lambda'$ is obtained from $\lambda$ by adding $l$ boxes, no more than one per column.

Given $n = \mu_1 + \cdots + \mu_j$, we can then start with the trivial representation  and induce along the chain
\begin{displaymath}
\hspace{1.5cm} 
S_{\mu_1} \times S_{\mu_2} \times \cdots S_{\mu_j} \hookrightarrow S_{\mu_1 + \mu_2} \times S_{\mu_3} \times \cdots S_{\mu_j} \hookrightarrow \cdots \hookrightarrow S_{n-\mu_j} \times S_{\mu_j} \hookrightarrow S_n.
\end{displaymath}
Generalizing Proposition \ref{direct_sum_gz}, this constructs a Young module $M^{\mu}$ and a decomposition of it into irreducible components indexed by {\em semistandard} tableaux $S$ of weight $\mu = (\mu_1, \mu_2, \ldots , \mu_j)$. As in the standard case, $S$ records the induction path leading to a particular irreducible submodule of $M^{\mu}$; specifically, we build up $S$ by writing $i$ in each of the $\mu_i$ boxes added to get to the diagram (isomorphism class of representations) at the $i$-th step of our induction path. The multiplicity of $V_{\lambda}$ in $M^{\mu}$, which is the number of such tableaux of shape $\lambda$, is known as the {\em Kostka number} $K_{\lambda \mu}$. We will refer to the irreducible modules in $M^{\mu}$ resulting from this construction as Gelfand-Tsetlin modules as well.

Like the regular representation, $M^{\mu}$ has a functional realization in terms of polynomials in $n$ variables, which is spanned by monomials in which $\mu_1$ variables have degree $0$, $\mu_2$ have degree $1$, and so on. $M^{\mu}$ is a quotient of $\mathbb{C}[S_n]$, which we can see from the functional realizations: for example, the projection $z_2 z_3^2 z_4^3 \mapsto z_3 z_4$, applied to the entire functional realization of $\mathbb{C}[S_4]$, maps $\mathbb{C}[S_4] \to M^{(2,2)}$. This can also be realized in the geometric picture above by setting some of the $a_i$ equal to each other.

The geometric counterparts of the Young modules are graded quotients $R[z_1, \ldots , R_n] / I_{\mu}$ (here $I_{\mu} \supseteq I^+$) of the space of coinvariants that have the same decomposition into irreducible representations of $S_n$ as the Young modules. These can be interpreted as the graded spaces associated to the coordinate rings of $S_n$-orbits of points $(a_1, \ldots , a_n)$ where some $a_i$ coincide (degenerate orbits), and also as cohomology rings of more general {\em Springer fibers}, which sit over non-zero nilpotent matrices in the resolution of singularities described above. In the geometric picture, the Kostka numbers becomes {\em Kostka polynomials} $K_{\lambda \mu}(q)$, which are the graded characters of these modules, where the (co)charge of a semistandard tableau still corresponds to degree in $R[z_1, \ldots , z_n]$. It would be a natural next step to generalize our construction to an isomorphism $M^{\mu} \to R[z_1, \ldots , R_n] / I_{\mu}$, as well as to connect it to the projection from the coordinate ring of a degenerate orbit to its associated graded space.

\item {\em Schur-Weyl modules and the $n \to \infty$ limit.} Let $V = \mathbb{C}^m$, with the usual left action of $\mathfrak{g} = \mathfrak{sl}_m(\mathbb{C})$. $S_n$ acts on the right on the tensor product $V^{\otimes n}$ by permuting factors, and this commutes with the action of $\mathfrak{g}$ there. The two actions centralize each other, and we have a decomposition
\begin{displaymath}
V^{\otimes n} = \bigoplus_{\lambda} W_{\lambda} \otimes V_{\lambda},
\end{displaymath}
where the $V_{\lambda}$ are irreducible $S_n$-modules as usual, $\lambda$ runs over partitions of $n$ with at most $m$ rows, and the $W_{\lambda}$ are the corresponding irreducible highest weight $\mathfrak{g}$-modules. Each $W_{\lambda}$ decomposes into weight spaces with respect to the Cartan subalgebra $\mathfrak{h} \subset \mathfrak{g}$ consisting of diagonal matrices, and if we consider the weight spaces for $\mathfrak{g}$, they are naturally isomorphic to Young modules. Thus, by the discussion above, $V^{\otimes n}$ has a decomposition into Gelfand-Tsetlin modules for $S_n$, and a simultaneous weight basis for both $\mathfrak{g}$ and $S_n$.

This setup also has a connection to the functional realization. Let $V(z) = V \otimes \mathbb{C}[z, z^{-1}]$, and consider 
\begin{equation} \label{eq-weyl-func}
V(z_1) \otimes \cdots \otimes V(z_n) = V^{\otimes n} \otimes R[z_1^{\pm 1}, \ldots , z_n^{\pm 1}].
\end{equation}
The motivation for adding the $z_i$ in this context is that enables us to extend the action of $\mathfrak{g}$ to an action of the affine Lie algebra $\widehat{\mathfrak{g}}$. If we let $v_1, \ldots , v_m$ be the usual basis for $V$ of eigenvectors of $\mathfrak{h} \subseteq \mathfrak{g}$, we define the action of the extra Serre generators $E_0$ and $F_0$ of $\widehat{\mathfrak{g}}$ on $V(z_i)$ by $E_0 \cdot z_i^k v_j = \delta_{j, 1} \cdot z_i^{k-1} v_n$, and $F_0 \cdot z_i^k v_j = \delta_{j, n} \cdot z_i^{k+1} v_1$. Further, $S_n$ acts on the right simultaneously in both $V^{\otimes n}$ and in $R[z_1, \ldots , z_n]$.

We can try to construct an analog of the charge map that projects onto the quotient of the right hand side of (\ref{eq-weyl-func}) by symmetric functions in the $z_i$, or, more properly, in the $z_i^{-1}$. This is of interest because in the limit $n \to \infty$, the quotient by symmetric functions of negative degree, along with additional wedge relations, becomes an irreducible level $1$ highest-weight $\widehat{\mathfrak{g}}$-module \cite{bib_KMS}. Taking the $n \to \infty$ limit back to $V^{\otimes n}$ by reversing the charge map could give a more sophisticated perspective on that limit in the simple finite dimensional setting as well.

\item {\em Centralizer algebras.} The right action of $S_n$ on $\mathbb{C}[S_n]$ permutes variables in $M^{(n-1)}[z_1, \ldots , z_n]$ in terms of their degrees and centralizes the left action of $S_n$. Passing to $R[z_1, \ldots , z_n]_{S_n}$ via the charge map, we lose this right action, as we no longer have $n$ distinct degrees for our $n$ variables. However, the irreducible decomposition into right $S_n$-modules remains, as the decomposition into left weight spaces $\overline{V}_{T, \bullet} = \Ch(V_{T, \bullet})$. This should help better understand the centralizer of the left action of $S_n$ on $R[z_1, \ldots , z_n]_{S_n}$, and could help in the limit $n \to \infty$ as well.

\item {\em Weight bases beyond $S_n$.} The symmetric group admits several classes of generalizations, including finite and affine Hecke algebras, complex reflection groups $G(m, 1, n) = (\mathbb{Z} / m \mathbb{Z})^n \rtimes S_n$ (which can be thought of as $n$-by-$n$ permutation matrices whose entries are $m$-th roots of unity), and cyclotomic Hecke algebras. Representations can be studied over $\mathbb{C}$ or over fields of finite characteristic. All these generalizations feature inductive chains of nested subalgebras as well as commutative subalgebras analogous to the one generated by the Jucys-Murphy elements, and can thus be decomposed into weight bases. For a sample of what is now an extensive literature, as well as additional references, see \cite{bib_K,bib_Ma}.

One reason for the wide interest in these algebras is their connection to representations of affine Lie algebras and their quantum groups via categorification. In a wide class of cases, when an affine or quantum affine Lie algebra acts on a category of Hecke algebra modules, the upper triangular part acts by restriction (generalizing the idea of removing a box from a Young diagram), and the lower triangular part acts by induction (generalizing the idea of adding a box to a Young diagram). Interpreting the weight bases for various generalizations of $S_n$ in terms of restricted and induced representations should give a direct algebraic realization of these constructions.

\end{itemize}

\end{document}